\documentclass{elsarticle}
\usepackage{graphicx}
\usepackage{subfigure}
\usepackage{color}
\usepackage{tikz}
\usepackage{newlfont}
\usepackage{multirow}
\usepackage{array}
\usepackage{longtable} %Allows for long tables that stretch across many pages
\usepackage{ifthen}
\usepackage{alltt}

\usepackage{enumerate}
 \usepackage[text={6.5in,8.5in},centering]{geometry} 
\newcolumntype{C}[1]{>{\centering\let\newline\\\arraybackslash\hspace{0pt}}m{#1}}
\newcolumntype{L}[1]{>{\raggedright\let\newline\\\arraybackslash\hspace{0pt}}m{#1}}
\newcolumntype{R}[1]{>{\raggedleft\let\newline\\\arraybackslash\hspace{0pt}}m{#1}}
\newcommand{\ba}{\begin{array} }
\newcommand{\ea}{\end{array} }
\newcommand{\bae}{\begin{eqnarray}}
\newcommand{\eae}{\end{eqnarray}}
\newcommand{\bea}{\begin{eqnarray*}}
\newcommand{\eea}{\end{eqnarray*}}
\newcommand{\be}{\begin{equation}}
\newcommand{\ee}{\end{equation}}

\newcommand{\modifyr}[1]{\textcolor{black}{#1}}

\newcommand{\pr}{{\bf Proof}~~}

\usepackage{amssymb}
\usepackage{amsthm}
\usepackage{amsmath}

\usepackage{graphicx}
\usepackage{color}
\usepackage[colorlinks]{hyperref}
\usepackage{lscape}
\usepackage{graphicx}
\usepackage{epstopdf}
\newtheorem{theorem}{\hskip\parindent\bf Theorem}[section]
\newtheorem{lemma}{\hskip\parindent\bf Lemma}[section]
\newtheorem{proposition}{\bf Proposition}[section]

\newtheorem{corollary}{\hskip\parindent\bf Corollary}[section]

\begin{document}
 \markboth{Co-evolutionary dynamics}{Collaborations}
\title{Co-evolutionary dynamics of a host-parasite interaction model: obligate versus  facultative social parasitism}
%a facultative predator and prey interaction model with a Holling Type II functional response}
\author{Yun Kang\footnote{Sciences and Mathematics Faculty, College of Letters and Sciences, Arizona State University, Mesa, AZ 85212, USA ({\tt yun.kang@asu.edu})} and Jennifer Harrison Fewell  \footnote{School of Life Sciences, Arizona State University, Tempe, AZ 85281, USA ({\tt j.fewell@asu.edu})}}%,
\begin{abstract}
Host-parasite co-evolution can have profound impacts on a wide range of ecological and evolutionary processes, including population dynamics, the maintenance of genetic diversity, and the evolution of recombination. To examine the co-evolution of quantitative traits in hosts and parasites, we present and study a co-evolutionary model of a social parasite-host system that incorporates (1) ecological dynamics that feed back into their co-evolutionary outcomes; (2) variation in whether the parasite is obligate or facultative; and (3) Holling Type II functional responses between host and parasite, which are particularly suitable for social parasites that face time costs for host location and its social manipulation. We perform local and global analyses for the co-evolutionary model and the corresponding ecological model. In the absence of evolution, our ecological model analyses imply that an extremely small value of the death rate for facultative social parasites, primarily due to hunting/searching for potential host species, can drive a host extinct globally under certain conditions, while an extremely large value of the death rate can drive the parasite extinct globally. The facultative parasite system can have one, two, or three interior equilibria, while the obligate parasite system can have either one or three interior equilibria. Multiple interior equilibria result in rich dynamics with multiple attractors. The ecological system, in particular, can exhibit bi-stability between the facultative-parasite-only equilibrium and the interior coexistence equilibrium when it has two interior equilibria. Our analysis on the co-evolutionary model provides important insights on how co-evolution can change the ecological and evolutionary outcomes of host-parasite interactions. Our findings suggest that: (a) The host and parasite can select different strategies that may result in local extinction of one species. These strategies can have convergence stability (CS), but may not be evolutionary stable strategies (ESS); (b) The host and its facultative (or obligate) parasite can have ESS that drive the host (or the obligate parasite) extinct locally; (c) Trait functions play an important role in the CS of both boundary and interior equilibria, as well as their ESS; and (d) A small variance in the trait difference that measures parasitism efficiency can destabilize the co-evolutionary system, and generate evolutionary arms-race dynamics with different host-parasite fluctuating patterns. 
\end{abstract}

\bigskip
\maketitle

\section{Introduction}
{Parasitism takes many forms in nature, but can be broadly defined as a symbiosis in which one member (the parasite) benefits from the use of resources gathered by the other member (the host). In most cases, the parasite lives on and/or feeds directly from the host.  In social parasitism, however, the parasite manipulates its host behaviorally rather than physiologically, and derives benefit from work provided by the host (Cini \emph{et al} 2015). The different forms of social parasitism exemplify the diversity of ways in which social parasites can cheat by manipulating the efforts of their hosts. In brood parasitism, a parent deceives or dominates other individuals into rearing their young (Field 1992; Brandt \emph{et al} 2005; Spottiswoode \emph{et al} 2012). For avian brood parasites, this generally occurs via deception, when females lay eggs on other birds's nests and the host parents respond as if the parasitic chick is one of their own (reviewed by Kruger 2007; Davies 2011). This deception is often continued by the parasitic chick, which acts as one of the host brood but may claim a disproportionate proportion of the food resources provisioned by the host parent (Kilner \emph{et al} 2004), or even remove the hosts own chicks from the nest (Spottiswoode \emph{et al} 2012).  \\}

{In another taxonomic realm, insect social parasites, or inquilines, gain acceptance into a social insect colony by mimicking the pheromonal signals of queens or workers; They then live within that colony, exploiting its resources for their own reproduction (Bourke and Franks 1991; Buschinger 2009). Colony parasitism can also occur via more direct aggression, as when queens of some wasp species directly displace other queens and usurp the colony as their own (Field 1992; Shorter and Tibbets 2009). In brood raiding ants, colonies invade nearby nests and retrieve brood, usually in pupal form just before they morph into new adult workers. The emerging workers behave as offspring of the invading colony, performing the same functions as if they were in their natal nest (Brandt \emph{et al} 2005; Holldobler and Wilson 1990 \&2009). In some of the most extreme examples of brood raiding, ant species (slave-making ants) become dependent on the raided workers, without which they cannot functionally maintain their own colonies. Their raiding efficiency and intensity has co-evolutionary consequences for both the host and raiding species (Foitzik \emph{et al} 2001). } \\

{Although diverse in their specifics, each of these cases takes a temporal progression in which the parasite must locate a host and then manipulate it into social acceptance of either the parasite or its offspring; it does so either by deception or aggression. It then establishes an ongoing relationship in which the host provides care, usually to the parasite's offspring, by providing resources, defending the parasitic offspring, and/or providing other parental care. These stages provide a scaffold around which we can model the interplay of ecological and evolutionary effects on parasites and their host. The dynamics within the social parasite-host system can also provide a framework to capture the dynamics inherent in parasite-host relationships more generally.}\\

{Interspecific social parasitism, in which one species parasitizes another, is rarer than intraspecific social parasitism, but can have profound effects on community dynamics. Its impact on the relationship between host and parasite is highly dependent on host number (Sorenson 1997; Brandt {\emph{et al}} 2005). When modeling the dyadic relationship between a parasitic species and a given host, we must consider whether that parasite is fully dependent on that specific host, i.e. obligate, or whether it can utilize other host species as well.  A parasite with multiple potential hosts (generalist parasites), essentially behaves as if it is facultative within the context of that dyadic relationship. For the purposes of this model, we focus on interspecific parasitism, and facultative and generalist parasites are considered as equivalent.}\\

{The ecological and evolutionary dynamics for facultative versus obligate social parasites are very different. A generalist parasite, such as the brown headed cowbird may parasitize tens to hundreds of species, and must flexibly deal with the associated variation in parental care and diets across those species (Rothstein 1975).  In contrast, the common cuckoo specializes on a few species, focusing primarily on the reed warbler. The reed warbler has evolved defensive strategies including recognition of cuckoos and their eggs (Davies and Brook 1988; Stoddard and Stevens 2011).  In return, the cuckoo uses a variety of mimicking strategies, from chick calls to egg color and markings to overcome the reed warbler's defenses (Davies and Brook 1988; Kilner \emph{et al} 1991; Davies 2011; Stoddard and Stevens 2011). The ecological and evolutionary drivers of obligate versus facultative parasitism are complex. Phylogenetic analyses suggest the number of hosts for cowbirds has increasingly expanded evolutionarily (Rothstein \emph{et al} 2002); however, this parasite is also increasing in range and number.  In contrast, the cuckoo host-parasite relationship suggests no clear pattern of expansion over evolutionary time, and variation in host number may be better explained by ecological conditions (Rothstein \emph{et al} 2002).  \\}

{Obligate social parasitism generates a potentially tight dynamical fitness relationship between parasite and host. As the host suffers fitness costs from the parasite, to avoid extinction it must evolve defensive strategies to counter the parasite (Bogusch \emph{et al} 2006). In turn, the parasite, dependent on the resources acquired from its host, is selected to overcome the host's defensive strategies (Poulin \emph{et al}. 2000). As a result, the continuous interactions between the parasite and its host lead to co-evolutionary dynamics, as illustrated by the evolutionary arms-race paradigm that has been used for many host-parasite systems (Anderson and May 1982; Thompson and Burdon 1992; Foitzik \emph{et al} 2003). Close co-evolutionary interactions in a stepwise fashion are especially likely to occur when host and parasite exhibit similar generation times and population sizes; this is generally the case for social parasites, because brood parasitism requires a match with host offspring developmental timing (Foitzik \emph{et al} 2003). In such situations, hosts are expected to more closely match the parasite in strategy evolution, and thus to evolve increased resistance strategies when parasite pressure is strong (Foitzik \emph{et al} 2001\& 2003). In support, recent studies on co-evolution in both avian and insect social parasites found indications of arms races and resistance strategies in highly parasitized host populations (Foitzik {\emph{et al}.} 2001; Kilner and Langmore 2011). \\}

{The evolutionary dynamics of host-parasite interactions can be influenced by multiple intersecting parameters. These may include: the efficiency of parasitism on a specific host; the genetic structure of host and parasite populations; migration rates of parasite and host, and; the degree of mutual specialization (Brandt {\emph{et al}} 2005). The dynamics among these variables can be complex, and their intensity may vary considerably between systems for which the parasite is obligatorily bound to a single host, versus those for which the parasite is facultative on a given host and/or can exploit multiple hosts. In this paper, we develop a simple co-evolutionary model by using evolutionary game theory (EGT), to investigate the ways in which these co-evolutionary dynamics can change the ecological and evolutionary outcomes of host-parasite symbioses.\\}

{Surprisingly, the mathematical models that examine the co-evolution of quantitative traits in hosts and parasites are relatively few in number (Hocherg and van Baalen 1998; van Baalen 1998; Gandon \emph{et al}. 2002; Koella and Bo\"ete 2003; Restif and Koella 2003; Bonds 2006; Zu \emph{et al} 2007; Best \emph{et al} 2009). However, they provide an important theoretical framework within which to consider the importance of co-evolutionary dynamics to the evolution of hosts and parasites. Hochberg and van Baalen (1998) used simulations to study host-parasite co-evolution in relation to spatial heterogeneities. Restif and Koella (2003) showed that when host defense occurs specifically through avoidance, the host-parasite relationship can move towards a single co-evolutionary stable state (CoESS), which cannot be invaded by other strategies. Van Baalen (1998) showed that resistance through recovery, rather than avoidance, can generate bistability such that resistance and virulence are either both low or both high. \\}

{A limitation of the above models is that they examined only the evolutionary stability (ES) of the outcomes, i.e. their ability to be invaded once reached. Eshel (1983) and Christiansen (1991) examined convergence stability (CS), the other important component of the evolutionary process, determining whether populations will move away from a CoESS over evolutionary time and whether evolutionary branching will occur. Dieckmann and Law (1996), and Marrow \emph{et al}. (1996) provided a general understanding of the impact of co-evolution on both ES and CS. In another approach, the work of Kisdi (2006) studied the impact of trade- offs on co-evolutionary dynamics through a general model of predator-prey co-evolution. Zu \emph{et al} (2007) expanded this to examine evolutionary dynamics of the host (as prey) and the parasite (as predator) with a Holling Type II function. Their results showed that branching in the host can induce secondary branching in the parasite, and that these evolutionary dynamics can produce a stable limit cycle. The work of Best {\emph{et al}.} (2009) also emphasized the importance of considering co-evolutionary dynamics and showed that certain highly virulent parasites may result from responses to host evolution.\\}

{These models collectively assume that the parasite (or predator) is obligate. However, facultative and/or generalist social parasites are actually the norm.  Further, many models assume a slow evolutionary timescale to allow application of a timescale separation argument and reduce dimensions and simplify mathematical analysis (Feng \emph{et al} 2004). However, in social parasitism the changes in allele frequencies (and associated phenotypes) likely occur at the same rate as changes in population densities or spatial distributions, which can alter the ecological processes driving changes in population densities or distribution (Gingerich 2009; Cortez and Weitz 2014). Motivated by these issues, we present a fully co-evolutionary model of a social host-parasite system that includes consideration of: (1) the ecological dynamics that feed back into the co-evolutionary outcome; and (2) variation in the form of social parasitism from obligate to facultative. The interaction between host and parasite is modeled by using Holling Type II functional responses. These are particularly suitable for brood parasitism, because they allow consideration of the need for parasites to search for and locate hosts, as well as time interacting with hosts. Our model allow us to explore the following questions: \\}

\begin{enumerate}
\item How do ecological dynamics change when a parasite transitions between being facultative and obligate?
\item Can the parasite or host have evolutionary stable strategies (ESS) that can drive each other extinct at local or global scales?
\item Can co-evolution rescue a host-parasite symbiosis from extinction?
\item What are the effects of trait functions on the ecological and evolutionary dynamics between the parasite and host?
\end{enumerate}

The remaining sections of this article are organized as follows: In Section 2, we provide the background of the EGT modeling approach, define CS and ESS mathematically, and derive our co-evolutionary host-parasite model to incorporate a parasite can be either obligate or facultative. In Section 3, we perform and compare local and global analyses of the ecological dynamics of the co-evolutionary model in the absence of evolution. Our results show that the ecological model of parasite being facultative can generate one, two or three interior equilibria with the possibility that the host goes to extinction either locally or globally. However, the obligate parasite model can have either one or three interior equilibria with the possibility that the obligate parasite goes extinct. In Section 4, we study the dynamics and related ESS of the fully co-evolutionary model of the host-parasite system, including ecological dynamics that feed back into the co-evolutionary outcomes. In particular,  when trait functions follow normal distributions, we derive sufficient conditions that the host-only equilibrium and the parasite-only equilibrium can have ESS, and study the boundary and interior dynamics. The work of the latter case shows that: (a) the co-evolutionary model can have multiple interior equilibria where each equilibrium potentially has two ESS; (b) evolution can save the host from extinction; and (c) a small variance in parasitism efficiency can destabilize the system, thus generating evolutionary arms-race dynamics with different host-parasite fluctuating patterns.
In Section 5, we conclude our study by providing a summary of our results and their potential biological implications. The last section gives the detailed proofs of our analytical results.\\

%%%%%%%%%%%%%%%%%%%%%%%%%%%%%%%%%%%%%%%%%%%%%%%%%%%%%%%%%%%%%%%%

%%%%%%%%%%%%%%%%%%%%%%%%%%%%%%%%%%%%%%%%%%%%%%%%%%%%%%%%%%%%%%%%
\section{Evolutionary game theory models: co-evolutionary host-parasite models}
In an EGT model, the co-evolution of hosts and their parasites can be considered as a mathematical game, where the host and parasite are players, with corresponding strategies, strategy sets, and pay-offs (Vincent and Brown 2005). The strategies represent phenotypic traits with heritable components, and the strategy set is the collection of all evolutionarily feasible strategies for a particular individual, with corresponding fitness pay-offs. Here we use EGT to investigate host-parasite co-evolution within a population dynamical framework, under the assumptions that the interaction between a host and parasite is modeled by a Holling Type II functional response, and that the reliance of the parasite on the identified host can range from completely obligate (the parasite is completely dependent on that host for any fitness pay-off) to facultative (the parasite can survive without parasitism and/or can parasitize other available hosts). \\

%%%%%%%%%%%%%%%%%%%%%%%%%%%%%%%%%%%%%%%%%%%%%%%%%%%%%%%%%%%%%%%%
\subsection{The modeling framework}

We follow the modeling methodology for EGT presented in Vincent and Brown (2005), (also see references Abrams \emph{et al} 1993a, 1993b; Rael \emph{\emph{et al}} 2011; Cushing and Hudson 2012; Kang and Udiani 2014; Kang \emph{et al} 2015). The methodology derives equations that describe the population dynamics of $n$ interacting species $x_i$ together with the dynamics of (mean) phenotypic traits $u_i$ (or strategies) which serve to characterize all individuals of a particular species and are assumed to have a heritable component. 
Let $x = [x_1, x_2]$ denote the vector of population densities of host, parasite and $u = [u_1, u_2]$ denote the vector of all strategies $u=[u_1,u_2] \in [\mathbb U_1,\mathbb U_2]=\mathbb U$ used by the species $i$, which are distinct and drawn from the same set $\mathbb U_i$ of evolutionarily feasible traits. The host-parasite co-evolution dynamic models in terms of differential equations have the following form:
\bae\label{E-pp}
\begin{array}{lcl}
\frac{dx_i}{dt}&=&x_i G_i(v_i,u,x)\vert_{v_i=u_i}=x_iH_i(u,x)\\\\
\frac{du_i}{dt}&=&\sigma_i^2 \frac{\partial G_i(v_i,u,x)}{\partial v_i}\big\vert_{v_i=u_i},\\\\
\end{array}
\eae where $G_i(v_i, u, x)$ {is the fitness of a focal individual that chooses (or inherits) trait $v_i$ when the population has mean
trait $u$ and density $x$}; $H_i(u, x)$ is considered as the fitness function for species $x_i$;
and $\sigma_i^2$ is the variance in traits (strategies) present in species $x_i$ about the mean trait $u_i$ which {therefore a measure of the ``evolutionary speed"}. The host-parasite population dynamics and its associated trait dynamics \eqref{E-pp} together constitute a dynamical system for the 4-vector
$$[x,u] = [x_1, x_2,u_1,u_2]\in \mathbb R^2_+\times \mathbb U_1 \times \mathbb U_2$$
whose dynamic describes an evolutionary process known as \emph{Darwinian dynamics} (Vincent and Brown 2005). These equations \eqref{E-pp}
allow us to investigate the role of evolution plays in determining, for example, the dynamical outcomes of host-parasite interactions. \\

\emph{Darwinian dynamics} will often possess a (locally asymptotically) stable equilibrium $[x^*,u^*]$. Suppose this equilibrium remains stable when additional species (with their associated traits) are added to the community. This means that the coalition of traits $u_i$ associated with those species present in the equilibrium (i.e. for which $x^*> 0$) is able to resist invasion by other species with their other traits. For example, the 4-dimensional equilibrium $[x^*,u^*]$ of \eqref{E-pp}, when embedded in the higher $2m$-dimensional space of the larger community ($m > 2$), remains (locally asymptotically) stable. {In this case, we say that this coalition of traits $u^*$ is an \emph{evolutionary stable strategy} (ESS). The ESS maximum principle provides a necessary condition for a coalition of traits in an equilibrium $[x^*,u^*]$ to be an ESS (Vincent and Brown 2005). According to this principle, if a coalition of strategies is an ESS, then the $G$-function evaluated at equilibrium conditions, i.e. $G(v, u^*, x^*)$, takes on an isolated global maximum with respect to $v$ at each trait $v_i = u_i^*$ in the coalition. Furthermore,this maximum value must equal 0 for the differential equation model.} \\

We say an equilibrium $x^*$ is an \emph{Ecologically Stable Equilibrium} (ESE) if $x^*$ is locally asymptotically stable for the ecological dynamics $\frac{dx_i}{dt}=x_i G_i(v_i,u,x)\vert_{v_i=u_i}=x_iH_i(u,x)$ for given values of traits $u_i$.  We say an equilibrium $(x^*, u^*)$ has \emph{Convergency Stability} if it is locally asymptotically stable in the co-evolutionary model \eqref{E-pp}. We say an equilibrium $(x^*, u^*)$ has $u^*$ as an \emph{Evolutionary Stable Strategy} (ESS) of the co-evolutionary model \eqref{E-pp} if $(x^*, u^*)$ has CS, and  it satisfies the following ESS maximum principle:
\bae\label{ESS-max}\max_{v_i\in \mathbb U_i}\{G_i(v_i, u^*,x^*)\}=G_i(u^*_i,u^*,x^*)= 0 \mbox{ for both } i=1,2 \mbox{ when } x_i^*>0, i=1,2 \mbox { or }\\
\max_{v_i\in\mathbb U_i }\{G_i(v_i, u^*,x^*)\}=G_i(u^*_i,u^*,x^*)\leq 0 \mbox{ for both }  i=1,2 \mbox{ when } x_1^*x_2^*=0
\eae where $u^*=(u^*_1,u^*_2),\,x^*=(x^*_1,x^*_2)$ and $\mathbb U_i$ is the feasible set of trait values $v_i$ which can be a bounded set or $\mathbb R$.\\

%%%%%%%%%%%%%%%%%%%%%%%%%%%%%%%%%%%%%%%%%%%%%%%%%%%%%%%%%%%%%%%%
\subsection{Derivation of a host-parasite co-evolutionary model}
We consider a host-parasite co-evolutionary model with Holling Type II functional responses where the parasite can be facultative or obligate. Let $G_i(v_i, u, x)$ be the fitness of a focal individual that inherits trait $v_i$ when the population has mean trait $u$ and density $x$. We can take $G_i$ in \eqref{E-pp} with the following forms:
\bae\label{E-pp-Gf}
\begin{array}{lcl}
G_1(v_1,u_2,x)&=&r_1(1-\frac{x_1}{K_1(v_1)})-\frac{a(v_1,u_2) x_2}{1+ h\,a(v_1,u_2)x_1}\\\\
G_2(v_2,u_1,x)&=&\frac{e\, a(u_1,v_2) x_1}{1+h\,a(u_1,v_2)x_1}-d(v_2)+r_2 (1-\frac{x_2}{K_2(v_2)})\end{array}
\eae where $a(v_1, v_2)$ denotes the parasitism efficiency of a parasite with phenotypic trait $v_2$ on host individuals with phenotypic trait $v_1$ with the assumption that the stronger host-parasite interactions are, the more similar host and parasite traits are; $K_i(v_i)$ is the carrying capacity of species $i$ individuals with phenotypic trait $v_i$; and $d(v_2)\geq 0$ is the death rate of parasitic individuals with phenotypic trait $v_2$ due to their hunting or attacking all potential hosts. For simplicity, we assume that other parameters (i.e., the host intrinsic growth rate $r_1>0$, the parasite intrinsic growth rate $r_2\geq 0$ in the absence of parasitism towards host species $x_1$, the parasite handling time $h$, and the conversion efficiency $e$) are not influenced by the quantitative traits.\\

As a consequence, the co-evolutionary dynamics of monomorphic resident host and parasite populations with traits $u_1$ and $u_2$ is given by the following set of nonlinear equations:
\bae\label{Ev-pp}
\begin{array}{lcl}
\frac{dx_1}{dt}&=&x_1G_1(v_1,u,x)\vert_{v_1=u_1}=x_1 H_1(u_1,u_2,x_1, x_2)=x_1 \left[r_1 \left(1-\frac{x_1}{K_1(u_1)}\right)-\frac{a(u_1,u_2) x_2}{1+h a(u_1,u_2)x_1}\right]\\\\
\frac{dx_2}{dt}&=&x_2G_2(v_2,u,x)\vert_{v_2=u_2}=x_2 H_2(u_1,u_2,x_1,x_2)=x_2\left[\frac{e a(u_1, u_2) x_1}{1+ha(u_1,u_2)x_1}-d+r_2 \left(1-\frac{x_2}{K_2(u_2)}\right)\right]\\\\
\frac{du_1}{dt}&=&\sigma_1^2\frac{\partial G_1(v_1,u,x)}{\partial v_1}\big\vert_{v_1=u_1}=\sigma_1^2 \frac{\partial H_1(u_1,u_2,x_1,x_2)}{\partial u_1}=\sigma_1^2\left[\frac{r_1 x_1 K_1'(u_1)}{K_1(u_1)^2}-\frac{x_2\frac{\partial a(u_1,u_2)}{\partial u_1}}{(1+h a(u_1,u_2)x_1)^2}\right]\\\\
\frac{du_2}{dt}&=&\sigma_2^2\frac{\partial G_2(v_2,u,x)}{\partial v_2}\big\vert_{v_2=u_2}=\sigma_2^2 \frac{\partial H_2(u_1,u_2,x_1,x_2)}{\partial u_2}=\sigma_2^2\left[\frac{ex_1\frac{\partial a(u_1,u_2)}{\partial u_2}}{(1+h a(u_1,u_2)x_1)^2}-d'(u_2)+\frac{r_2 x_2 K_2'(u_2)}{K_2(u_2)^2}\right]
\end{array}
\eae where trait functions $K_i(u_i),\,i=1,2$ and  $a(u_1,u_2)$ are positive, bounded and smooth functions in $\mathbb U_1\times \mathbb U_2$. The co-evolutionary model \eqref{Ev-pp} incorporates the following ecological assumptions:\\
\begin{enumerate}
\item  \modifyr{In the absence of other species, the population of each species follows a logistic growth function. When a host species $x_1$ is absent, the population dynamics of the parasite $x_2$ depends on the death rate $d$ caused by that parasite on all potential host species, coupled with the parasite's intrinsic growth rate $r_2 >0$, in the absence of parasitism on host species $x_1$. For a chosen trait value $u_2$ of parasite $x_2$, if $d(u_2) > r_2$, then the parasite population with trait $u_2$ will go extinct in the absence of host $x_1$. In such a case, we consider parasite $x_2$ to be obligate. If $d(u_2) < r_2$, then the parasite population with trait $u_2$ can persist in the absence of host $x_1$, and for purposes of this model is called facultative. Thus, in the case of obligate parasitism, the parasite species is unable to persist in the absence of the host.  Under the case of facultative parasitism, the parasite can either (a) function in the absence of any parasitism (fitting the classic definition of facultatively parasitic), or they can survive using other host species.  Assuming parasite-host dynamics are evolutionary labile, this formulation allows us to investigate conditions under which parasites might move evolutionarily from facultative to obligate and/or from specialist to generalist}. \\

\item The interactions between host and parasite can be described using a Holling Type II functional response equation to describe the average attacking rate of a parasite, time spent searching for and handling a host, and time spent on other activities associated with the host-parasite relationship (Anderson and May 1978; Skalsii and Gilliam 2001).The Holling Type II functional response is particularly suitable for modeling social parasitism, because it allows consideration of the main challenges social parasites face: finding a host or host colony, overcoming the defenses of that host, and exploiting its resources for the parasite's own reproduction (Cini \emph{et al} 2011). As an example, in the case of avian brood parasites, overcoming defenses could include the tactics employed by a parasitic adult female and in relation to the host pair, while exploitation would consider the relationship between the chick and the host parents.\\

\item \modifyr{In the absence of co-evolution, i.e., $\sigma_1^2=\sigma_2^2=0$, the co-evolutionary model \eqref{Ev-pp} reduces to the following host-parasite ecological model:
$$\begin{array}{lcl}
\frac{dx_1}{dt}&=&x_1 \left[r_1 \left(1-\frac{x_1}{K_1}\right)-\frac{a x_2}{1+h ax_1}\right]\\\\
\frac{dx_2}{dt}&=&x_2\left[\frac{e a x_1}{1+ha x_1}-d+r_2 \left(1-\frac{x_2}{K_2}\right)\right]\end{array}
 $$ which can be host v.s. facultative or host v.s. obligate parasites depending on the the values of $r_2$ and $d$: parasite $x_2$ is facultative if $r_2<d$ while it is facultative if $r_2>d$. }\\
 \item \modifyr{Social parasitism also occurs facultatively within species.  However, for purposes of this model, we assume that the host and parasite are different species; this allows us to consider them ecologically removed that there is negligible competition for food, space or other resources beyond the social parasitism relationship. }\\
\end{enumerate}

In the following section, we will explore the ecological dynamics of the co-evolutionary model \eqref{Ev-pp} when $\sigma_1^2=\sigma_2^2=0$ (in the absence of evolution).
%%%%%%%%%%%%%%%%%%%%%%%%%%%%%%%%%%%%%%%%%%%%%%%%%%%%%%%%%%%%%%%%
%%%%%%%%%%%%%%%%%%%%%%%%%%%%%%%%%%%%%%%%%%%%%%%%%%%%%%%%%%%%%
\section{Ecological dynamics of a host-parasite model: facultative v.s. obligate parasites}
Let $\sigma_1^2=\sigma_2^2=0$, then the co-evolutionary dynamics of model \eqref{Ev-pp} is reduced to the following ecological model \eqref{pp-g} with fixed traits $u_i, i=1,2$ of host and parasite, respectively:
\bae\label{pp-g}
\begin{array}{lcl}
\frac{dx_1}{dt}&=&x_1G_1(v,u,x)\vert_{v=u_1}=x_1 H_1(u_1,u_2,x_1, x_2)=x_1 \left[r_1 \left(1-\frac{x_1}{K_1(u_1)}\right)-\frac{a(u_1,u_2) x_2}{1+h a(u_1,u_2)x_1}\right]\\\\
\frac{dx_2}{dt}&=&x_2G_2(v,u,x)\vert_{v=u_2}=x_2 H_2(u_1,u_2,x_1,x_2)=x_2\left[\frac{e a(u_1, u_2) x_1}{1+ha(u_1,u_2)x_1}-d(u_2)+r_2 \left(1-\frac{x_2}{K_2(u_2)}\right)\right]\end{array}
\eae whose state space is defined as $X=\{(x_1,x_2)\in \mathbb R^2_+\}$. In addition, we define the following two sets:
$$X_1=\{(x_1,0)\in X\} \mbox{ and } X_2= \{(0,x_2)\in X\}.$$
For convenience, we let $K_1(u_1)=K_1, \,K_2(u_2)=K_2, \,a(u_1,u_2)=a, \mbox{ and } d(u_2)=d$ for short in this section. We have the following lemma regarding the basic dynamical properties of Model  \eqref{pp-g}:\\
%%%%%%%%%%%%%%%%%%%%%%%%%%%%%%%%%%%
\begin{lemma}[Basic dynamical properties]\label{l1:dp} Let $K^M_1=\sup_{u_1\in \mathbb U_1}\{K_1(u_1)\}$ and $K^M_2=\sup_{u_2\in\mathbb U_2}\{K_2(u_2)\}$. Then
the host-parasite model \eqref{pp-g} is positively invariant and every trajectory starting in $X$ attracts to the compact set $C=[0, K_{1}^M]\times[0, \frac{K_{02} (r_2+\frac{e }{h})}{r_2}]$.\\
\end{lemma}

\noindent\textbf{Notes:} The results of Lemma \ref{l1:dp} suggest that our co-evolutionary model \eqref{E-pp-g} is biologically well-defined. 
The host-parasite model \eqref{pp-g} always has the following two boundary equilibria
$E_{00}=(0,0), \, E_{10}=(K_1,0)$, and it has an additional boundary equilibrium $E_{01}=\left(0,(1-\frac{d}{r_2})K_2\right)$ on the $y$-axis if $r_2>d$. The results on the boundary equilibria of Model \eqref{pp-g} are summarized with the following proposition.\\
%%%%%%%%%%%%%%%%%%%%%%%%%%%%%%%%%%%
\begin{proposition}[Boundary equilibria]\label{p1:be}
The host-parasite model \eqref{pp-g} can have two or three boundary equilibria where their existence and stability is listed in Table \ref{t1:be}:
\begin{table}[ht]
\centering
\begin{tabular}{|C{3.5cm}|L{11.5cm}|}
	\hline
Boundary Equilibria &    Stability Condition           \\\hline
$E_{00}$    &  Saddle if $r_2<d$; Source if $r_2>d$.
\\\hline

$E_{10}$      &  Saddle if $r_2+\frac{eaK_1}{1+ahK_1}>d$; Locally asymptotically stable if $r_2+\frac{eaK_1}{1+ahK_1}<d$  \\\hline
 $E_{01}$      &  Requires $r_2>d$ for the existence: it is a saddle if $\frac{r_1}{a }>K_2\left(1-\frac{d}{r_2}\right)>0$ while it is locally asymptotically stable if $\frac{r_1}{a }<K_2\left(1-\frac{d}{r_2}\right)$. \\\hline

         \end{tabular}
\caption{The local stability of boundary equilibria for Model \eqref{pp-g}}
\label{t1:be}
\end{table}
\end{proposition}

\noindent\textbf{Notes:} Proposition \ref{p1:be} implies that for the host-only equilibrium $(K_1,0)$ to be locally stable (i.e. local extinction of parasite) requires that the  death rate of the parasite due to searching for all  potential host species must be larger than the sum of the parasite's intrinsic growth rate $r_2$ in the absence of parasitism towards to $x_1$ and benefits from the parasitism of $x_1$, i.e. $d>r_2+\frac{eaK_1}{1+ahK_1}$. In this case, the parasite is obligate. However, for the parasite-only equilibrium $\left(0,K_2\left(1-\frac{d}{r_2}\right)\right)$ to be locally stable (i.e. local extinction of host) requires that the death rate of parasite due to searching for all potential host species must be smaller than the ratio of sum of the host's intrinsic growth rate $r_1$ to the parasitism efficiency $a$, i.e. $\frac{r_1}{a }>K_2\left(1-\frac{d}{r_2}\right)>0$. This is the case when the parasite is facultative. \\

%%%%%%%%%%%%%%%%%%%%%%%%%%%%%%%%%%%
\begin{theorem}[Persistence of the system]\label{th:persistence}
The host  $x_1$ of Model \eqref{pp-g} is persistent in $X$ if $$d>r_2\mbox{ or  }\,\frac{r_1}{a }>K_2\left(1-\frac{d}{r_2}\right)>0.$$
The parasite $x_2$ of Model \eqref{pp-g} is persistent in $X$ if $$r_2+\frac{eaK_1}{1+ahK_1}>d.$$
Moreover, Model \eqref{pp-g} is permanent in $X$ if either $$ r_2+\frac{eaK_1}{1+ahK_1}>d>r_2$$ or
$$\frac{r_1}{a }>K_2\left(1-\frac{d}{r_2}\right)>0.$$
\end{theorem}

\noindent\textbf{Notes:} Theorem \ref{th:persistence} indicates that host $x_1$ is always persistent if $r_2<d$, while the parasite is always persistent if $r_2>d$. The condition  $r_2>d$ indicates that the parasite is facultative, i.e., it is able to survive in the absence of host $x_1$. The coexistence of both host and parasite requires proper values of $d$ regardless of the parasite being facultative or obligate.\\

%%%%%%%%%%%%%%%%%%%%%%%%%%%%%%%%%%%
 Let $F(x_1)=c_3 x_1^3+c_2 x_1^2+c_1 x_1+c_0$ where
$$\begin{array}{lcl}
c_3=r_1r_2 (a h)^2&&c_2=ahr_1r_2(2-ahK_1)\\
c_1=a^2hK_1K_2\left(r_2+\frac{e}{h}-d\right)+r_1r_2(1-2ahK_1)&&c_0=a r_2K_1\left[K_2\left(1-\frac{d}{r_2}\right)-\frac{r_1}{a}\right]\end{array}$$ and 
$$x_c^1=\frac{-c_2-\sqrt{c_2^2-3c_1c_3}}{3c_3}<x_c^2=\frac{-c_2+\sqrt{c_2^2-3c_1c_3}}{3c_3}.$$
%%%%%%%%%%%%%%%%%%%%%%%%%%%%%%%%%%%
\begin{theorem}[Interior equilibria]\label{th:interior}
Model \eqref{pp-g} can have none, one, two or three interior equilibria depending on the values of parameters. The sufficient conditions for the existence of interior equilibria are listed in Table \ref{t2:interior}. \begin{table}[ht]
\centering
\begin{tabular}{|C{4cm}|L{10cm}|}%{|P{3.5cm}|P{11.5cm}|}
	\hline
Number of Interior Equilibria &    Sufficient Conditions          \\\hline
None    &  One of the following conditions holds: (1) $d>r_2+\frac{aeK_1}{ahK_1+1}$; or (2) $\frac{r_1}{a }<K_2\left(1-\frac{d}{r_2}\right)$ and $a^2 r_2 K_1K_2\left(h+\frac{d}{r_2}-\frac{e}{r_2}\right)>\frac{r_1r_2(ahK_1+1)^2}{3}$; or (3) $ K_1<\frac{1}{ah}$ and $\frac{r_1}{a }<K_2\left(1-\frac{d}{r_2}\right)$; or (4) $ K_1>\frac{1}{ah}$ and $\frac{r_1(1+ahK_1)^2}{4a^2K_1K_2}<1-\frac{d}{r_2}.$
\\\hline

One      &  One of the following conditions holds: (1)$r_2<d<r_2+\frac{aeK_1}{ahK_1+1}$; or (2) $\frac{r_1}{a }>K_2\left(1-\frac{d}{r_2}\right)$ and $\left(h+\frac{d}{r_2}-\frac{e}{r_2}\right)>\frac{r_1(ahK_1+1)^2}{3a^2 K_1K_2}$;  or (3) $\frac{r_1}{a}>K_2\left(1-\frac{d}{r_2}\right)$,\,$\left(h+\frac{d}{r_2}-\frac{e}{r_2}\right)<\frac{r_1(ahK_1+1)^2}{3a^2  K_1K_2}$ and $F(x^2_c)<0$.
\\\hline
Two     &   $\frac{r_1}{a }<K_2\left(1-\frac{d}{r_2}\right)$,\,$\left(h+\frac{d}{r_2}-\frac{e}{r_2}\right)<\frac{r_1(ahK_1+1)^2}{3a^2  K_1K_2}$ and $F(x^2_c)<0$ with $x^2_c>0$.
\\\hline
Three     &  $\frac{r_1}{a }>K_2\left(1-\frac{d}{r_2}\right)$,\,$\left(h+\frac{d}{r_2}-\frac{e}{r_2}\right)<\frac{r_1(ahK_1+1)^2}{3a^2  K_1K_2}$, $F(x^1_c)>0$,\, and $F(x^2_c)<0$ with $x^1_c>0$. \\\hline

         \end{tabular}
\caption{The existence condition of interior equilibria for Model \eqref{pp-g}.}
\label{t2:interior}
\end{table}Let $(x_1^*,x_2^*)$ be an interior equilibrium of Model \eqref{pp-g}, then it is locally asymptotically stable if the following conditions are satisfied
\bae\label{stable}\frac{r_2}{K_2}>\frac{a}{2} \mbox{ and } \frac{r_1r_2}{K_1K_2}+\frac{ea^2}{(1+ah K_1)^3}>a^2hK_{2} (r_2+\frac{e }{h}).\eae\end{theorem}

\noindent\textbf{Notes:} Theorem \ref{th:interior} indicates that the ecological model  \eqref{pp-g} can have one, two, or three interior equilibria when parasite is facultative, i.e., $r_2>d$; while  \eqref{pp-g} can have one,  or three interior equilibria when parasite is obligate, i.e., $r_2<d$. Theoretically, it is possible that  \eqref{pp-g} has two interior equilibria in the singular cases when $r_2<d$. We did not investigate such singular cases. The result for ignoring these singular cases is that host persists in Model \eqref{pp-g} whenever $r_2<d$. And this is different from the case when $r_2>d$. The detailed results on the persistence have been provided in Theorem \ref{th:persistence}. In addition, we can conclude that  the inequalities \eqref{stable} hold if either $\frac{r_2}{K_2}$ is large enough or $a$ is small enough; however, the value of $d$ is not presented in \eqref{stable}.  In the case that $r_2>d$, then small values of $K_2$ can lead to a unique locally stable interior equilibrium. The existences of two or three interior equilibria suggest that under certain conditions, Model  \eqref{pp-g} can have multiple attractors and generate complicated dynamics.\\

%%%%%%%%%%%%%%%%%%%%%%%%%%%%%%%%%%%

\begin{theorem}[Global stability and extinction of one species]\label{th:global}
Model \eqref{pp-g} has global stability at $E_{10}=(K_1, 0)$ whenever the inequality $d>r_2+\frac{eaK_1}{1+ahK_1}$ holds. 
Model \eqref{pp-g} has global stability at $E_{01}=\left(0,(1-\frac{d}{r_2})K_2\right)$ if $\frac{r_1}{a }<K_2\left(1-\frac{d}{r_2}\right)$ and one of the following conditions is satisfied
\begin{enumerate}
\item  $a^2 r_2 K_1K_2\left(h+\frac{d}{r_2}-\frac{e}{r_2}\right)>\frac{r_1r_2(ahK_1+1)^2}{3}$; or 
\item $ K_1<\frac{1}{ah}$; or 
\item $ K_1>\frac{1}{ah}$ and $\frac{r_1(1+ahK_1)^2}{4a^2hK_1K_2}<1-\frac{d}{r_2}.$\\
\end{enumerate}
\end{theorem}

\noindent\textbf{Notes:} Theorem \ref{th:global} implies that the necessary condition for the global stability of $E_{10}$ is $r_2<d$. This indicates that the parasite $x_2$ cannot be facultative if host species $x_1$ persists. If the parasite is facultative (i.e., $r_2>d$), then the parasite $x_2$ can drive the host $x_1$ extinct locally when $K_2$ is large enough and $K_1$ is small enough such that the inequalities $K_1r_1<\frac{r_1}{a}<K_2\left(1-\frac{d}{r_2}\right)$ hold; local extinction of the host can also occur when $K_2$ is large enough and $K_1ah>1$ such that the inequality $K_2\left(1-\frac{d}{r_2}\right)> \frac{r_1(1+ahK_1)^2}{4a^2K_1}$ holds.\\\\
%%%%%%%%%%%%%%%%%%%%%%%%%%%%%%%%%%%%%%%%%%%%%%%%%%%%%%%%%%%%%

\noindent\textbf{Ecological dynamics of facultative versus obligate parasites:} The condition of the parasite $x_2$ being either facultative or obligate on host $x_1$ has strong effects on the ecological dynamics of the host-parasite interaction model  \eqref{pp-g}. These impacts have been listed in Table \ref{t3:comparison}. We can summarize the main differences as follows: 
\begin{enumerate}
\item Within the constraints of this model, a facultative parasite is always persistent, because it is not fully dependent on the dynamics with its host. If the parasite is obligate, the host always persists. A facultative parasite can cause the extinction of the host, while an obligate parasite can go extinction itself under certain conditions. 
\item When the parasite is facultative, either large values of $\frac{r_1}{aK_2}$ or small values of $r_2-d$ can result in its permanence. If the parasite is obligate, both the host and the parasite can persist under intermediate values of the parasite death rate $d$.  
\item When the parasite is facultative, the system is prone to coexistence of both species and can produce one, two, or three interior equilibria. When the parasite is obligate, the system can have one or three interior equilibria. If the system has two interior equilibria, it exhibits
bistability between the parasite-only boundary attractor and the  interior coexistence attractor. If the system has three interior equilibria, it may have two interior coexistence attractors.

\begin{table}[ht]
\centering
\begin{tabular}{|C{4cm}|L{5cm}|L{6cm}|}
	\hline
Cases &    facultative parasite ($r_2>d$)& obligate parasite ($r_2\leq d$)         \\\hline
Stability of $E_{00}$&Source& Saddle\\\hline
Stability of $E_{10}$&Saddle& Stable if $r_2+\frac{eaK_1}{1+ahK_1}<d$\\\hline
Stability of $E_{01}$&A saddle if $\frac{r_1}{a }>K_2\left(1-\frac{d}{r_2}\right)>0$; A sink if $\frac{r_1}{a K_2}<1-\frac{d}{r_2}$& Does not exist \\\hline
Persistence of host $x_1$&$\frac{r_1}{a}>K_2\left(1-\frac{d}{r_2}\right)>0$& Always \\\hline
Persistence of parasite $x_2$&Always& $r_2+\frac{eaK_1}{1+ahK_1}>d>r_2$\\\hline
Permanence of the system&$\frac{r_1}{a}>K_2\left(1-\frac{d}{r_2}\right)>0$&$r_2+\frac{eaK_1}{1+haK_1}>d>r_2.$ \\\hline
Extinction of host $x_1$&$\frac{r_1}{a}>K_2\left(1-\frac{d}{r_2}\right)$ and one of the following conditions holds: (1)$ K_1<\frac{1}{ah}$ or (2)$ K_1>\frac{1}{ah}$ and $K_2\left(1-\frac{d}{r_2}\right)> \frac{r_1(1+ahK_1)^2}{4a^2hK_1}$ or (3)$\left(h+\frac{d}{r_2}-\frac{e}{r_2}\right)>\frac{r_1(ahK_1+1)^2}{3a^2  K_1K_2}$& Never\\\hline
Extinction of parasite $x_2$&Never & $d>r_2+\frac{eaK_1}{1+ahK_1}$.\\\hline
Number of interior equilibria&None, one, two, or three; This suggests the possibility of two types of multiple attractors: (a) two interior attractors; and (b) a boundary attractor and an interior attractor.& None, one, or three; This suggests the possibility of two interior attractors. \\\hline

         \end{tabular}
\caption{The comparison of sufficient conditions that lead to different dynamical outcomes of a parasite being facultative versus obligate for  Model \eqref{pp-g}.}
\label{t3:comparison}
\end{table}
\end{enumerate}

\noindent\textbf{Notes:} In the case that the parasite is facultative, i.e., $d<r_2$, its population model can be rewritten as
$$\frac{dx_2}{dt}=x_2\left[\frac{e a x_1}{1+ha x_1}+(r_2-d) (1-\frac{x_2}{K_2(1-\frac{d}{r_2})})\right]$$ which can be rewritten as
$$\frac{dx_2}{dt}=x_2\left[\frac{e a x_1}{1+ha x_1}+r_2 (1-\frac{x_2}{K_2})\right]$$ by letting
$$r_2\rightarrow r_2-d \mbox{ and } K_2\rightarrow K_2(1-\frac{d}{r_2}).$$
This is equivalent to the case when $d=0$ in the host-parasite model. Thus, we can obtain the simplified host-parasite model \eqref{pp-d0} with Holling-Type II functional responses when the parasite is facultative 
\bae\label{pp-d0}
\begin{array}{lcl}
\frac{dx_1}{dt}&=&x_1 \left[r_1 (1-\frac{x_1}{K_1})-\frac{a x_2}{1+h a x_1}\right]\\\\
\frac{dx_2}{dt}&=&x_2\left[\frac{e a x_1}{1+hax_1}+r_2 (1-\frac{x_2}{K_2})\right].\end{array}
\eae All the analytic results presented in this section can be applied to Model \eqref{pp-d0} by letting $d=0$. \\
%%%%%%%%%%%%%%%%%%%%%%%%%%%%%%%%%%%%%%%%%%%%%%%%%%%%%%%%%%%%%

\noindent\textbf{Numerical investigations and bifurcation diagrams:}  According to Theorem \ref{th:interior}, the ecological model \eqref{pp-g} can have none, one, two or three interior equilibria depending on parameter values. Unfortunately, the explicit forms of these interior equilibria are far too complicated. Our analytical results suggest that the values of $r_2$ and $d$ have profound effects on the persistence of the system (see Theorem \ref{th:persistence}-\ref{th:global}) and the existence of the number of  interior equilibria as well as the stability of these equilibria (see Theorem \ref{th:interior}). For these reasons, we use bifurcation diagrams to explore the effects of $r_2$ and $d$ on the system's dynamical patterns.\\

\begin{figure}[ht]
\centering
   \includegraphics[scale =.45] {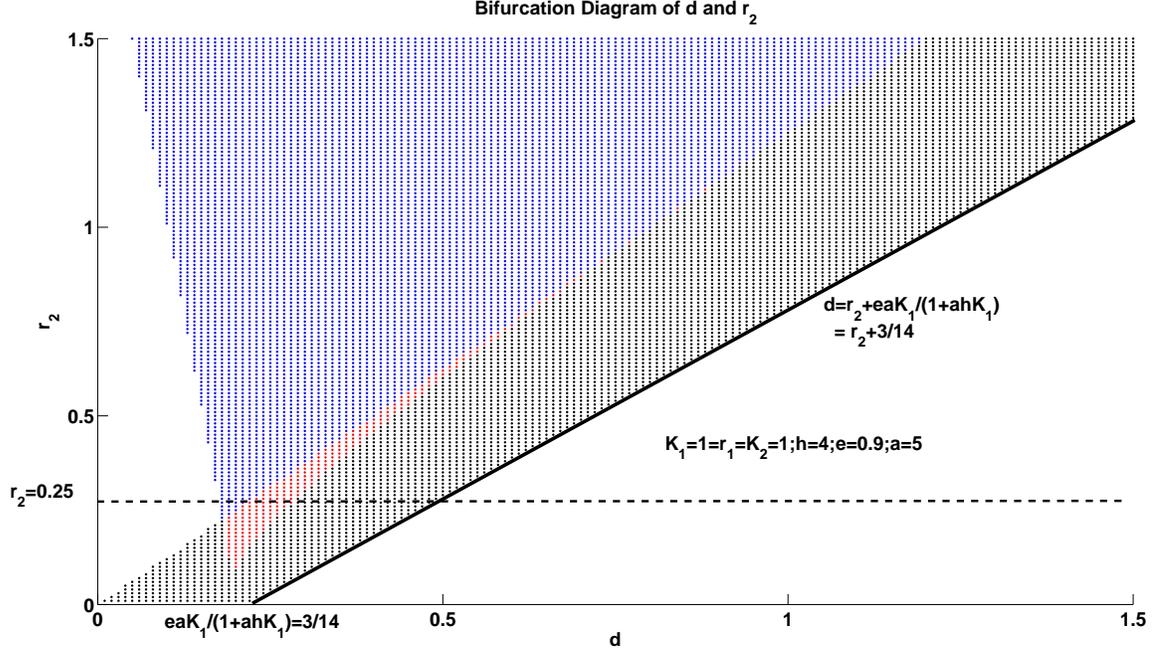}
 \caption{Multiple attractors: Bifurcation diagram of $d$ v.s. $r_2$ for the number of interior equilibrium when $K_1=r_1=K_2=1; h=4; e=0.9;a=5$: the while region means \textbf{no interior equilibrium}; the black region means \textbf{one interior equilibrium}; the blue region means \textbf{two interior equilibria}; and the red region means \textbf{three interior equilibria}. The solid line is $d=r_2+\frac{eaK_1}{1+ahK_1}=r_2+3/14$. \modifyr{When the parasite is facultative, i.e., $d<r_2$, the ecological system \eqref{pp-g} can process none, one, two, or three interior equilibria while parasite is obligate, i.e., $d>r_2$, the system can process none, one, or three interior equilibria}.} \label{fig:interior}
\end{figure}

 For convenience, we fix $K_1=r_1=K_2=1; h=4; e=0.9;a=5$ as a typical example. \modifyr{Figure \ref{fig:interior} illustrates the number of interior equilibria under this set of parameter values by varying both $r_2$ and $d$ from 0.01 to 1.5 where the white region means no interior equilibrium; the black region means one interior equilibrium; the blue region means two interior equilibria; and the red region means three interior equilibria. The solid black line is $d=r_2+\frac{e a K_1}{1+ ah K_1}$ below which there is no interior equilibrium. This is supported by Theorem \ref{th:global}: \emph{Model \eqref{pp-g} has global stability at $E_{10}=(K_1, 0)$ whenever the inequality $d>r_2+\frac{eaK_1}{1+ahK_1}$ holds}. In Figure \ref{fig:interior}, we can observe that Model \eqref{pp-g} has none, one , two, or three interior equilibria when $r_2>d$ while the system has none, one, or three interior equilibria when $r_2<d$. This confirms our results in Theorem \ref{th:interior} and Table \ref{t2:interior}}. \\

\begin{figure}[ht]
\centering
   \includegraphics[scale =.45] {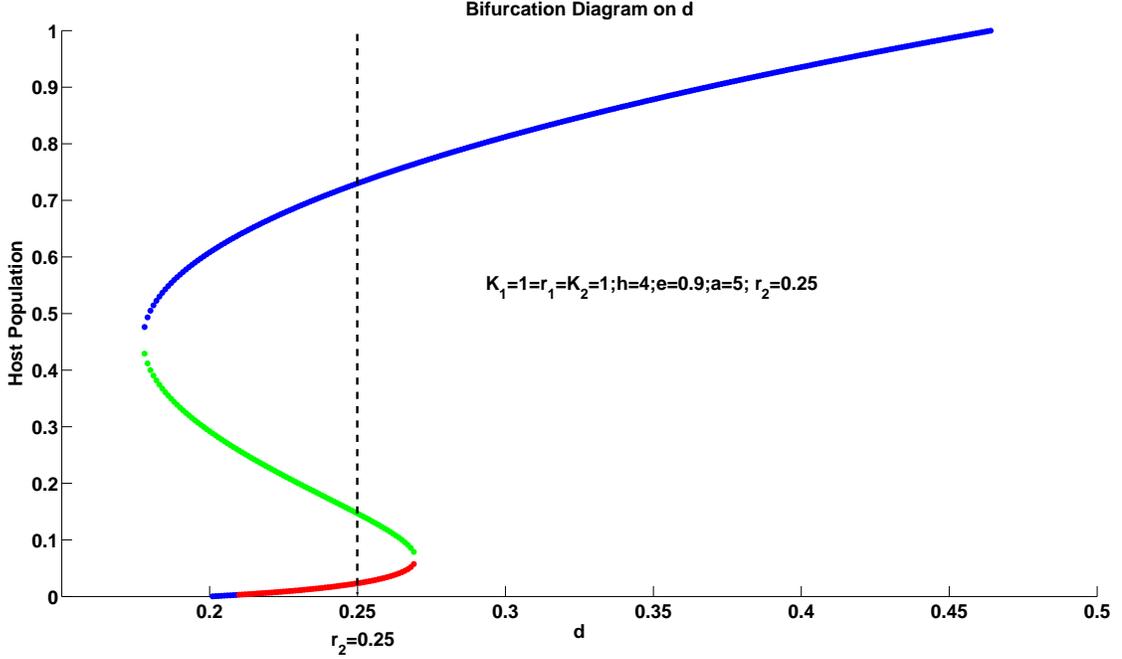}
 \caption{Bifurcation diagrams when $K_1=r_1=K_2=1; r_2=0.2; h=4; e=0.9;a=5$: blue means locally asymptotically stable; green means saddle; and red means source. \modifyr{When the parasite is facultative, i.e., $d<r_2$, the ecological system \eqref{pp-g} can process none, one, two, or three interior equilibria, while the parasite is obligate, i.e., $d>r_2$, the system can process none, one, or three interior equilibria. In addition, large values of $d$ can destabilize the system and lead to the extinction of the parasite.}} \label{fig:d}
\end{figure}

 Let $r_2=0.25$ (see the dashed black line in  Figure \ref{fig:interior}); we perform an additional bifurcation diagram of the stability of interior equilibria with respect to $d$, as shown in Figure \ref{fig:d}, under the same set of parameter values:  blue means locally asymptotically stable; green means saddle; and red means source. \modifyr{Figure \ref{fig:d} shows that when $d$ is small (e.g.,$d\in (0, 0.1)$), there is no interior equilibrium; as $d$ increases (e.g.,$d\in (0.1, 0.2)$), Model \eqref{pp-g} goes through a saddle node bifurcation (around $d=0.1$) which gives one saddle interior equilibrium and one stable interior equilibrium along with the local stable boundary equilibrium $\left(0,K_2\left(1-\frac{d}{r_2}\right)\right)$; $d$ continues to increase (e.g.,$d\in (0.2, 0.215)$), there are three interior equilibria where one is a sink, the second is a saddle and the last is also a sink; further increasing $d$ destabilizes the third interior equilibrium (e.g., there is a hopf-bifurcation occurring at the third interior equilibrium around $d=0.215$); larger $d$ makes the system go through a cusp bifurcation (occurring around $d=0.27$) which leads to only one stable interior equilibrium; at the extreme value of $d$ (i.e., $d>r_2+\frac{e a K_1}{1+ ah K_1}$), the system has global stability at the host-only boundary equilibrium $(K_1,0)$}. \\
 
 One of the interesting questions addressed by the model is what happens to the local stable boundary equilibrium $\left(0,K_2\left(1-\frac{d}{r_2}\right)\right)$ when the system has two interior equilibria after we allow evolution to occur  (i.e., $\sigma_i>0, i=1,2$). Can evolution destabilize this boundary equilibrium and save the host from extinction? We will explore this analytically and numerically in the next section on the co-evolutionary dynamics of \eqref{Ev-pp}.\\

%%%%%%%%%%%%%%%%%%%%%%%%%%%%%%%%%%%%%%%%%%%%%%%%%%%%%%%%%%%%%
%%%%%%%%%%%%%%%%%%%%%%%%%%%%%%%%%%%%%%%%%%%%%%%%%%%%%%%%%%%%%
\section{Co-evolutionary dynamics}

Define $$a_{u_1}=\frac{\partial a(u_1,u_2)}{\partial u_1},\, a_{u_2}=\frac{\partial a(u_1,u_2)}{\partial u_2},\,\, a_{u_1u_1}=\frac{\partial^2 a(u_1,u_2)}{\partial u_1^2},\,\,$$\mbox{ and }\, $$  a_{u_2u_2}=\frac{\partial^2 a(u_1,u_2)}{\partial u_2^2},\,\,  a_{u_1u_2}=\frac{\partial^2 a(u_1,u_2)}{\partial u_1\partial u_2}=\frac{\partial^2 a(u_1,u_2)}{\partial u_1\partial u_2}.$$
Assume that $u_1^*, u_2^*$ are trait values such that
$$K_1'(u_1^*)=K_2'(u_2^*)=a_{u_1}(u_1^*,u_2^*)=a_{u_2}(u_1^*,u_2^*)=d'(u_2^*)=0$$ and let
$$K_1=K_1(u_1^*),\,K_2=K_2(u_2^*),\, a=a(u_1^*,u_2^*),\, \mbox{ and } d=d(u^*_2).$$
Then according to Proposition \ref{p1:be} and Theorem \ref{th:interior}, the co-evolutionary model \eqref{Ev-pp} can have the boundary equilibria $(0,0,u_1^*,u_2^*), E_{x_10u_1u_2}=(K_1, 0,u_1^*,u_2^*),\,E_{0x_2u_1u_2}=\left(0, K_2\left(1-\frac{d}{r_2}\right), u_1^*,u_2^*\right)$ and potentially multiple interior equilibria $(x_1^*, x_2^*,u_1^*,u_2^*)$ depending on the values of $K_i, r_i, i=1,2$ and $d, a$. In general, the following theorem provides the existence and local stability results on the facultative-parasite-only equilibrium $E_{0x_2u_1u_2}$ and the host-only equilibrium $E_{x_10u_1u_2}$ for the co-evolutionary model \eqref{Ev-pp}:\\
%%%%%%%%%%%%%%%%%%%%%%%%%%%%%%%%%%%%%%%%%%%%
\begin{theorem}\label{p2:be-ev}
The host-only equilibrium $E_{x_10u_1u_2}=(K_1, 0,u_1^*,u_2^*)$ exists if 
$$K_1'(u_1^*)=0\,\, \mbox{ and }\,\, d'(u_2^*)=\frac{eK_1(u_1^*)a_{u_2}(u_1^*,u_2^*)}{(1+h a(u_1^*,u_2^*)K_1(u_1^*))^2}$$ 
and it is locally asymptotically stable if the following inequalities hold
$$0<\frac{ea(u_1^*,u_2^*)K_1(u_1^*)}{1+ha(u_1^*,u_2^*)K_1(u_1^*)}<d(u_2^*)-r_2,\,\,K_1''(u_1^*)<0,\,\,$$ \mbox{ and } $$\frac{eK_1(u_1^*)\left[a_{u_2u_2}(u_1^*,u_2^*)(1+haK_1(u_1^*))-2 K_1(u_1^*)a^2_{u_2}(u_1^*,u_2^*)\right]}{(1+h a(u_1^*,u_2^*)K_1(u_1^*))^3}<d''(u_2^*).$$
The parasite-only equilibrium $E_{0x_2u_1u_2}=\left(0, K_2\left(1-\frac{d}{r_2}\right)u_1^*,u_2^*\right)$ exists if 
$$ a_{u_1}(u_1^*,u_2^*)=0,\,d(u_2^*)<r_2\mbox{ and }d'(u_2^*)=\frac{r_2K'_2(u_2^*)}{K_2(u_2^*)}$$ and it is locally asymptotically stable if the following inequalities hold
$$\frac{r_1}{a(u_1^*,u_2^*)}<K_2(u_2^*)\left(1-\frac{d(u_2^*)}{r_2}\right),\,\,a_{u_1u_1}(u_1^*,u_2^*)>0,\,\,$$ \mbox{ and } $$\frac{r_2\left(1-\frac{d(u_2^*)}{r_2}\right)\left[K_2(u_2^*) K''_2(u_2^*)-2(K'_2(u_2^*))^2\right]}{K_2(u_2^*)^2}<d''(u_2^*).$$

\end{theorem}
%%%%%%%%%%%%%%%%%%%%%%%%%%%%%%%%%%%%%%%%%%%%%%%%%%%%%
\noindent\textbf{Notes:} According to Proposition \ref{p1:be}, the extinction equilibrium $E_{00}$ is alway unstable for the ecological model \eqref{pp-g}. The host-only equilibrium $E_{10}$ is locally asymptotically stable when $r_2+\frac{eaK_1}{1+ahK_1}<d$. This condition implies that social parasite should be obligate, i.e., $r_2<d$. The facultative-social-parasite-only equilibrium $E_{01}$ is locally asymptotically stable when the inequality $\frac{r_1}{a K_2}<1-\frac{d}{r_2}$ holds.  The results of Theorem \ref{p2:be-ev} imply that the local stability of the facultative-social-parasite-only equilibrium $E_{0x_2u_1u_2}$ and the host-only equilibrium $E_{x_10u_1u_2}$ are determined by the concavity of the trait function $K_i(v_i), i=1,2$ and $a(v_1,v_2)$ evaluated at the equilibrium under conditions that the equilibrium is Ecologically Stable (i.e., ES). If the  death  rate due to parasitism $d$ is independent of the trait value $u_2$, then we have the following corollary: \\

%%%%%%%%%%%%%%%%%%%%%%%%%%%%%%%%%%%
\begin{corollary}\label{cp2:be-ev}
The host-only equilibrium $E_{x_10u_1u_2}$ exists if $K_1'(u_1^*)=a_{u_2}(u_1^*,u_2^*)=0$,
and it is locally asymptotically stable if the following inequalities hold
$$0<\frac{ea(u_1^*,u_2^*)K_1(u_1^*)}{1+ha(u_1^*,u_2^*)K_1(u_1^*)}<d-r_2,\,\,K_1''(u_1^*)<0,\,\, \mbox{ and } a_{u_2u_2}(u_1^*,u_2^*)<0.$$
The parasite-only equilibrium $E_{0x_2u_1u_2}$ exists if $K_2'(u_2^*)=a_{u_1}(u_1^*,u_2^*)=0$,
 and it is locally asymptotically stable if the following inequalities hold
$$\frac{r_1}{a(u_1^*,u_2^*) }<K_2(u_2^*)\left(1-\frac{d}{r_2}\right),\,\,a_{u_1u_1}(u_1^*,u_2^*)>0,\,\, \mbox{ and }K''_2(u_2^*)<0.$$

\end{corollary}
%%%%%%%%%%%%%%%%%%%%%%%%%%%%%%%%%%%%%%%%%%%%%%%%%%%%%%%%%%%%%%%%%%%%%%%%
%\noindent\textbf{Applications:} 
\subsection{Boundary equilibria and ESS}
 \modifyr{The trait functions can have many forms, such as Gaussian distributions, polynomial, exponential functions (Abrams 1990; Bergelson \emph{et al} 2001; Mostowy and Engelst\"adter 2011; Nuismer \emph{et al} 2012; Landi \emph{et al} 2013). In this subsection, we apply the results of Theorem \ref{p2:be-ev} and its Corollary \ref{cp2:be-ev} to some specific trait functions to explore how these functions affect whether or not the boundary equilibrium $E_{x_10u_1u_2}$ or  $E_{0x_2u_1u_2}$ can have ESS. More specifically, we assume that 
$$K_i(v_i)=K_{0i} e^{-\frac{v_i^2 (v_i-c)^2}{2\sigma_{K_i}^2}}, \,\,a(v_1,v_2)= a_0e^{-\frac{(v_1-v_2+c)^2(v_1-v_2-c)^2}{2\sigma_a^4}}, \mbox{ and } d(v_2)=d_{0} e^{-\frac{(v_2-c)^2 (v_2+c)^2}{2\sigma_{d}^2}}, i=1,2$$where $c>0$ and $v_i\in [0,c]$. These chosen trait functions are modified from Gaussian distributions. }\\
%%%%%%%%%%%%%%%%%%%%%%%%%%%%%%%%%%%%%%%%%%%%%%%%%%%%%%%%%%%%%%

\subsubsection{The fixed parasite death rate\\ }
In this subsection, we apply the results of Corollary \ref{cp2:be-ev} to particular trait functions when the death rate $d$ of the parasite due to attacking all potential hosts is independent of the trait $u_2$. 
Let trait values set be $\mathbb U_i=[0,c], i=1,2.$ We assume that $K_i(v_i)=K_{0i} e^{-\frac{v_i^2 (v_i-c)^2}{2\sigma_{K_i}^2}}, i=1,2$ and $a(v_1,v_2)= a_0e^{-\frac{(v_1-v_2+c)^2(v_1-v_2-c)^2}{2\sigma_a^4}}$ which gives follows:
$$K_i'0)=K_i(c)=K_{01},\,\,a(0,0)=a(c,c)=a_0 e^{-\frac{c^4}{4\sigma_a^2}},\,\,a(0,c)=a(c,0)=a_0,i=1,2$$
$$K_i'(0)=K_i'(c)=a_{v_i}(0,0)=a_{v_i}(c,c)=a_{v_i}(0,c)=a_{v_i}(c,0),i=1,2$$  
$$K_i''(0)=K_i''(c)=-\frac{K_{0i}c^2}{\sigma_{K_i}^2}<0,\,\,a_{v_iv_i}(0,0)=a_{v_iv_i}(c,c)=\frac{a_0c^2 e^{-\frac{c^2}{\sigma_a^4}}}{\sigma_a^2}>0,\,\,$$\mbox{ and }$$a_{v_iv_i}(0,c)=a_{v_iv_i}(c,0)=-\frac{2a_0c^2}{\sigma_a^2}<0.$$
This implies that the trait dynamics of the co-evolutionary model \eqref{Ev-pp} are positively invariant in the trait space $\mathbb U_1\times\mathbb U_2=[0,c]\times[0,c]$.\\

\noindent\textbf{The parasite-only equilibrium:} According to Corollary \ref{cp2:be-ev}, the parasite-only equilibrium $E_{0x_2u_1u_2}=\left(0,K_{02}\left(1-\frac{d}{r_2}\right),u^*_1,u^*_2*\right)$ is locally asymptotically stable (i.e., CS) if the following conditions hold
{\small\bae\label{exESS1}K_2'(u^*_2)=a_{u_1}(u^*_1,u^*_2)=0,\,\,\frac{r_1}{1-\frac{d}{r_2}}<a_0K_{02} e^{-\frac{c^4}{4\sigma_a^4}}, a_{u_1u_1}(u^*_1,u^*_2)=\frac{a_0c^2 e^{-\frac{c^2}{\sigma_a^4}}}{\sigma_a^2}>0,\, \mbox{ and }K''_2(u^*_2)<0\eae} where $(u^*_1,u^*_2)=(0,0)$ or $(u^*_1,u^*_2)=(c,c)$. Then the fitness functions of host and parasite at $E_{0x_2u_1u_2}$ are represented as follows:
$$\begin{array}{lcl}
G_1(v_1,u^*,x^*)&=&r_1-a(v_1,u_2^*) K_{02}\left(1-\frac{d}{r_2}\right)=r_1-a_0K_{02}\left(1-\frac{d}{r_2}\right)e^{-\frac{(v_1-u_2^*+c)^2(v_1-u_2^*-c)^2}{4\sigma_a^2}}\\
G_2(v_2,u^*,x^*)&=&-d+r_2 \left(1-\frac{K_{02}\left(1-\frac{d}{r_2}\right)}{K_2(v_2)}\right)=-d+r_2 \left(1-\left(1-\frac{d}{r_2}\right)e^{\frac{v_2^2(v_2-c)^2}{2\sigma_{K_2}^2}}\right)\end{array}$$which gives
$$\max_{v_1\in\mathbb U_1}\{G_1(v_1,u^*,x^*)\}= G_1(u_1^*,u^*,x^*)=r_1-a_0K_{02}\left(1-\frac{d}{r_2}\right)e^{-\frac{c^4}{4\sigma_a^2}}<0$$\mbox{ and }
$$\max_{v_2\in\mathbb U_2}\{G_2(v_2,u^*,x^*)\}= G_2(u_2^*,u^*,x^*)=-d+r_2 \left(1-\left(1-\frac{d}{r_2}\right)\right)=0.$$
Therefore, according to the definition of ESS (i.e., \eqref{ESS-max}), we can conclude that the two parasite-only equilibrium $\left(0,K_{02}\left(1-\frac{d}{r_2}\right),0,0\right)$ and $\left(0,K_{02}\left(1-\frac{d}{r_2}\right),c,c\right)$ of the co-evolutionary model \eqref{Ev-pp} are ESS.\\

\noindent\textbf{Biological scenarios:} The case studied above suggest that  the parasite-only-equilibrium $E_{0x_2u_1u_2}$ can have two ESS strategies $(u^*_1,u^*_2)=(0,0)$ and $(u^*_1,u^*_2)=(c,c)$ when the parasite is facultative (i.e., $d<r_2$). This may describe the case of facultative slavemakers, such as  \emph{Formica subnuda} whose colonies can survive as slaveless (Savolainen and Deslippe 1996). \\\\%Another example could be facultative social parasitism upon bumble bee (\emph{Bombus spp.}) that enter a nest and kill the resident queen may stay and adopt the brood (Voveikov 1953; Fisher 1987).\\\\

\noindent\textbf{The host-only equilibrium:} According to Corollary \ref{cp2:be-ev}, the host-only equilibrium $E_{x_10u_1u_2}=\left(K_{01},0,u^*_1,u^*_2\right)$ is locally asymptotically stable if the following conditions hold
{\small\bae\label{exESS2}K_1'(u^*_1)=a_{u_2}(u^*_1,u^*_2)=0,\,\,\frac{ea_0K_{01}}{1+ha_0K_{01}}<d-r_2 , a_{u_2u_2}(u^*_1,u^*_2)=-\frac{2a_0c^2}{\sigma_a^2}<0,\, \mbox{ and }K''_2(u^*_2)<0\eae} where $(u^*_1,u^*_2)=(0,c)$ or $(u^*_1,u^*_2)=(c,0)$. Then the fitness functions of host and parasite at $E_{x_10u_1u_2}$ are represented as follows:
$$\begin{array}{lcl}
G_1(v_1,u^*,x^*)&=&r_1\left(1-\frac{K_{01}}{K_{01}e^{-\frac{v_1^2(v_1-c)^2}{2\sigma_{K_1}^2}}}\right)=r_1\left(1-e^{\frac{v_1^2(v_1-c)^2}{2\sigma_{K_1}^2}}\right)\\
G_2(v_2,u^*,x^*)&=&\frac{ea(u_1^*,v_2)K_{01}}{1+ha(u_1^*,v_2)K_{01}}-d+r_2 =-d+r_2+\frac{eK_{01}a_0e^{-\frac{(u_1^*-v_2+c)^2(u_1^*-v_2-c)^2}{2\sigma_a^4}}}{1+hK_{01}a_0e^{-\frac{(u_1^*-v_2+c)^2(u_1^*-v_2-c)^2}{2\sigma_a^4}}}\end{array}$$which gives
$$\max_{v_1\in\mathbb U_1}\{G_1(v_1,u^*,x^*)\}= G_1(u_1^*,u^*,x^*)=0$$\mbox{ and }
$$\max_{v_2\in\mathbb U_2}\{G_2(v_2,u^*,x^*)\}= G_2(u_2^*,u^*,x^*)=-d+r_2+\frac{ea_0K_{01}}{1+ha_0K_{01}}<0.$$
Therefore, according to \eqref{ESS-max}, the two host-only equilibria $\left(K_{01},0,0,c\right)$ and $\left(K_{01},0,c,0\right)$ of the co-evolutionary model \eqref{Ev-pp} are ESS.\\

\noindent\textbf{Biological scenarios:} The cases studied above suggest that for the trait functions given, the host-only equilibrium $E_{x_10u_1u_2}$ of the co-evolutionary model \eqref{Ev-pp} can have two ESS $(u^*_1,u^*_2)=(0,c)$ and $(u^*_1,u^*_2)=(c,0)$ when the parasite is obligate (i.e., $d>r_2$). This can be classified as one co-evolutionary outcome when the host successfully resists invasion by the parasite, resistance occurs via effective front-line defenses (Kilner and Langmore 2011). For example, Ortolani and Cervo (2010) show that Polistes dominulus foundresses are now so large and  aggressive that they consistently defend their nests from attack by the brood parasite P. sulcifer, and that parasitism is rarely seen.\\%  The extinction of the large blue butterfly in 1979 in Britain may be an another example of the host-only equilibrium having ESS strategies. The large blue butterfly exhibits a unique obligate parasitic relationship with \emph{Myrmica} ants, specifically with \emph{Myrmica sabuleti} (Als \emph{\emph{et al}} 2004). \\

%%%%%%%%%%%%%%%%%%%%%%%%%%%%%%%%%%%%%%%%%%%%%%%%%%%%%%%%%%%%%%

\subsubsection{The parasite death rate depending on its trait\\}
Our  application of Corollary \ref{cp2:be-ev} in the previous subsection indicates that both boundary equilibria $E_{x_10u_1u_2}$ and $E_{0x_2u_1u_2}$ cannot simultaneously  have ESS strategies, due to the fact that $d$ and $r_2$ are fixed. The interesting questions become, if the death rate $d$ of the parasite also depends on its trait $u_2$, \textbf{can Model  \eqref{Ev-pp} have both boundary equilibria $E_{x_10u_1u_2}$ and $E_{0x_2u_1u_2}$ being locally asymptotically stable (i.e., CS) under different trait values of $u_2$? If it can,  is it possible for these two boundary equilibria to have ESS?}\\

\begin{figure}[ht]
\centering
  \subfigure[The fitness function of host (i.e., $G_1(v_1,u^*,x^*)$) with $u^*=(0,c)$ and $x^*=(K_{01},0)$]{
   \includegraphics[scale =.4] {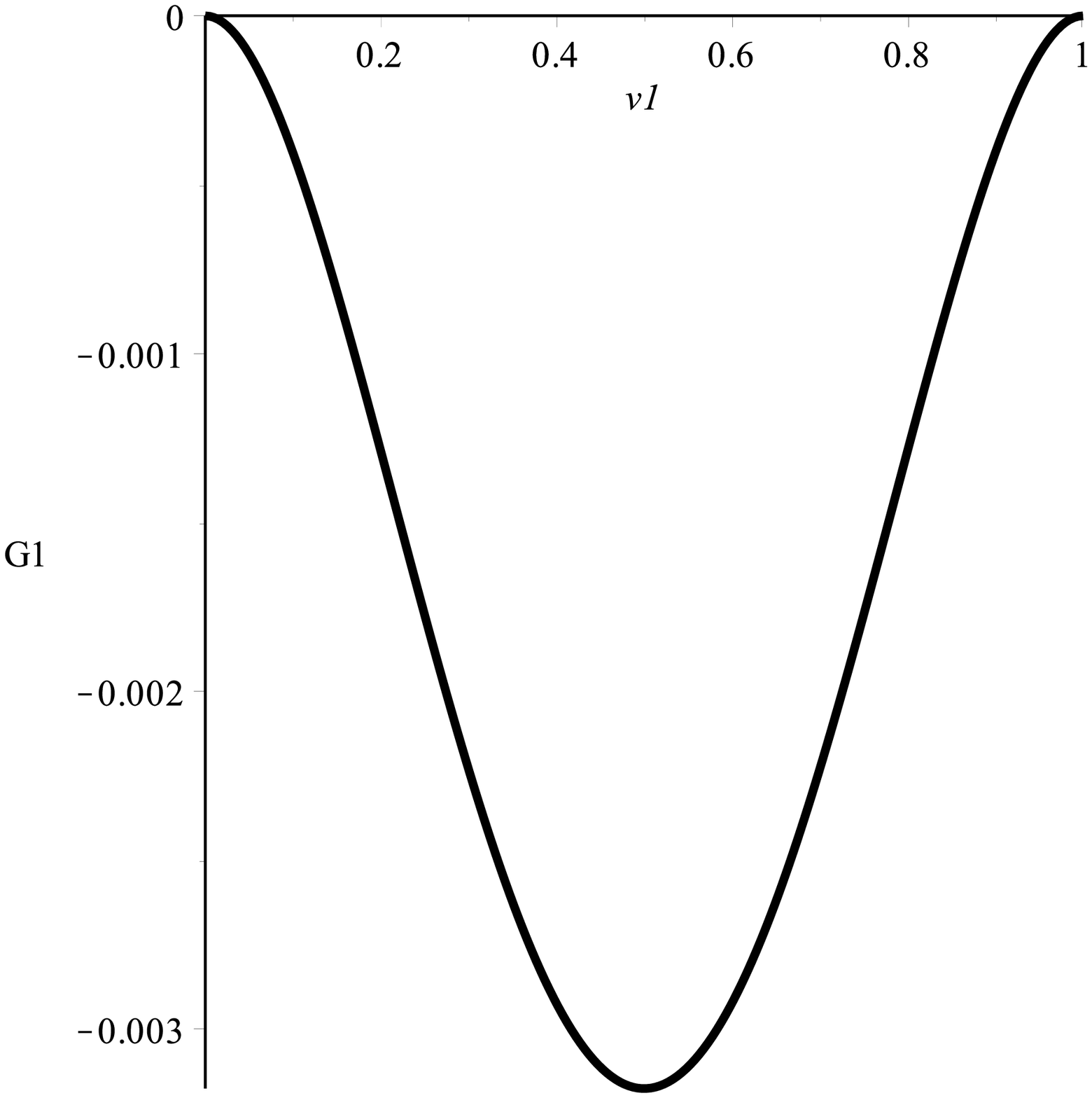}
 \label{fig:G1K1}}
  \subfigure[The fitness function of parasite (i.e., $G_2(v_2,u^*,x^*)$) with $u^*=(0,c)$ and $x^*=(K_{01},0)$]{
  \includegraphics[scale =.4] {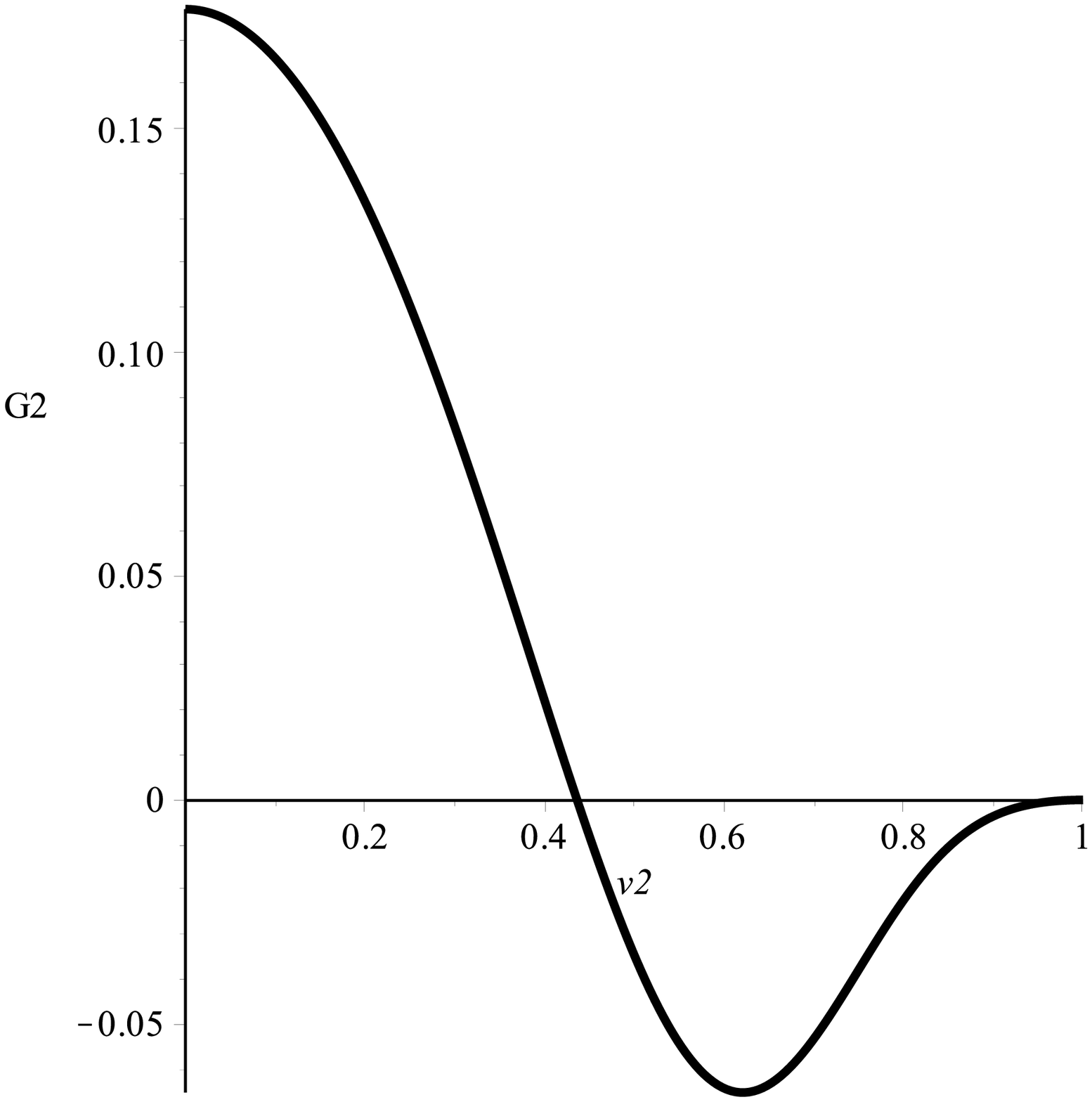}
 \label{fig:G2K1}}
 \caption{The fitness functions of host (i.e., $G_1(v_1,u^*,x^*)$) and of parasite (i.e., $G_2(v_2,u^*,x^*)$) for the co-evolutionary model \eqref{Ev-pp} when $r_1=0.1,\,\,r_2=1.5<d_0=2.1,\,\,a_0=2,\,\,e=0.9,\,\,h=c=\sigma_{K_1}=\sigma_{K_2}=K_{01}=1,\,\,K_{02}=100,\,\,\sigma_d=1.05,\,\,\sigma_a=0.56$. \modifyr{The strategy $u^*=(0,0)$ is not an ESS, because there are trait values $v_2$ such that $G_2(v_2,u^*,x^*)>0$. However, the co-evolutionary system has local stability at $(u^*,x*)$}. }\label{fig1:G1G2}
\end{figure}

\begin{figure}[ht]
\centering
  \subfigure[The fitness function of host (i.e., $G_1(v_1,u^*,x^*)$) with $u^*=(0,0)$ and $x^*=\left(0,K_{02}\left(1-\frac{d_0e^{-\frac{c^4}{2\sigma_d^2}}}{r_2}\right)\right)$]{
 \includegraphics[scale =.4] {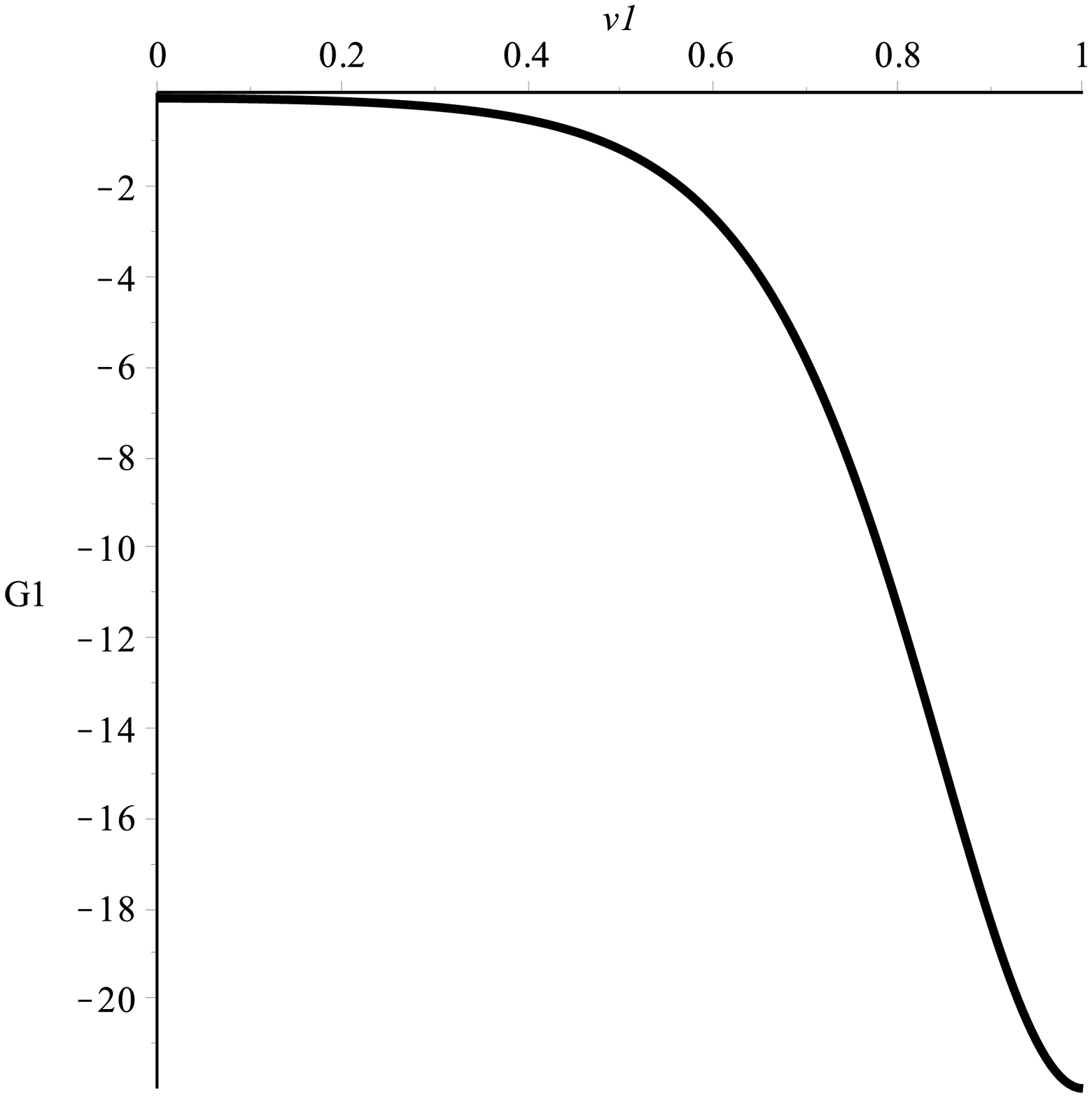}
 \label{fig:G1K2}}
  \subfigure[The fitness function of parasite (i.e., $G_2(v_2,u^*,x^*)$) with $u^*=(0,0)$ and $x^*=\left(0,K_{02}\left(1-\frac{d_0e^{-\frac{c^4}{2\sigma_d^2}}}{r_2}\right)\right)$]{
   \includegraphics[scale =.4] {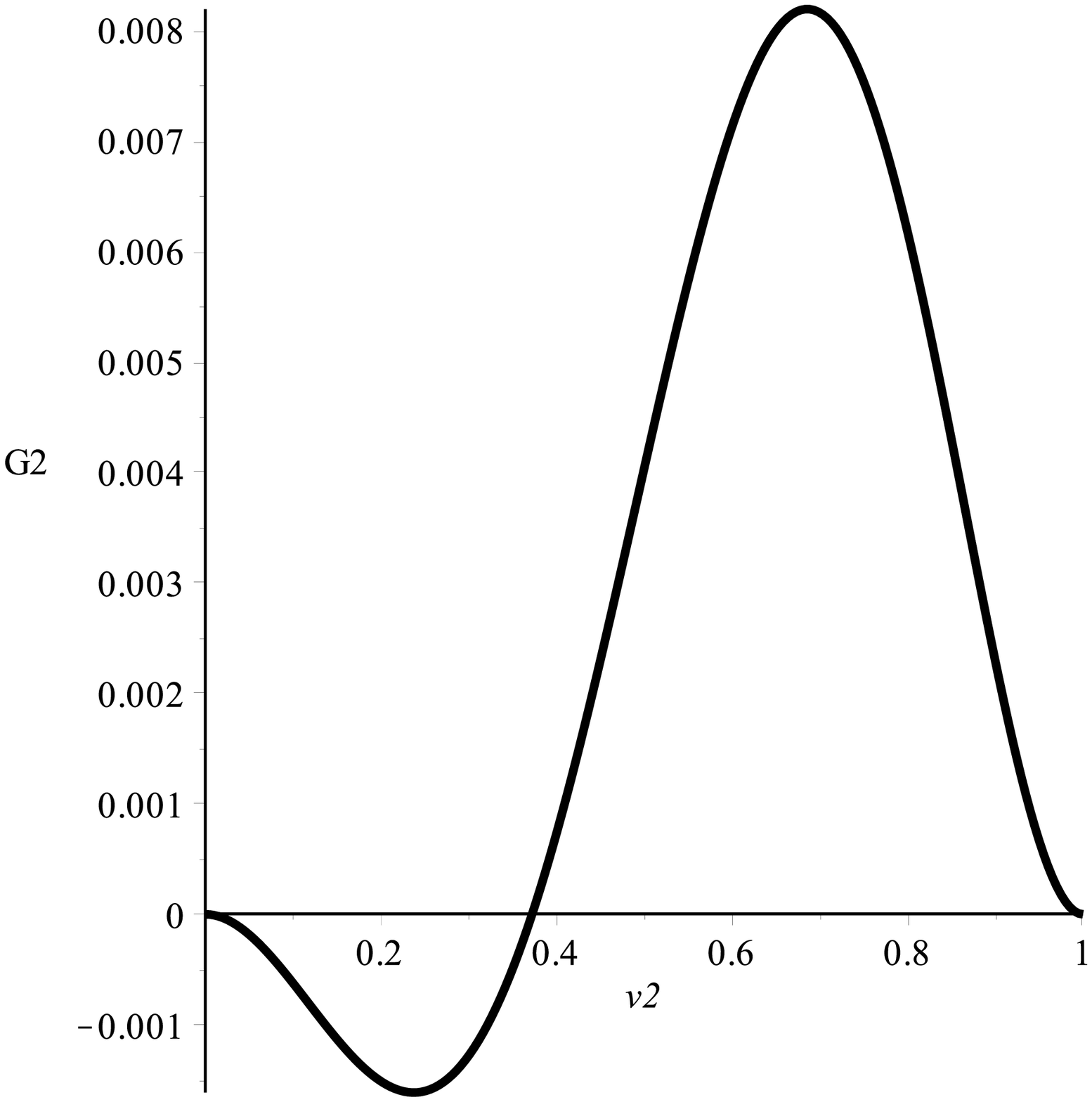}
 \label{fig:G2K2}}
 \caption{The fitness functions of host (i.e., $G_1(v_1,u^*,x^*)$) and of parasite (i.e., $G_2(v_2,u^*,x^*)$) for the co-evolutionary model \eqref{Ev-pp} when $r_1=0.1,\,\,r_2=1.5<d_0=2.1,\,\,a_0=2,\,\,e=0.9,\,\,h=c=\sigma_{K_1}=\sigma_{K_2}=K_{01}=1,\,\,K_{02}=100,\,\,\sigma_d=1.05,\,\,\sigma_a=0.56$. \modifyr{The strategy $u^*=(0,c)$ is not ESS for since since there are trait values $v_2$ such that $G_2(v_2,u^*,x^*)>0$. However, the co-evolutionary system has local stability at $(u^*,x*)$}. }\label{fig2:G1G2}
\end{figure}

To investigate these questions, we should apply the results of Theorem \ref{p2:be-ev} by letting $d$ depend on the trait value $u_2$. To continue our study, we choose $d(v_2)=d_{0} e^{-\frac{(v_2-c)^2 (v_2+c)^2}{2\sigma_{d}^2}}$ which has the following properties:
$$d(0)= d_{0}e^{\frac{-c^4}{2\sigma_{d}^2},}\,\, d(c)=d_{0},\,\,d'(0)=d'(c)=0, \,\,d''(0)=\frac{2d_0c^2 e^{-\frac{c^2}{\sigma_d^2}}}{\sigma_d^2}>0,\,\, \mbox{ and }d''(c)=-\frac{4d_0c^2}{\sigma_d^2}<0.$$
If we assume that $K_i(v_i),i=1,2$ and $a(v_1,v_2)$ have the same trait functions as before, then we can conclude that the trait dynamics of the co-evolutionary model \eqref{Ev-pp} are positively invariant in the trait space $\mathbb U_1\times\mathbb U_2=[0,c]\times[0,c]$; and the equilibrium trait values $(u^*_1,u^*_2)$ are $u^*_i=0$ or $c$ for both $i=1,2$. According to Theorem \ref{p2:be-ev}, we can conclude that both boundary equilibria $E_{x_10u_1u_2}=(K_{01},0,0,c)$ and $E_{0x_2u_1u_2}=\left(0,K_{02}\left(1-\frac{d_0}{r_2}\right),0,0\right)$ have convergence stability if the following conditions hold:
$$\frac{d_0\sigma_a^4 (1+a_0 h K_{01})}{\sigma_d^2}<\frac{e a_0 K_{01}}{1+a_0 h K_{01}}<d_0-r_2 \mbox{ and } r_1<a_0 K_{02} e^{-\frac{c^4}{2\sigma_a^4}}\left(1-\frac{d_0e^{-\frac{c^4}{2\sigma_d^2}}}{r_2}\right).$$
 For example, if we let 
 $$r_1=0.1,\,\,r_2=1.5<d_0=2.1,\,\,a_0=2,\,\,e=.9,\,\,h=c=\sigma_{K_1}=\sigma_{K_2}=K_{01}=1,\,\,K_{02}=100,\,\,\sigma_d=1.05,\,\,\sigma_a=.56,$$
 then \modifyr{both $E_{x_10u_1u_2}=(K_{01},0,0,c)$ and $E_{0x_2u_1u_2}=\left(0,K_{02}\left(1-\frac{d_0}{r_2}\right),0,c\right)$ are locally asymptotically stable (CS). However,  the equilibrium trait $(u^*_1,u^*_2)=(0,c)$ is not ESS for $E_{x_10u_1u_2}=(K_{01},0,0,c)$; and $(u^*_1,u^*_2)=(0,0)$ is not ESS for $E_{0x_2u_1u_2}=\left(0,K_{02}\left(1-\frac{d_0e^{-\frac{c^4}{2\sigma_d^2}}}{r_2}\right),0,0\right)$. See the fitness functions $G_i(v_i, u^*,x^*), i=1,2$ for these two boundary equilibria shown in Figure \ref{fig1:G1G2}-\ref{fig2:G1G2}:  Figure \ref{fig:G1K1} and \ref{fig:G1K2} show that the strategies $(u^*_1,u^*_2)=(0,c)$ and $(u^*_1,u^*_2)=(0,0)$ are ESS for  the host at the boundary equilibrium $(K_{01},0)$, $\left(0,K_{02}\left(1-\frac{d_0e^{-\frac{c^4}{2\sigma_d^2}}}{r_2}\right)\right)$, respectively. However,  Figure \ref{fig:G2K1} and \ref{fig:G2K2} show that these strategies are not ESS for the parasite at these boundary equilibria, since there are trait values $v_2$ such that $G_2(v_2,u^*,x^*)>0$}.\\\\

 \noindent\textbf{Biological scenarios and implications:} The case studied above suggests that a host and parasite can employ multiple strategies that generate convergence stability but not ESS. This may be due to the interactions of different strategies  and counter-strategies acquired at different stages of co-evolution (Kilner and Langmore 2011). For example, parasites may initially be facultative, then acquire host-specific signatures, but through counter-selection by hosts may subsequently revert to  genetic signatures expressed before parasitism. As one possible example, among the insects, selection to become chemically insignificant may have facilitated the ability to acquire host hydrocarbon signatures previously after parasitism (Kilner and Langmore 2011). Such a case could be applied to slave-making ants, which can be either obligate  social parasites, depending on enslaved hosts ants throughout their whole lives (Topoff and Zimmerli 1991; Ruano \emph{\emph{et al}.} 2013) or alternatively facultative slave-makers. Facultative slave-making ants, like those in the \emph{Formica sanguinea} complex, may engage in slave making, but individual colonies are able to revert to producing their own workers if parasitized workers are removed (Topoff and Zimmerli 1991). They could represent an intermediate parasitic group, between freeliving species on the one hand, and obligate slave-making species on the other.\\\\% In laboratory tests, slaves were removed from colonies of \emph{Formica sanguinea} and \emph{Polyergus rufescens}: The behavior of \emph{F. sanguinea} changed dramatically within 30 days of slave removal, with workers becoming self-sufficient at feeding and brood care; Workers of \emph{Polyergus}, by contrast, were unable to care for their brood, and experienced high mortality (Topoff and Zimmerli 1991).\\\\

  %%%%%%%%%%%%%%%%%%%%%%%%%%%%%%%%%%%%%%%%%%%%%%%%%%%%%%%%%%%%%
\subsection{Trait functions follow Gaussian distributions\\}
To continue our analytical study, we focus on the co-evolutionary dynamics of Model \eqref{Ev-pp} for chosen trait functions of $K_i(v_i)$ and $a(v_1,v_2)$. For the rest of the section, we let the  death rate $d$ of the parasite $x_2$ be independent of the trait $u_2$; and we assume the following functional forms for the carrying capacity $K_i(v_i)$ and the parasitism efficiency $a(v_1, v_2)$:
\bae\label{E-pp-ka}
\begin{array}{lcl}
K_i(v_i)&=&K_{0i} e^{-\frac{v_i^2}{2\sigma_{K_i}^2}}\\\\
a(v_1,v_2)&=& a_0e^{-\frac{(v_1-v_2)^2}{2\sigma_a^2}}
\end{array}
\eae $a(v_1, v_2)$ denotes the parasitism efficiency of a parasite with phenotypic trait $v_2$ on host individuals with phenotypic trait $v_1$ with the assumption that the stronger host-parasite interactions are, the more similar host and parasite traits are. The symmetric form of $a$ in \eqref{E-pp-ka} has been previously used in the study of character displacement by Taper and Case (1992), which assumes that the parasitism efficiency is normally distributed around a maximum value of $a_0$ with a variance $\sigma_a^2$ as a function of trait difference in the parasite and host $v_1-v_2$. The larger value of $\sigma_a^2$, the greater sensitivity of the parasitism efficiency $a$ changes with respect to the changes in trait difference $v_1-v_2$. There are other alternative functions of $a$ (see Dieckmann and Marrow 1995; Doebeli and Dieckmann 2000; Zu \emph{\emph{et al}} 2007), for example,  the asymmetrical predation efficiency $a$ which has been previously used in the study of character displacement can be described as
\bae\label{as-a}
\begin{array}{lcl}
a(v_1,v_2)&=& a_0e^{\frac{(\sigma_a^2\beta)^2}{2}}e^{-\frac{(v_1-v_2+\sigma_a^2\beta)^2}{2\sigma_a^2}}
\end{array}
\eae  where $\beta\neq 0$ (Doebeli and U. Dieckmann 2000).  We also assume that the resource availability for the host and parasite varies with their phenotypic trait $v_i$, which follows a Gaussian distribution $N(0, \sigma_{K_i}^2)$. This distribution assumes that the maximum inherent equilibrium level $K_{0i}$ for each individual single species $x_i$, using strategy $v_i$, in the absence of other species, is attained at trait $v_i = 0$, and that $K_i(v_i)$ is normally distributed around $v_i = 0$ with variance $\sigma_{K_i}^2$.  The larger variance $\sigma_{K_i}^2$, the greater sensitivity of the carrying capacity $K_i(v_i)$ changes with respect to the changes in the trait $v_i$. \modifyr{Fix the trait value of the parasite $u_2$, the formulations of $K_{1}(v_1)$ and $a(v_1,u_2)$ indicate that increases in the host's trait value $\vert v_1\vert$ result in the decreased values of $K_1(v_1)$ and $a(v_1,u_2)$ when $v_1<0$ or $v_1>u_2$ but can result in the increased values of $a$ when $v_1\in (0, u_2)$. This implies that there is a trade-off between the carrying capacity of the host and the parasitism efficiency for a certain range of $v_1$ when $v_1<0$ or $v_1>u_2$. On the other hand, for a fixed host's trait value $u_1$, there is no such trade-off for the parasite when $v_2<0$ or $v_2>u_1$ but a trade-off exists when $v_2\in (0,u_1)$}.\\

As a consequence, the co-evolutionary dynamics of the monomorphic resident host and parasite populations with residence traits of $u_1$ and $u_2$ in the host, parasite, respectively, are given by the following set of nonlinear equations:
\bae\label{E-pp-g}
\begin{array}{lcl}
\frac{dx_1}{dt}&=&x_1G_1(v_1,u,x)\vert_{v_1=u_1}=x_1 H_1(u_1,u_2,x_1, x_2)=x_1 \left[r_1 (1-\frac{x_1}{K_1(u_1)})-\frac{a(u_1,u_2) x_2}{1+h a(u_1,u_2)x_1}\right]\\\\
\frac{dx_2}{dt}&=&x_2G_2(v_2,u,x)\vert_{v_2=u_2}=x_2 H_2(u_1,u_2,x_1,x_2)=x_2\left[\frac{e a(u_1, u_2) x_1}{1+ha(u_1,u_2)x_1}-d+r_2 (1-\frac{x_2}{K_2(u_2)})\right]\\\\
\frac{du_1}{dt}&=&\sigma_1^2\frac{\partial G_1(v_1,u,x)}{\partial v_1}\big\vert_{v_1=u_1}=\sigma_1^2 \frac{\partial H_1(u_1,u_2,x_1,x_2)}{\partial u_1}=\sigma_1^2\left[-\frac{r_1 x_1 u_1}{\sigma_{K_1}^2K_1(u_1)}+\frac{(u_1-u_2)x_2a(u_1,u_2)}{\sigma_a^2(1+h a(u_1,u_2)x_1)^2}\right],\\\\
\frac{du_2}{dt}&=&\sigma_2^2\frac{\partial G_2(v_2,u,x)}{\partial v_2}\big\vert_{v_2=u_2}=\sigma_2^2 \frac{\partial H_2(u_1,u_2,x_1,x_2)}{\partial u_2}=\sigma_2^2\left[-\frac{r_2 x_2 u_2}{\sigma_{K_2}^2K_2(u_2)}+\frac{e(u_1-u_2)x_1a(u_1,u_2)}{\sigma_a^2(1+h a(u_1,u_2)x_1)^2}\right],
\end{array}
\eae where we take $a$ with the symmetric form given in \eqref{E-pp-ka}. Since the death rate of the parasite $x_2$ does not change over time, thus the difference of $r_2$ and $d$ determines whether the parasite is facultative or obligate: the parasite is facultative when $r_2>d\geq 0$, and is obligate when $d>r_2\geq 0$. Our model \eqref{E-pp-g} allows us to investigate the ecological and evolutionary dynamics of the host and parasite interactions, where the parasite can be either obligate or facultative. In particular, Model \eqref{E-pp-g} includes the special case which is studied by \modifyr{Zu \emph{\emph{et al}.} (2007)}.\\

Assume that $(x^*_1, x^*_2, u^*_1,u^*_2)$ is an equilibrium of the co-evolutionary host-parasite model \eqref{E-pp-g}. We say $(x^*_1, x^*_2, u^*_1,u^*_2)$ is a boundary equilibrium if $x^*_1x^*_2=0$ while it is an interior equilibrium if $x^*_1x^*_2>0$.  The equilibrium $(x^*_1, x^*_2, u^*_1,u^*_2)$ satisfies the following equations
\bae\label{eq-x1}
\frac{dx_1}{dt}&=&0\Rightarrow x_1=0 \mbox{ or } r_1 \left(1-\frac{x_1}{K_1(u_1)}\right)=\frac{a(u_1,u_2) x_2}{1+h a(u_1,u_2)x_1}\\
\label{eq-x2}
\frac{dx_2}{dt}&=&0\Rightarrow x_2=0 \mbox{ or } \frac{e a(u_1, u_2) x_1}{1+ha(u_1,u_2)x_1}+r_2 \left(1-\frac{x_2}{K_2(u_2)}\right)=d\\
\label{eq-u1}
\frac{du_1}{dt}&=&0\Rightarrow \frac{r_1 x_1 u_1}{\sigma_{K_1}^2K_1(u_1)}=\frac{(u_1-u_2)x_2a(u_1,u_2)}{\sigma_a^2(1+h a(u_1,u_2)x_1)^2},\\
\label{eq-u2}
\frac{du_2}{dt}&=&0\Rightarrow \frac{r_2 x_2 u_2}{\sigma_{K_2}^2K_2(u_2)}=\frac{e(u_1-u_2)x_1a(u_1,u_2)}{\sigma_a^2(1+h a(u_1,u_2)x_1)^2}.
\eae We have the following proposition regarding the boundary equilibria of \eqref{E-pp-g} and their stability.\\
%%%%%%%%%%%%%%%%%%%%%%%%%%%%%%%%%%%
\begin{theorem}\label{th:be}[Equilibria of the evolutionary model \eqref{E-pp-g}] The evolutionary model \eqref{E-pp-g} always has the following two boundary equilibria $E_{0000}=(0,0,0,0),\,\, E_{x_1000}=(K_{01},0,0,0)$ where $E_{x_1000}$ is locally asymptotically stable if $r_2+\frac{ea_0K_{01}}{1+ahK_{01}}<d$ and it is a saddle if $r_2+\frac{ea_0K_{01}}{1+ahK_{01}}>d$. In addition, the following statements are true:
\begin{enumerate}
\item If $r_2>d$, then Model \eqref{E-pp-g} has the third boundary equilibrium  $E_{0x_200}=(0,K_{02}\left(1-\frac{d}{r_2}\right),0,0)$. Both $E_{0000}$ and $E_{0x_200}$ are always saddles. 
\item Model \eqref{E-pp-g} can have one, two, or three interior equilibria $(x^*_1, x^*_2,0,0)$ where $(x^*_1, x^*_2)$ is an interior equilibrium of the ecological model \eqref{pp-g} by letting $K_1=K_{01},\, K_2=K_{02} \mbox{ and } a=a_0$. 
\item The interior equilibrium $(x^*_1, x^*_2,0,0)$ is locally asymptotically stable if the following conditions hold:
\bae\label{stable-evo}\begin{array}{lcl}\min\Big\{\frac{2r_2}{K_{02}},\frac{r_1}{\sigma^2_{K_1}K_{01} }\Big\}>a_0,\,\, \frac{r_1r_2}{K_{01}K_{02}}+\frac{ea_0^2}{(1+a_0h K_{01})^3}>a_0^2hK_{02} (r_2+\frac{e }{h}),\,\, \mbox{ and }\\
\max\Big\{\frac{\sigma_{K_2}}{\sigma_{K_1}}\sqrt{\frac{r_1eK_{02}}{r_2K_{01}}}, \frac{r_1\sigma_{K_2}^2\sigma_a^2}{a_0K_{01}}\big\}>\frac{x^*_2}{x_1^*}\end{array}.\eae 
\end{enumerate}
\end{theorem}

\noindent\textbf{Notes:} According to Corollary \ref{cp2:be-ev}, one necessary condition for $E_{0x_2u_1u_2}$ being locally asymptotically stable is that $a_{u_1u_1}(u_1^*,u^*_2)>0$. The results of Theorem \ref{th:be} imply that the boundary equilibrium $E_{0x_200}=(0,K_{02}\left(1-\frac{d}{r_2}\right),0,0)$ cannot be locally stable, since $a_{u_1u_1}(0,0)<0$. Considering the result in Theorem \ref{p2:be-ev} indicates that the trait function of parasitism efficiency $a(u_1,u_2)$ plays an important role in determining whether the boundary equilibrium $E_{0x_200}$ can be locally stable or not. In addition, Condition \eqref{stable-evo} in Theorem \ref{th:be} indicates that a large ratio of the host (or parasite) intrinsic growth rate $r_1$ to the maximum carrying capacity of host (or parasite ) $K_{01}$,  a large ratio of variance of the carrying capacity of parasite to host, i.e. $\frac{\sigma_{K_2}}{\sigma_{K_1}}$, and large values of the variance of the trait difference of parasitism efficiency $a$, can lead to local stability of the coexistence of host and parasite $(x^*_1, x^*_2,0,0)$.  According to Theorem \ref{th:persistence}, for fixed values of $r_i, K_{0i}, a_0, e$, the proper value of $d$ can guarantee the permanence of the ecological system, and thus guarantee the existence of the interior equilibrium $(x^*_1, x^*_2)$. 
In the rest of this subsection, we will focus on the ESS of the host-only boundary equilibrium $E_{x_1,0,0,0}$ and the interior equilibrium $(x^*_1, x^*_2,0,0)$.\\

%%%%%%%%%%%%%%%%%%%%%%%%%%%%%%%%%%%
\begin{theorem}\label{th:ess-b}[ESS of co-evolutionary model \eqref{E-pp-g}]
The strategy $(u^*_1,u^*_2)=(0,0)$ of the boundary equilibrium $E_{x_1000}$ is an ESS if the following inequality holds
$$r_2+\frac{ea_0K_{01}}{1+a_0hK_{01}}<d.$$
The strategy $(u^*_1,u^*_2)=(0,0)$ of the interior equilibrium $(x^*_1, x^*_2, 0,0)$ is an ESS if Condition \eqref{stable-evo} holds and the following inequalities hold
\bae\label{ESS0}
\begin{array}{lcl}
\frac{r_1 }{\sigma_{K_1}^2K_{01}}>\frac{\max\Big\{\frac{\sigma_{K_2}}{\sigma_{K_1}}\sqrt{\frac{r_1eK_{02}}{r_2K_{01}}}, \frac{r_1\sigma_{K_2}^2\sigma_a^2}{a_0K_{01}}\big\}a_0}{\sigma_a^2}
\end{array}.
\eae
\end{theorem}
\noindent\textbf{Notes:}  By simple rearrangements, we have the following two equivalent equations: 
$$\frac{r_1 }{\sigma_{K_1}^2K_{01}}>\frac{\frac{\sigma_{K_2}}{\sigma_{K_1}}\sqrt{\frac{r_1eK_{02}}{r_2K_{01}}}a_0}{\sigma_a^2}\Leftrightarrow \frac{\sigma_a^2}{\sigma_{K_1}\sigma_{K_2}}>\sqrt{\frac{eK_{01}K_{02}}{r_1r_2}}a_0$$ and 
$$\frac{r_1 }{\sigma_{K_1}^2K_{01}}>\frac{\frac{r_1\sigma_{K_2}^2\sigma_a^2}{a_0K_{01}}a_0}{\sigma_a^2}\Leftrightarrow \sigma_{K_1}\sigma_{K_2}<1.$$
Thus, according to the proof of Theorem \ref{th:be} and \ref{th:ess-b}, we have  the following corollary:
\begin{corollary}\label{c:ess-b}
Assume the following inequalities hold $$\min\Big\{\frac{2r_2}{K_{02}},\frac{r_1}{\sigma^2_{K_1}K_{01} }\Big\}>a_0,\,\, \frac{r_1r_2}{K_{01}K_{02}}+\frac{ea_0^2}{(1+a_0h K_{01})^3}>a_0^2hK_{02} (r_2+\frac{e }{h}).$$
The strategy $(u^*_1,u^*_2)=(0,0)$ of the interior equilibrium $(x^*_1, x^*_2, 0,0)$ is an ESS if one of the following conditions holds
\begin{enumerate}
\item  $\frac{\sigma_{K_2}}{\sigma_{K_1}}\sqrt{\frac{r_1eK_{02}}{r_2K_{01}}}>\frac{x^*_2}{x_1^*}$ and 
$ \frac{\sigma_a^2}{\sigma_{K_1}\sigma_{K_2}}>\sqrt{\frac{eK_{01}K_{02}}{r_1r_2}}a_0.$
\item $\frac{r_1\sigma_{K_2}^2\sigma_a^2}{a_0K_{01}}>\frac{x^*_2}{x_1^*}$ and $\sigma_{K_1}\sigma_{K_2}<1.$

\end{enumerate}
\end{corollary}
\noindent\textbf{Notes:} By comparing the results of Theorem \ref{th:be} with Theorem \ref{th:ess-b} and its corollary \ref{c:ess-b}, we can conclude that the variances of the proposed trait functions, i.e., $\sigma_{K_1},\sigma_{K_2}$ and $\sigma_{d}$, play essential roles in guaranteeing the local stable interior equilibrium  $(x^*_1, x^*_2, 0,0)$ having an ESS. More specifically, large variances of the parasite carrying capacity $\sigma_{K_2}$ and of the parasitism efficiency $\sigma_{a}$ are required to make sure that the strategy $(u^*_1,u^*_2)=(0,0)$ of $(x^*_1, x^*_2, 0,0)$ is an ESS but it may not be global ESS because it is possible that the system has other ESS strategies for coexistence. This can be considered as a \emph{tolerance of the parasite} co-evolutionary outcome, where hosts not only concede to the parasite and accept it in their nests, but also make adjustments to their life history (or other traits) to minimize the negative effects of parasitism on their fitness (Kilner and Langmore 2011). This co-evolutionary outcome is more likely be the case that complete parasitic control of the co-evolutionary trajectory. This could fit with the ecological findings of imperfect recognition of parasite eggs, and a low but positive level of acceptance of brood parasitism by some cuckoo and cowbird host species (Davies \emph{et al} 1996; reviewed by Kruuger 2007).\\

%%%%%%%%%%%%%%%%%%%%%%%%%%%%%%%%%%%
\begin{theorem}\label{th:unique}[The unique ESS of the co-evolutionary model \eqref{E-pp-g}]
The strategy $(u^*_1,u^*_2)=(0,0)$ of the boundary equilibrium $E_{x_1000}$ is the unique ESS of the co-evolutionary host-parasite model \eqref{E-pp-g} if the following inequality holds
$$r_2+\frac{ea_0K_{01}}{1+ahK_{01}}<d.$$
Assume that $\sigma_{K_2}>\sigma_{K_1}$, then the strategy $(u^*_1,u^*_2)=(0,0)$ of the interior equilibrium $(x^*_1, x^*_2, 0,0)$ is the unique ESS of the co-evolutionary host-parasite model \eqref{E-pp-g} whenever the following conditions hold \bae\label{ess-unique}\begin{array}{lcl}
K_{02} a_0<\frac{r_1}{1-\frac{d}{r_2}},\,K_{01}a_0<\frac{1}{h}\\
\min\Big\{\frac{2r_2}{K_{02}},\frac{r_1}{\sigma^2_{K_1}K_{01} }\Big\}>a_0,\,\, \frac{r_1r_2}{K_{01}K_{02}}+\frac{ea_0^2}{(1+a_0h K_{01})^3}>a_0^2hK_{02} (r_2+\frac{e }{h}) \mbox{ and }\\
\frac{\sigma_{K_2}}{\sigma_{K_1}}\sqrt{\frac{r_1eK_{02}}{r_2K_{01}}}>\frac{x^*_2}{x_1^*},\,\frac{\sigma_a^2}{\sigma_{K_1}\sigma_{K_2}}>\sqrt{\frac{eK_{01}K_{02}}{r_1r_2}}a_0.\end{array}\eae
\end{theorem}

\noindent\textbf{Notes:}  Theorem \ref{th:unique} provides sufficient conditions when Model \eqref{E-pp-g} has a unique ESS. One direct implication is that if $(u^*_1,u^*_2)=(0,0)$ is an ESS of the host-only boundary equilibrium $E_{x_1000}$, then Model \eqref{E-pp-g} cannot have other locally asymptotically stable equilibrium. Thus a large death rate due to the parasite hunting/attacking the host can lead to the extinction of the parasite, and lead to a \emph{successful resistance by hosts} evolutionary outcome (Kilner and Langmore 2011).\\\vspace{5pt}

\noindent On the other hand, when the ratio of the variance of the parasite to the host is larger than 1, i.e., $\frac{\sigma_{K_2}}{\sigma_{K_1}}>1$, small values of $a_0$ can lead to the interior equilibrium $(x^*_1, x^*_2, 0,0)$ being the only equilibrium  of the co-evolutionary model \eqref{E-pp-g} having ESS. This type co-evolutionary outcome is referred  to as \emph{acceptance of the parasite}, which can be considered an adaptive strategy for hosts when the costs of rearing a parasite are, on average, lower than any recognition costs (reviewed by Kruuger 2007;  Kilner and Langmore 2011).\\\vspace{5pt}

Our analysis suggests that it is possible for Model \eqref{E-pp-g} to have multiple interior equilibria with potential ESS in the following scenarios:
\begin{enumerate}
\item The condition $\sigma_{K_2}>\sigma_{K_1}$ does not hold in Theorem \ref{th:unique}.
\item The ecological model can have multiple stable equilibria which also can have ESS.
\end{enumerate}

%%%%%%%%%%%%%%%%%%%%%%%%%%%%%%%%%%%%%%%%%%%%%%%%%%%%%%%%%%%%%
%%%%%%%%%%%%%%%%%%%%%%%%%%%%%%%%%%%%%%%%%%%%%%%%%%%%%%%%%%%%%

We have been focusing on the equilibrium $(x^*_1,x^*_2,u^*_1,u^*_2)$ of the co-evolutionary model \eqref{E-pp-g} when $u^*_1=u^*_2=0$. It is possible that our model \eqref{E-pp-g} has an equilibrium with $u^*_1u^*_2\neq 0$. See the following proposition:
%%%%%%%%%%%%%%%%%%%%%%%%%%%%%%%%%%%%%%%%%%%%%%%%%%%%%%%%%%%%%

\begin{proposition}\label{p3:ess} If $(x^*_1,x^*_2,u^*_1,u^*_2)$ is an equilibrium of the co-evolutionary model \eqref{E-pp-g} with $u^*_1u^*_2\neq 0$, then we have $x^*_1x^*_2>0$, and $x^*_1,x^*_2,u^*_1,u^*_2$ can be solved from the following equations:
$$x_1^*=f_1(u^*_1,u^*_2)=\frac{x^*_2}{g(u^*_1,u^*_2)}=f_2(u^*_1,u^*_2)\mbox{ and } \frac{r_1 u_1^*}{\sigma_{K_1}^2K_1(u_1^*)}=\frac{(u_1^*-u_2^*)g(u_1^*,u^*_2)a(u^*_1,u^*_2)}{\sigma_a^2\left(1+h a(u^*_1,u^*_2)f_1(u^*_1,u^*_2)\right)^2}$$where
$$\begin{array}{lcl}
g(u_1,u_2)&=&\sqrt{\frac{e r_1\sigma_{K_2}^2u_1K_2(u_2)}{r_2 \sigma_{K_1}^2u_2K_1(u_1)}}\\
f_1(u_1,u_2)&=&\frac{K_1}{2}-\frac{r_1+a K_1 g}{2r_1a h}+\frac{\sqrt{\left(r_1+a K_1g\right)^2+r_1a h K_1 \left(r_1 a h K_1+2 r_1-2  K_1g\right)}}{2r_1a h}\\
f_2(u_1,u_2)&=&\frac{K_2}{2g}+\frac{a e K_2}{2r_2a hg}-\frac{1}{2a h}+\frac{\sqrt{\left(r_2g-a e K_2\right)^2+r_2a h K_2 \left(r_2 a h K_2+2 r_2g+2 e a K_2\right)}}{2r_1a h}.
\end{array}$$
\end{proposition}
\noindent\textbf{Notes: } Since our trait functions are even function, so are $f_i (u_1, u_2)$ and $g(u_1, u_2)$, i.e.,
$$f_1(u^*_1,u^*_2)=f_1(-u^*_1,-u^*_2), f_2(u^*_1,u^*_2)=f_2(-u^*_1,-u^*_2) \mbox{ and } g(u^*_1,u^*_2)=g(-u^*_1,-u^*_2).$$
Therefore, if $(u^*_1,u^*_2, x^*_1, x^*_2)$ is an interior equilibrium of Model \eqref{E-pp-g}, then so is $(- u^*_1,-u^*_2, x^*_1, x^*_2)$. Proposition \eqref{p3:ess} suggests that the co-evolutionary model \eqref{E-pp-g} could have multiple interior equilibria, and have potential multiple interior ESS, with different trait values. \modifyr{Since we are not able to solve the interior equilibrium $(x^*_1, x^*_2,u^*_1,u^*_2)$ explicitly for the case of $u^*_1u^*_2\neq 0$, we use numerical simulations to illustrate typical scenarios of  the interior equilibrium $(x^*_1, x^*_2,u^*_1,u^*_2)$ when $u^*_1u^*_2\neq 0$. We are particularly interested in the co-evolutionary outcomes when the parasite is facultative, i.e., $r_2>d$. As an example, we fix the values of parameters as follows $$K_{01}=r_1=K_{02}=1;\, a_0=5;\, h=4; \,e=0.9; \,r_2=0.25>d=0.185.$$ We would like to note that the chosen parameter values represent typical dynamics of Model \eqref{E-pp-g} when the ecological model \eqref{pp-g}  processes two interior equilibria, e.g., the blue region of Figure \ref{fig:interior}.} \\

According to Figure \ref{fig:d}, we can see that, under this set of parameter values, the ecological model \eqref{pp-g} has one saddle interior equilibrium; one stable interior equilibrium $(0.5428, 1.0841)$; and one stable boundary equilibrium $(0,0.26)$.  Depending on initial conditions, the trajectory of \eqref{pp-g} will  converge to either $(0.5428, 1.0841)$ or $(0,0.26)$. For example, if we take $(x_1(0), x_2(0))=(0.5,2)$ as the initial conditions, the ecological model \eqref{pp-g} converges to the parasite-only boundary equilibrium $(0,0.26)$. After we turn on evolution, we have the co-evolutionary model \eqref{E-pp-g}. According to Corollary \ref{cp2:be-ev} and Theorem \ref{th:be}, the parasite-only boundary equilibrium $(0,0.26,0,0)$ is unstable. \\

To investigate the co-evolutionary outcome, we fix the initial conditions $(x_1(0), x_2(0), u_1(0),u_2(0))=(0.5,2,1,0.1)$, let $\sigma_{K_1}=\sigma_{K_2}=1$, and vary the value of $\sigma_a$. More specifically, we simulate the dynamics of the co-evolutionary model \eqref{E-pp-g} when $\sigma_a=2,\,1,\,0.5,\,0.15,\,0.05,\,0.005$. The simulations are shown in Figure \ref{fig:sigma2}-\ref{fig:sigma005} where the host $x_1$ is red; the parasite $x_2$ is blue; the trait value for host $u_1$ is green; and the trait value for parasite $u_2$ is cyan: When $\sigma_a=2$, the co-evolutionary model \eqref{Ev-pp} converges to the stable interior equilibrium $(0.0102,    0.3452,   2.3846,    0.1287)$ locally (see Figure \ref{fig:sigma2}); When $\sigma_a=1$, the co-evolutionary model \eqref{Ev-pp} converges to the stable interior equilibrium $(0.0619,    0.4984,    1.8290,    0.4813)$ locally (see Figure \ref{fig:sigma1}); When $\sigma_a=0.5$, the co-evolutionary model \eqref{Ev-pp} converges to the stable  interior equilibrium $(0.5428,    1.0841,   0, 0)$ (see Figure \ref{fig:sigma05}). However, when $\sigma_a=0.15,\,0.05$, the co-evolutionary model \eqref{Ev-pp} has fluctuating dynamics (see Figure \ref{fig:sigma015}-\ref{fig:sigma005}). And when the value of $\sigma_a$ is small enough, e.g., $\sigma_a=0.005$, the co-evolutionary model \eqref{Ev-pp} has fluctuating dynamics first but eventually converges to the stable  interior equilibrium $(0.5428,    1.0841,   0, 0)$ (see Figure \ref{fig:sigma0005}). These simulations suggest follows:

 \begin{figure}[ht]
\centering
  \includegraphics[scale =.45] {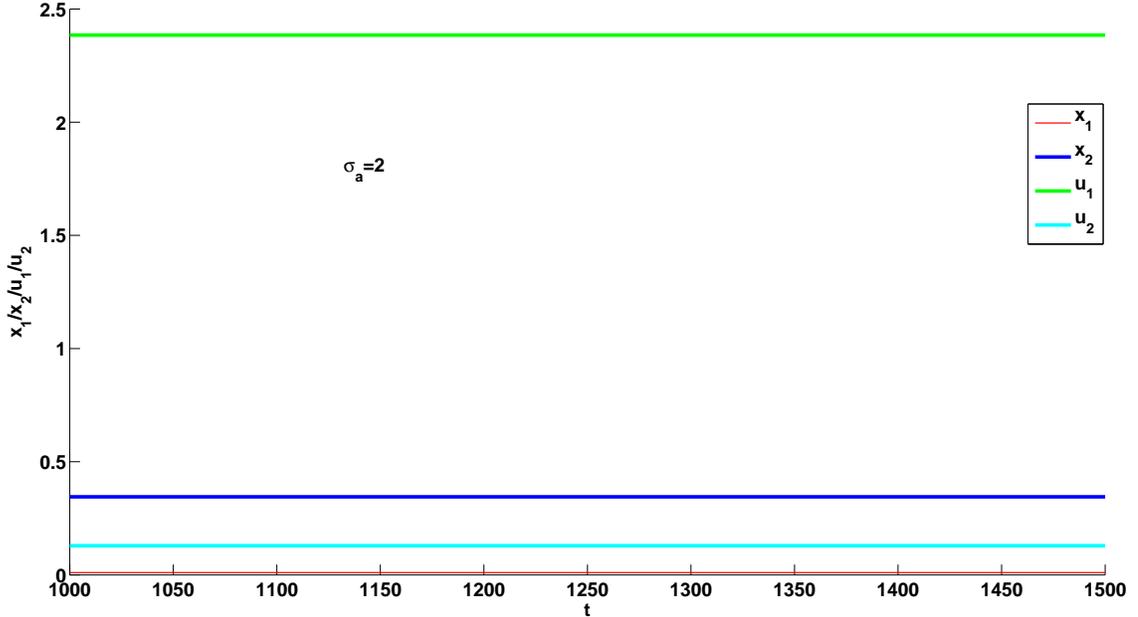}
 \caption{\modifyr{The stable interior equilibrium is $(0.0102,    0.3452,   2.3846,    0.1287)$}. The dynamics of the co-evolutionary model \eqref{Ev-pp} when $\sigma_{K_1}=\sigma_{K_2}=K_{01}=r_1=K_{02}=1,\,\,r_2=0.25>d=0.185,\,\,a_0=5,\,\,e=0.9,\,\,h=4,\,\,\sigma_a=2$ with initial values $(x_1(0), x_2(0), u_1(0),u_2(0))=(0.5,2,1,0.1)$: the host $x_1$ is red; the parasite $x_2$ is blue; the trait value for the host $u_1$ is green; and the trait value for the parasite $u_2$ is cyan.} \label{fig:sigma2}
\end{figure}

 \begin{figure}[ht]
\centering
  \includegraphics[scale =.45] {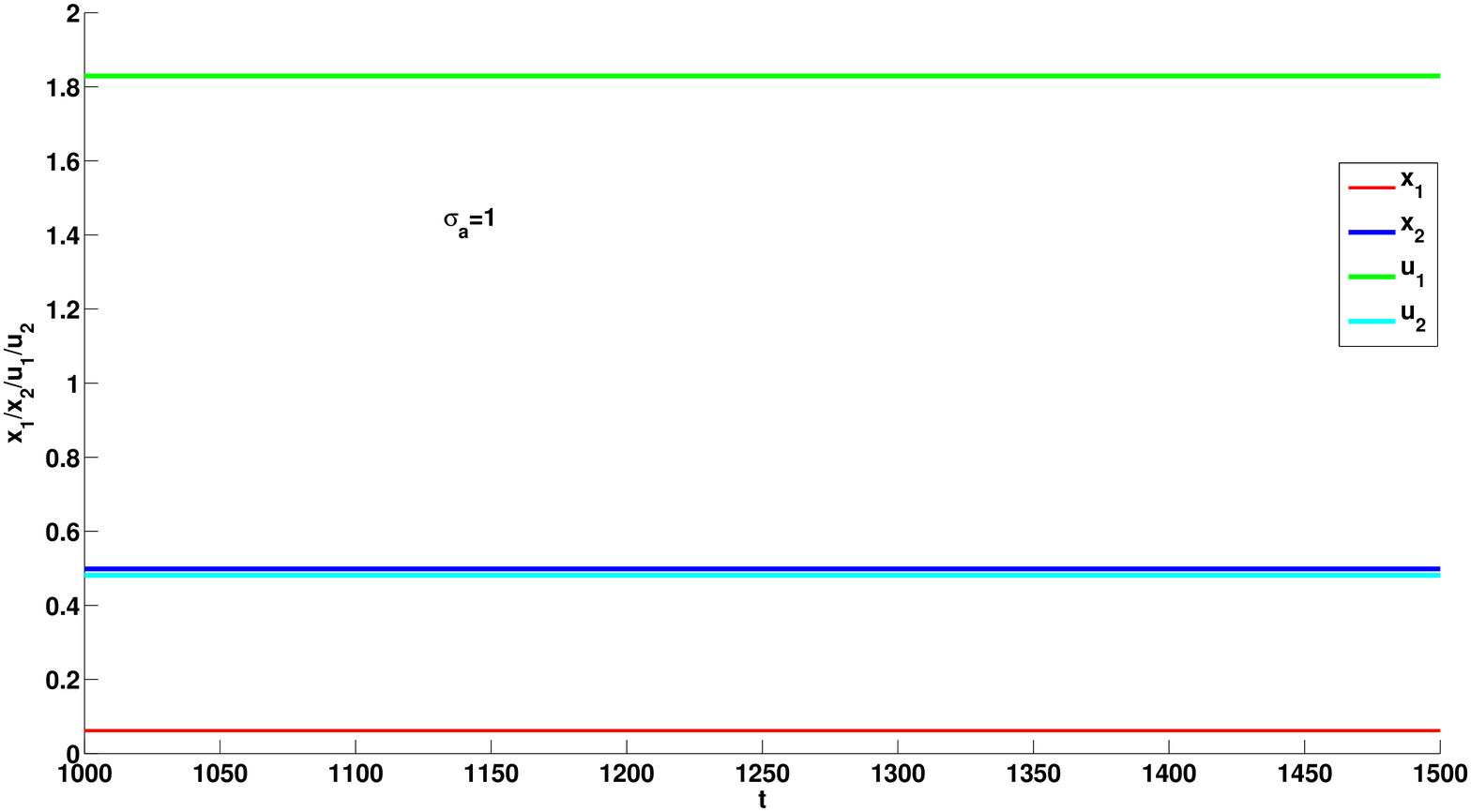}
 \caption{\modifyr{The stable interior equilibrium is  $(0.0619,    0.4984,    1.8290,    0.4813)$}. The dynamics of the co-evolutionary model \eqref{Ev-pp} when $\sigma_{K_1}=\sigma_{K_2}=K_{01}=r_1=K_{02}=1,\,\,r_2=0.25>d=0.185,\,\,a_0=5,\,\,e=0.9,\,\,h=4,\,\,\sigma_a=1$ with initial values $(x_1(0), x_2(0), u_1(0),u_2(0))=(0.5,2,1,0.1)$: the host $x_1$ is red; the parasite $x_2$ is blue; the trait value for host $u_1$ is green; and the trait value for the parasite $u_2$ is cyan.} \label{fig:sigma1}
\end{figure}

 \begin{figure}[ht]
\centering
  \includegraphics[scale =.45] {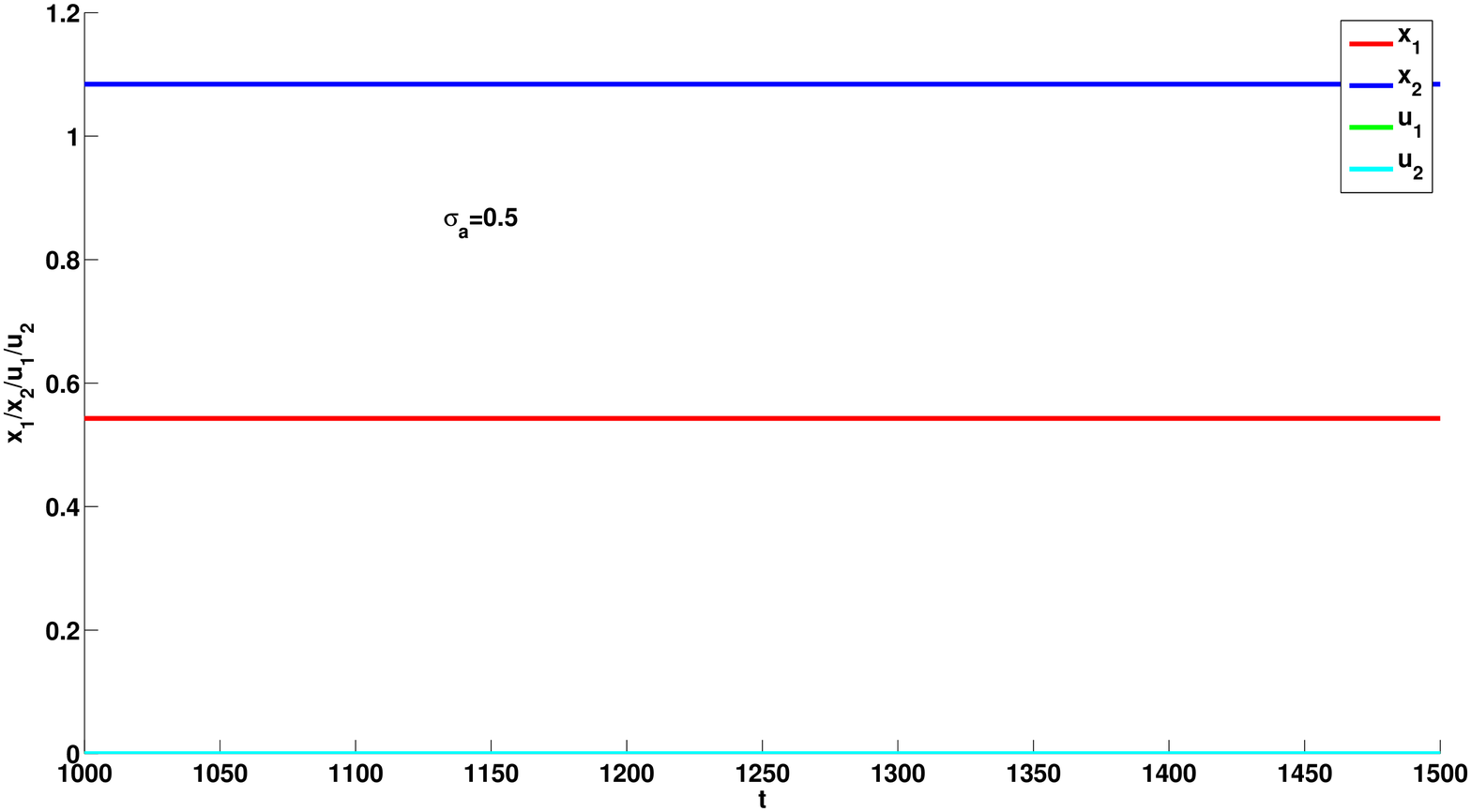}
 \caption{\modifyr{The stable interior equilibrium is  $(0.5428,    1.0841,   0, 0)$}. The dynamics of the co-evolutionary model \eqref{Ev-pp} when $\sigma_{K_1}=\sigma_{K_2}=K_{01}=r_1=K_{02}=1,\,\,r_2=0.25>d=0.185,\,\,a_0=5,\,\,e=0.9,\,\,h=4,\,\,\sigma_a=0.5$ with initial values $(x_1(0), x_2(0), u_1(0),u_2(0))=(0.5,2,1,0.1)$: the host $x_1$ is red; the parasite $x_2$ is blue; the trait value for the host $u_1$ is green; and the trait value for the parasite $u_2$ is cyan.} \label{fig:sigma05}
\end{figure}

 \begin{figure}[ht]
\centering
  \includegraphics[scale =.45] {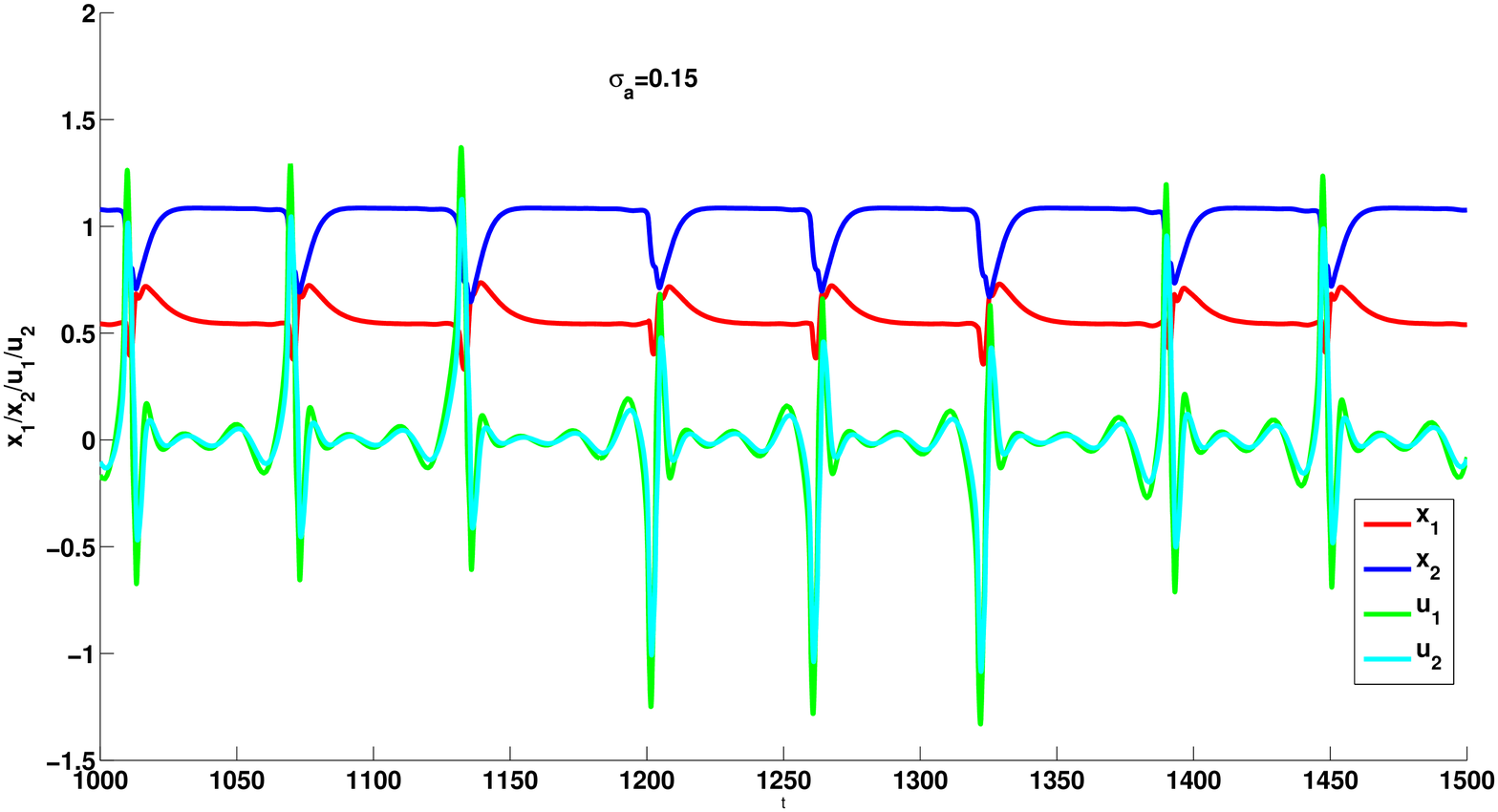}
 \caption{\modifyr{The oscillating dynamics when $\sigma_a$ is small}. The dynamics of the co-evolutionary model \eqref{Ev-pp} when $\sigma_{K_1}=\sigma_{K_2}=K_{01}=r_1=K_{02}=1,\,\,r_2=0.25>d=0.185,\,\,a_0=5,\,\,e=0.9,\,\,h=4,\,\,\sigma_a=0.15$ with initial values $(x_1(0), x_2(0), u_1(0),u_2(0))=(0.5,2,1,0.1)$: the host $x_1$ is red; the parasite $x_2$ is blue; the trait value for the host $u_1$ is green; and the trait value for the parasite $u_2$ is cyan.} \label{fig:sigma015}
\end{figure}

\begin{figure}[ht]
\centering
  \includegraphics[scale =.45] {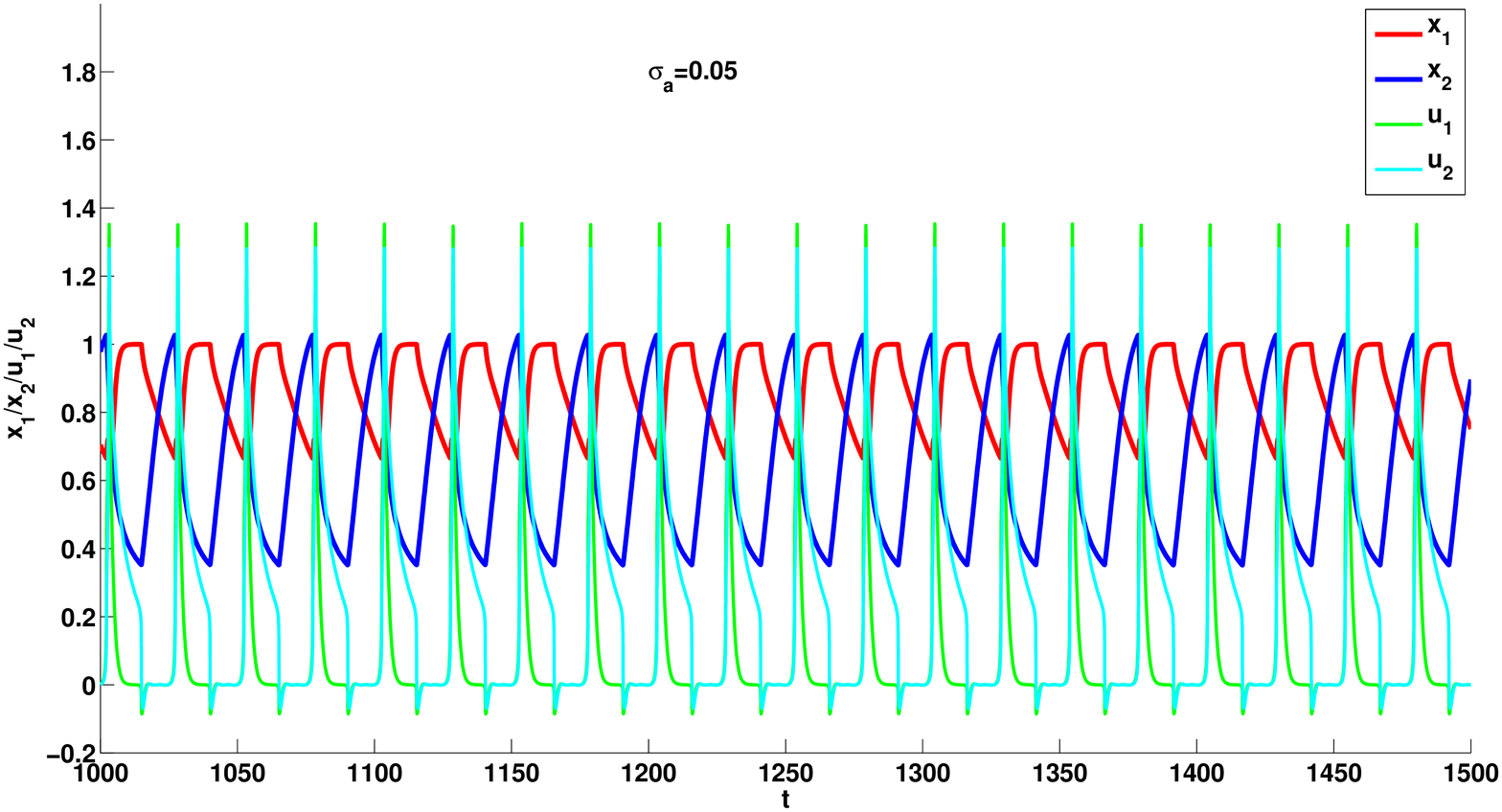}
 \caption{\modifyr{The larger oscillating dynamics when $\sigma_a$ is smaller}. The dynamics of the co-evolutionary model \eqref{Ev-pp} when $\sigma_{K_1}=\sigma_{K_2}=K_{01}=r_1=K_{02}=1,\,\,r_2=0.25>d=0.185,\,\,a_0=5,\,\,e=0.9,\,\,h=4,\,\,\sigma_a=0.05$ with initial values $(x_1(0), x_2(0), u_1(0),u_2(0))=(0.5,2,1,0.1)$: the host $x_1$ is red; the parasite $x_2$ is blue; the trait value for host $u_1$ is green; and the trait value for parasite $u_2$ is cyan.} \label{fig:sigma005}
\end{figure}
\begin{figure}[ht]
\centering
  \includegraphics[scale =.45] {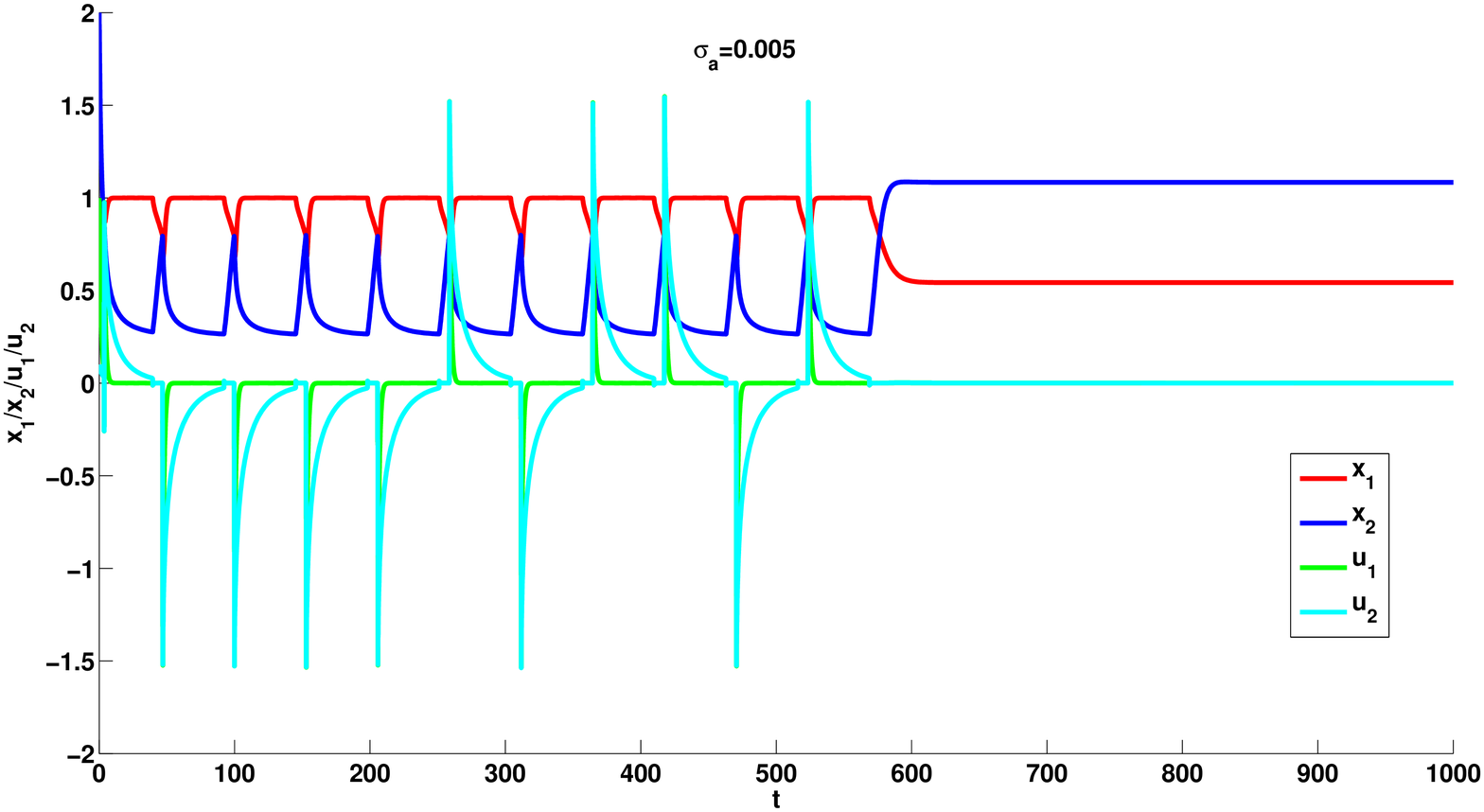}
 \caption{\modifyr{The oscillating dynamics converge to the stable interior equilibrium is $(0.5428,    1.0841,   0, 0)$  when $\sigma_a=0.005$}. The dynamics of the co-evolutionary model \eqref{Ev-pp} when $\sigma_{K_1}=\sigma_{K_2}=K_{01}=r_1=K_{02}=1,\,\,r_2=0.25>d=0.185,\,\,a_0=5,\,\,e=0.9,\,\,h=4,\,\,\sigma_a=0.05$ with initial value $(x_1(0), x_2(0), u_1(0),u_2(0))=(0.5,2,1,0.1)$: the host $x_1$ is red; the parasite $x_2$ is blue; the trait value for host $u_1$ is green; and the trait value for parasite $u_2$ is cyan.} \label{fig:sigma0005}
\end{figure}

\begin{enumerate}
\item \textbf{The stable interior equilibrium $(x^*_1, x^*_2,u^*_1,u^*_2)$ with $u^*_1u^*_2\neq 0$:} The co-evolutionary model \eqref{E-pp-g} can indeed have locally stable interior equilibrium $(x^*_1, x^*_2,u^*_1,u^*_2)$ with $u^*_1u^*_2\neq 0$ (see Figure \ref{fig:sigma2} and \ref{fig:sigma1}). This occurs when the ecological dynamics \eqref{pp-g} has two attractors in the absence of evolution: the parasite-only equilibrium $\left(0,K_{02}\left(1-\frac{d}{r_2}\right)\right)$ and the interior attractor where both the host and parasite coexist.

\item \textbf{ESS or not:} \modifyr{Calculations show that: (a) when $\sigma_a=2,\,1$, the co-evolutionary model \eqref{Ev-pp} converges to the stable interior equilibrium $(0.0102,    0.3452,   2.3846,    0.1287)$, $(0.0619,    0.4984,    1.8290,    0.4813)$ locally, respectively, which are not ESS; and (b) when $\sigma_a=0.5,\,0.005$, the co-evolutionary model \eqref{Ev-pp} converges to the stable  interior equilibrium $(0.5428,    1.0841,   0, 0)$ where  the strategy $(0,0)$ is an ESS.}

\item \textbf{Evolution can save the host from local extinction by creating new stable interior equilibria (see Figure \ref{fig:sigma2}-\ref{fig:sigma1}) or shifting the basins of attractions (see Figure \ref{fig:sigma05})}: In the absence of evolution, under the condition of $K_{01}=r_1=K_{02}=1;\, a_0=5;\, h=4; \,e=0.9; \,r_2=0.25>d=0.185$ and $(x_1(0), x_2(0))=(0.5,2)$, the ecological model \eqref{pp-g} converges to the parasite-only equilibrium $(0,0.26)$; However, in the presence of evolution,  the co-evolutionary model \eqref{Ev-pp} can converge to the stable interior equilibrium locally where both the host and parasite can coexist.  

\item \textbf{Effects of the variance $\sigma_a$ of the trait difference in the parasitism efficiency $a(u_1,u_2)$}: Small values of $\sigma_a$ can destabilize the evolutionary dynamics and generate oscillations between population and trait dynamics (see Figure \ref{fig:sigma015}-\ref{fig:sigma005}), while extremely small values of $\sigma_a$ can cause \emph{catastrophe events} where the large oscillations collapse with the boundary of the basins attractions of the stable interior equilibrium $(x^*_1, x^*_2,0,0)$ with the consequence that the trajectory converges to $(x^*_1, x^*_2,0,0)$ (see Figure \ref{fig:sigma0005}).\\

\item \textbf{Host-parasite cycles when $\sigma_a$ is small}: Compare the direction of host-parasite cycles in Figure \ref{fig:sigma015HP} to Figure \ref{fig:sigma005HP}: For $\sigma_a=0.15$, both host and parasite exhibit their peaks at almost the same time, but when $\sigma_a=0.05$, the host peaks at the lowest population of the parasite and the parasite peaks at the lowest population of the host. These simulation results suggest that the ability of  the host and parasite to adapt, as measured by the variance of the trait differences in the parasitism efficiency $a$, can alter the community dynamics of natural systems, leading to novel dynamics including antiphase and cryptic cycles. Such cycles, induced by the variance of the  trait difference measured in the parasitism efficiency $a$, could be another potential signature of host-parasite co-evolution and reveal that host-parasite co-evolution can shape, and possibly reverse, community dynamics. These results offer another perspective on the recent work by Cortez and Weitza (2014), which uses an eco-co-evolutionary model to show that predator-prey co-evolution can drive population cycles where the opposite of canonical Lotka-Volterra oscillations occurs; predator peaks precede prey peaks. \\
 \begin{figure}[ht]
\centering
 \subfigure[The oscillations  when $\sigma_a=0.15$]{
  \includegraphics[scale =.45] {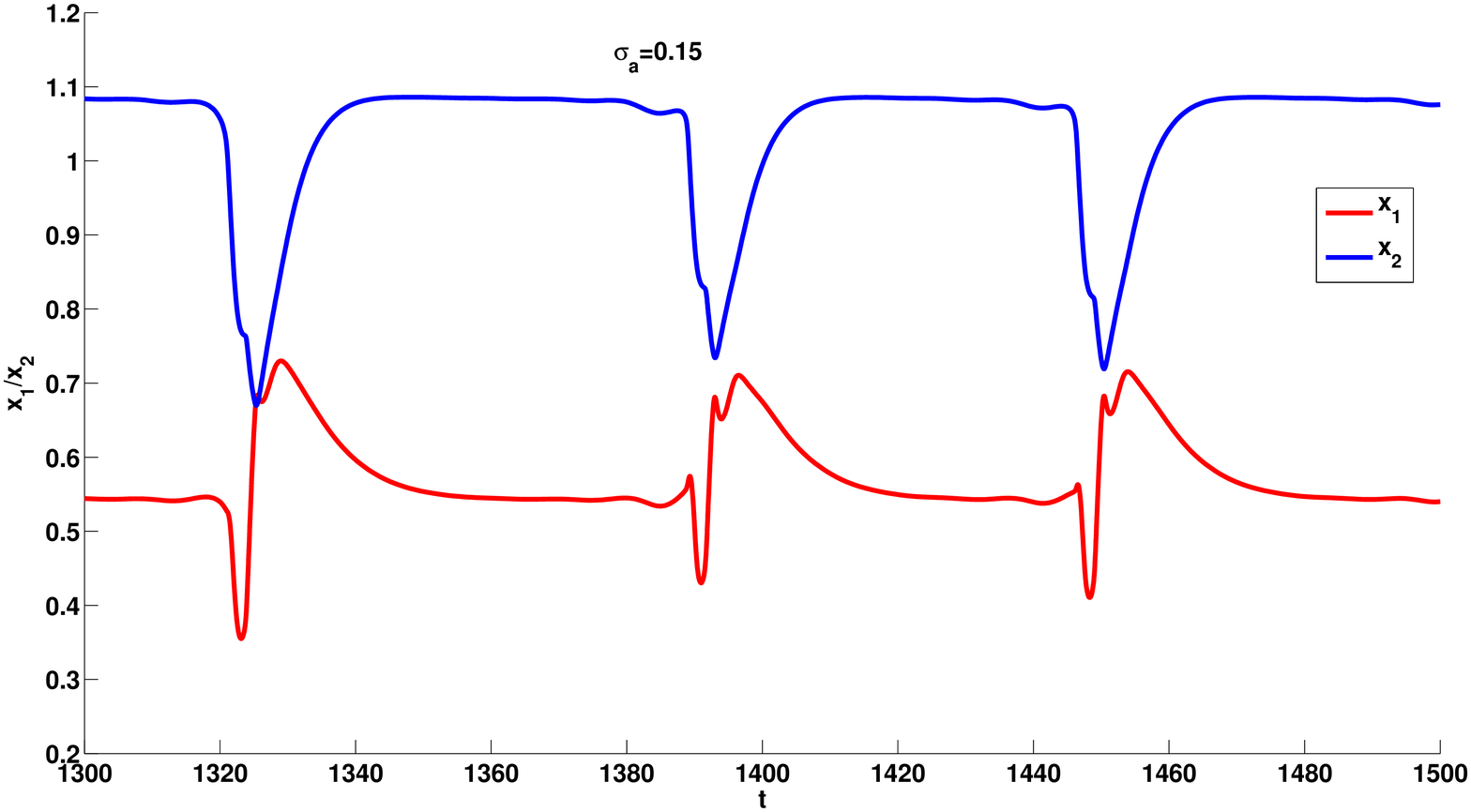}
 \label{fig:sigma015HP}}
   \subfigure[The larger and faster oscillations  when $\sigma_a=0.05$]{
   \includegraphics[scale =.45] {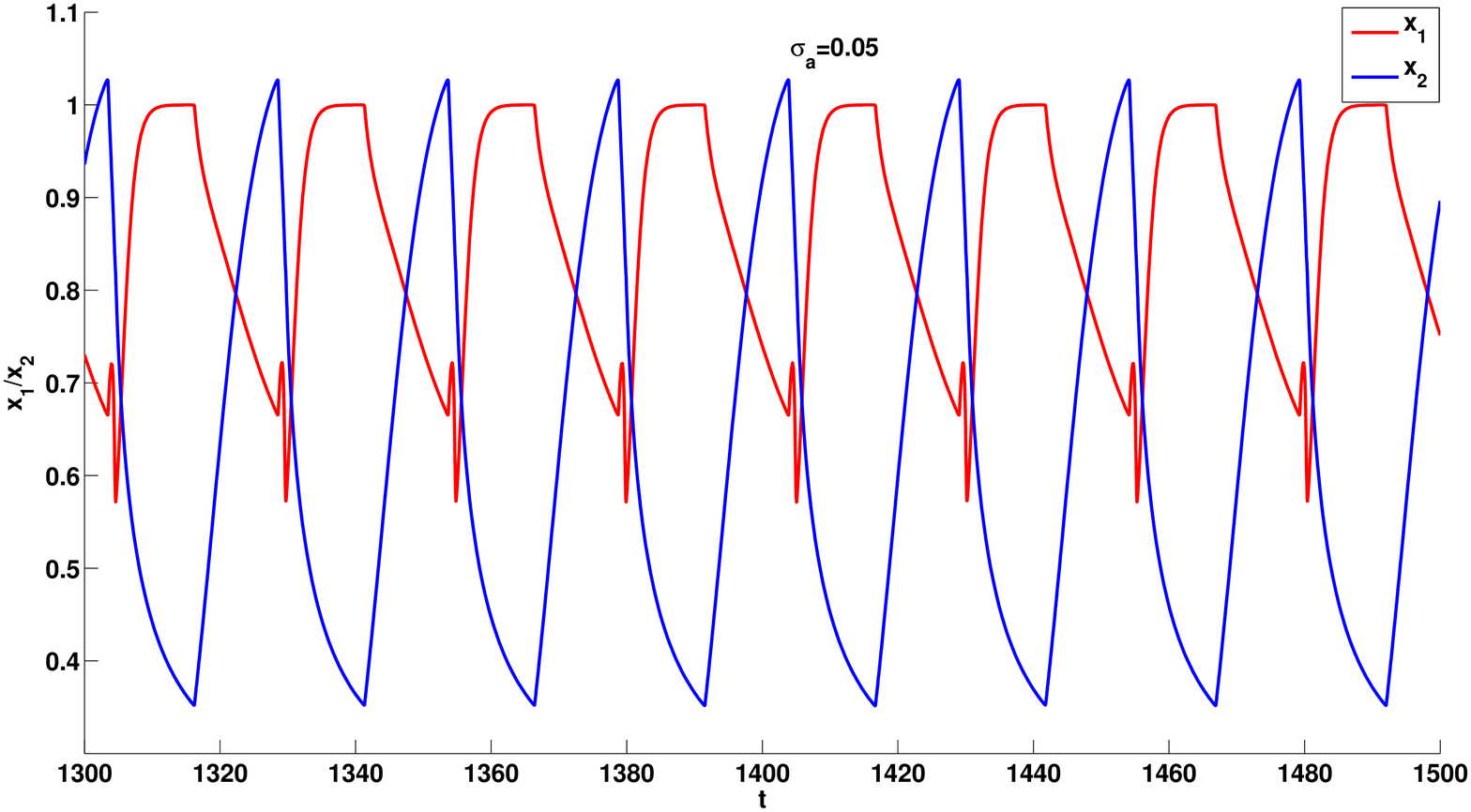}
 \label{fig:sigma005HP}}
  \caption{The population dynamics of the host and parasite for the co-evolutionary model \eqref{Ev-pp} when $\sigma_{K_1}=\sigma_{K_2}=K_{01}=r_1=K_{02}=1,\,\,r_2=0.25>d=0.185,\,\,a_0=5,\,\,e=0.9,\,\,h=4$ with initial values $(x_1(0), x_2(0), u_1(0),u_2(0))=(0.5,2,1,0.1)$, while the value of $\sigma_a$ is 0.15 (Figure \ref{fig:sigma015HP}) and 0.05 (Figure \ref{fig:sigma005HP}): the host $x_1$ is red; and the parasite $x_2$ is blue}\label{fig:evoHP}
\end{figure}

\end{enumerate}
%%%%%%%%%%%%%%%%%%%%%%%%%%%%%%%%%%%%%%%%%%%%%%%%%%%%%%%%%%%%%
%%%%%%%%%%%%%%%%%%%%%%%%%%%%%%%%%%%%%%%%%%%%%%%%%%%%%%%%%%%%%
%\subsection{Comparisons}
%%%%%%%%%%%%%%%%%%%%%%%%%%%%%%%%%%%%%%%%%%%%%%%%%%%%%%%%%%%%%
%\subsection{Comparison of co-evolutionary dynamics when parasite is facultative v.s. obligate}
%%%%%%%%%%%%%%%%%%%%%%%%%%%%%%%%%%%%%%%%%%%%%%%%%%%%%%%%%%%%%
%%%%%%%%%%%%%%%%%%%%%%%%%%%%%%%%%%%%%%%%%%%%%%%%%%%%%%%%%%%%%

\section{Conclusion}

Host-parasite co-evolution, or the reciprocal evolution of host defenses and parasite counter-strategies, can affect a range of ecological and evolutionary processes (Thompson 2005; G{\'o}mez {\emph{et al}.} 2014), from population dynamics (Yoshida {\emph{et al}.} 2003), to the maintenance of genetic variation (Clark {\emph{et al}} 2007).  Despite this, there are relatively few mathematical models examining the co-evolution of quantitative traits in hosts and parasites  (Best {\emph{et al}} 2009). Although most models  assume that the parasite is obligate, cases of facultative and/or generalist parasites are common and ecologically significant (Spottiswoode {\emph{et al}} 2012). \\

In this paper, we have presented a co-evolutionary model of a social parasite-host system that includes (1) ecological dynamics that feed back into the co-evolutionary outcomes;  (2) consideration of both obligate and facultative parasitic strategies; and (3) Holling Type II functional responses between the host and parasite. The analytical study on the proposed model provides insightful information on the impacts of host-parasite co-evolution on varied ecological and evolutionary processes. Here we summarize our main results as follows.\\

Recall that $d$ is the death rate of a parasite due to searching/hunting for all potential hosts, and $r_2\geq 0$ is the intrinsic growth rate of the parasite without parasitizing a given  host $x_2$. If $d>r_2$, the parasite is obligate; and if $d<r_2$, the parasite is facultative.  In the absence of evolution, we performed local and global analyses to investigate the ecological outcomes when parasites range from facultative to obligate.\\

 When we fix other parameters and let $d$ vary, our proposed ecological model can exhibit a wide range of dynamics illustrating how ecological dynamics change when the parasite makes the transition from facultative to obligate. The typical dynamics are shown Figure \ref{fig:interior}-\ref{fig:d}. As an example displayed in Figure \ref{fig:d}, we can see that: when $d$ is extremely small, the host goes extinct globally; when $d$ is small, the system can exhibit bi-stability between the parasite-only boundary equilibrium and the coexistence interior attractor; when $d$ is in the intermediate range of values, the system is permanent and can process three interior equilibria that lead to two distinct coexistence attractors; when $d$ is large, the system is permanent with only one coexistence attractor; however, when $d$ is extremely large, the parasite goes extinct, resulting in global stability of the host-only boundary equilibrium. More specifically, our analytical results imply that: \\%Our results show that  Figure \ref{fig:bifur}
\begin{enumerate}
\item Facultative parasitism can drive the host extinct locally (Proposition \ref{p1:be}) or globally (Theorem \ref{th:global}), while obligate parasitism can only generate global extinction of parasites (Proposition \ref{p1:be} and Theorem \ref{th:global}).\\
\item \modifyr{When the parasite is obligate, i.e., $d>r_2$, the host always persists. However, when the parasite is facultative, i.e., $d<r_2$, the host can go to extinction under certain conditions, while the facultative parasite always persists} (see Theorem \ref{th:persistence}).\\
\item The host-facultative parasite model can have rich dynamics that process one, two, and three interior equilibria with consequences of two or three attractors; the host-obligate parasite model can have either one or three interior equilibria (Theorem \ref{th:interior}). When the system has two interior equilibria (only for the facultative parasite) the system has a parasite-only boundary attractor and a coexistence interior attractor; when the system has three interior equilibria the system has two distinct coexistence interior attractors.\\

%\item 
\end{enumerate}

Host-parasite co-evolution plays an essential role in both ecological and evolutionary processes. Our work on the co-evolutionary dynamics confirms this, and addresses the importance of trait function effects on evolutionary outcomes. More specifically, our main findings are:\\
\begin{enumerate}
 \item When the death rate $d$ of the parasite depends on its trait value, under proper trait functions (Theorem \ref{p2:be-ev}), the parasite can choose different strategies such that it can be facultative at some trait values while it is obligate at others. However, these strategies are not ESS (see Figure \ref{fig1:G1G2}-\ref{fig2:G1G2}). A potential biological example supporting these results may be slave-making ants in the \emph{Formica sanguinea} complex. These species display the ability to vary behaviorally between reliance on parasitism of other ants, and producing their own worker offspring. They have been suggested to represent an intermediate parasitic stage, between those species without predominant brood raiding, and obligate slave-making species (Topoff and Zimmerli 1991; Ruano \emph{\emph{et al}.} 2013). \\
 
\item  When the death rate $d$ of the parasite is independent of its trait value, Corollary \ref{cp2:be-ev}  implies that the host-only equilibrium $E_{x_10u_1u_2}$ can have ESS when the parasite is obligate; and the facultative parasite-only equilibrium $E_{0x_2u_1u_2}$ can also have ESS.\\

 \item Let the death rate $d$ of the parasite be independent of its trait value, and other trait functions follow normal distributions. Our results (Theorem \ref{th:ess-b} and its corollary \ref{c:ess-b}; Theorem \ref{th:unique}) show that: (1) co-evolution can save the host from extinction by destabilizing the facultative-parasite only boundary equilibrium and generating a new locally stable coexistence interior equilibrium (see Theorem \ref{th:be} and Figure \ref{fig:sigma2}-\ref{fig:sigma0005});  (2) the variances of the proposed trait functions, i.e., $\sigma_{K_1},\sigma_{K_2}$ and $\sigma_{d}$, play essential roles in guaranteeing the local stable interior equilibrium  $(x^*_1, x^*_2, 0,0)$ having ESS; (3) when the parasite death rate $d$ is large, the host can have ESS that drive the parasite extinct globally; and (4) Large variances of the parasite carrying capacity $\sigma_{K_2}$ and of the parasitism efficiency $\sigma_{a}$ are required to make sure of $(x^*_1, x^*_2, 0,0)$ being ESS. In addition, we have following interesting results:\\
 \begin{enumerate}
 \item  It is possible to have locally stable coexistence interior equilibrium $(x^*_1, x^*_2, u^*_1, u^*_2)$ with $u^*_1, u^*_2>0$ (see Proposition \eqref{p3:ess} and Figure \ref{fig:sigma2}-\ref{fig:sigma1}) but these strategies may not be ESS. This can occur when the ecological system has two interior equilibria with the facultative parasite-only equilibrium being locally stable; and co-evolution destabilizes this facultative-parasite-only boundary equilibrium and generates a new locally stable coexistence interior equilibrium .\\
 %\item Figure \ref{fig:evopp} - \ref{fig:evoHP}
 \item   The variance of the trait difference in parasitism efficiency has huge impacts on the population dynamics. More specifically, small values of the variance $\sigma_a$ can destabilize the system, thus generate evolutionary arms-race dynamics with different host-parasite fluctuating patterns (see Figure \ref{fig:sigma2} - \ref{fig:evoHP}).\\
\end{enumerate}
\end{enumerate}

\modifyr{Our theoretical work on the ecological and co-evolutionary dynamics of the host-parasite system show interesting matches with recent ecological considerations of social parasitism (Foitzik \emph{et al.} 2001 \& 2003; reviewed by Kruger 2007; Davies 2011; Kilner and Langmore 2011) in the following ways:  
%greatly support the recent work by Kilner and Langmore (2011):
\begin{enumerate}
\item The study in Section 4.1 of the host-only equilibrium suggests that  the host-only equilibrium $E_{x_10u_1u_2}$ of the co-evolutionary model \eqref{Ev-pp} can have two ESS  $(u^*_1,u^*_2)=(0,c)$ and $(u^*_1,u^*_2)=(c,0)$ when the parasite is obligate (i.e., $d>r_2$). This can be classified as one of the co-evolutionary outcomes when the host successfully resists invasion by the parasite, and when resistance is due to effective front-line defenses (Foitzik \emph{et al.} 2001 \& 2003; reviewed by Kruger 2007; Davies 2011; Kilner and Langmore 2011). 
\item The results in Theorem \ref{th:be}, Theorem \ref{th:ess-b} and its corollary \ref{c:ess-b} imply that large variances of the parasite carrying capacity $\sigma_{K_2}$ and of parasitism efficiency $\sigma_{a}$ are required to make sure of the strategy $(u^*_1,u^*_2)=(0,0)$ of $(x^*_1, x^*_2, 0,0)$ being an ESS but it may not be a global ESS, because of the possibility that the system has other ESS for coexistence. This can be considered as equivalent to a \emph{tolerance of the parasite} co-evolutionary outcome where hosts to some degree concede to the parasite, and accept it into their nests, but also make adjustments to their life history (or other traits) to minimize the negative effects of parasitism on their fitness (Foitzik \emph{et al.} 2001 \& 2003; reviewed by Kruger 2007; Davies 2011; Kilner and Langmore 2011). This co-evolutionary outcome is more likely to be the case that parasitic control of the co-evolutionary trajectory. 
\item Theorem \ref{th:unique}: A large death rate due to the parasite overhunting/attacking the host can lead to the extinction of the parasite, and lead to a \emph{successful resistance by hosts} evolutionary outcome (Foitzik \emph{et al.} 2001 \& 2003; reviewed by Kruger 2007; Davies 2011; Kilner and Langmore 2011). 
\item Theorem \ref{th:unique}: when the ratio of the variance of the parasite to the host is larger than 1, i.e., $\frac{\sigma_{K_2}}{\sigma_{K_1}}>1$, small values of $a_0$ can lead to the interior equilibrium $(x^*_1, x^*_2, 0,0)$ being the only equilibrium that has ESS in the co-evolutionary model \eqref{E-pp-g}. This type co-evolutionary outcome is referred to as \emph{acceptance of the parasite}. It can become an adaptive strategy for hosts when the costs of rearing a parasite are, on average, lower than the costs of recognition and rejection (Foitzik \emph{et al.} 2001 \& 2003; reviewed by Kruger 2007; Davies 2011; Kilner and Langmore 2011).
\end{enumerate}}
%The recent work by Kilner and Langmore (2011)

%%%%%%%%%%%%%%%%%%%%%%%%%%%%%%%%%%%%%%%%%%%%%%%%%%%%%%%%%%%%%%%%
%%%%%%%%%%%%%%%%%%%%%%%%%%%%%%%%%%%%%%%%%%%%%%%%%%%%%%%%%%%%%%%%
\section{Proofs}
\subsection*{Proof of Lemma \ref{l1:dp}}
\begin{proof} Both sets $X_1$ and $X_2$ are invariant for the ecological model \eqref{pp-g} since we have 
$$\frac{dx_1}{dt}\big\vert_{x_1=0}= \frac{dx_2}{dt}\big\vert_{x_2=0} =0.$$
The continuity of the expressions of both $\frac{dx_1}{dt}$ and $\frac{dx_2}{dt}$ implies that Model \eqref{pp-g} is positively invariant in $X$.

Let $K^M_1=\sup_{u_1\in \mathbb R}\{K_1(u_1)\}$ and $K^M_2=\sup_{u_2\in \mathbb R}\{K_2(u_2)\}$. From  \eqref{pp-g} and its positively invariant property in $X$, we have the following two inequalities
$$\begin{array}{lcl}
\frac{dx_1}{dt}&=&x_1 \left[r_1 (1-\frac{x_1}{K_1})-\frac{a x_2}{1+h a x_1}\right]\leq  r_1 x_1 (1-\frac{x_1}{K_1})\leq  r_1 x_1 (1-\frac{x_1}{K_{1}^M}) \\\\
\frac{dx_2}{dt}&=&x_2\left[\frac{e a x_1}{1+ha x_1}-d+r_2 (1-\frac{x_2}{K_2})\right]\leq \frac{e x_2}{h}+ r_2x_2 (1-\frac{x_2}{K_2^M}) =(r_2+\frac{e }{h})x_2\left(1-\frac{x_2}{\frac{K_2^M (r_2+\frac{e }{h})}{r_2}}\right)\\\\
&\leq& (r_2+\frac{e }{h})x_2\left(1-\frac{x_2}{\frac{K_{2}^M (r_2+\frac{e }{h})}{r_2}}\right)
\end{array}
$$ which implies that 
$$\limsup_{t\rightarrow\infty} x_1(t)  \leq K_{1}^M \mbox{ and } \limsup_{t\rightarrow\infty} x_2(t)  \leq \frac{K_{2}^M (r_2+\frac{e }{h})}{r_2}.$$ 
Therefore, the statement of Lemma \ref{l1:dp} holds.
\end{proof}
%%%%%%%%%%%%%%%%%%%%%%%%%%%%%%%%%%%%%%%%%%%%%%%%%%%%%%%%%%%%%
\subsection*{Proof of Proposition \ref{p1:be}}
\begin{proof} The Jacobian matrix of the ecological model \eqref{pp-g} evaluated at its equilibrium $(x_1^*,x_2^*)$ can be represented as follows:
\bae\label{J}
J_{(x_1^*,x_2^*)}=\left[\begin{array}{cc}
r_1 (1-\frac{x_1^*}{K_1})-\frac{a x_2^*}{1+h a x_1^*}+x_1^*\left(-\frac{r_1x_1^*}{K_1}+\frac{e a x_2^*}{(1+a hx_1^*)^2}\right)&-\frac{a x_1^*}{1+h a x_1^*}\\
\frac{ea x^*_2}{(1+a hx_1^*)^2}&r_2(1-\frac{x_2^*}{K_2}+\frac{ea x_2^*}{1+h a x_1^*})-d-\frac{r_2x_2^*}{K_2}
\end{array}\right].
\eae
Therefore, we can have the following three cases:
\begin{enumerate}
\item The eigenvalues of the extinction equilibrium $E_{00}$ are $\lambda_1=r_1$ and $\lambda_2=r_2-d$. Thus, $E_{00}$ is a source if $r_2>d$ while it is a saddle if $r_2<d$.
\item The eigenvalues of $E_{10}$ are $\lambda_1=-r_1$ and $\lambda_2=r_2-d+\frac{eaK_1}{1+ahK_1}$. Thus, $E_{10}$ is a sink if $r_2+\frac{eaK_1}{1+ahK_1}<d$ while it is a saddle if $r_2+\frac{eaK_1}{1+ahK_1}>d$.
\item The equilibrium $E_{01}$ exists if $r_2>d$. In this case, $E_{00}$ is a source and $E_{10}$ is a saddle. The eigenvalues of $E_{01}$ are $\lambda_1=r_1+\frac{ad K_2}{r_2}-aK_2$ and $\lambda_2=d-r_2<0$. Therefore, $E_{01}$ is a sink if $r_1+\frac{ad K_2}{r_2}<aK_2$ (i.e., $\frac{r_1}{aK_1}<1-\frac{d}{r_2}$) while it is a saddle if $r_1+\frac{ad K_2}{r_2}>aK_2$ (i.e., $\frac{r_1}{aK_1}>1-\frac{d}{r_2}$).
\end{enumerate}
Thus, the statement of Proposition \ref{p1:be} holds.
%\modifyr{will insert here soon by looking the Jacobian matrix evaluated these equilibria}

\end{proof}
%%%%%%%%%%%%%%%%%%%%%%%%%%%%%%%%%%%%%%%%%%%%%%%%%%%%%%%%%%%%%
\subsection*{Proof of  Theorem \ref{th:persistence}}
\begin{proof} 
According to Lemma \ref{l1:dp} and Proposition \ref{p1:be}, we can conclude that (i) Model \eqref{pp-g} is positively invariant in both $X_1$ and $X_2$; (ii) the omega limit set of $X_1$ is $E_{10}$; and (iii) the omega limit set of $X_2$ is $E_{00}$ if $r_2<d$ while the omega limit set of $X_2$ is $E_{01}$ if $r_2>d$.

The results (Theorem 2.5 and its corollary) in Hutson (1984) guarantee that the persistence of the host $x_1$ is determined by the sign of
$$\frac{dx_1}{x_1dt}\Big\vert_{X_1}=\left\{\begin{array}{ll}\frac{dx_1}{x_1dt}\Big\vert_{E_{00}} \mbox{ if } r_2<d\\\\
\frac{dx_1}{x_1dt}\Big\vert_{E_{01}} \mbox{ if } r_2>d\end{array}\right\}.$$ And the persistence of the parasite $x_2$ is determined by the sign of
$$\frac{dx_2}{x_2dt}\Big\vert_{X_1}=\frac{dx_2}{x_2dt}\Big\vert_{E_{10}}.$$
Therefore, we can conclude the following:
\begin{enumerate}
\item The prey $x_1$ is persistent in $X$ if the following inequality holds
\bae\label{x1p}
\frac{dx_1}{x_1dt}\Big\vert_{X_2}&=&\frac{dx_1}{x_1dt}\Big\vert_{E_{10}}=\left\{\begin{array}{ll}\frac{dx_1}{x_1dt}\Big\vert_{E_{00}}=r_1>0 \mbox{ if } r_2<d\\\\
\frac{dx_1}{x_1dt}\Big\vert_{E_{01}} =r_1-a(1-\frac{d}{r_2})K_2>0\mbox{ if } r_2>d\end{array}\right\}
\eae This implies that $x_1$ is persistent if $r_2<d$ or  $\frac{r_1}{aK_2}>1-\frac{d}{r_2}>0$.
\item The predator $x_2$ is persistent in $X$ if the following inequality holds
\bae\label{x2p}
\frac{dx_2}{x_2dt}\Big\vert_{X_1}&=&\frac{dx_2}{x_2dt}\Big\vert_{E_{10}}=r_2-d+\frac{eaK_1}{1+haK_1}>0
\eae This implies that $x_2$ is persistent if $r_2+\frac{eaK_1}{1+haK_1}>d$.
\item The discussions of two items above imply that Model \eqref{pp-g} is permanent (i.e., both prey $x_1$ and predator $x_2$ are persistent in $X$) if either
$$\frac{r_1}{a K_2}>1-\frac{d}{r_2}>0$$
or
$$r_2+\frac{eaK_1}{1+haK_1}>d>r_2.$$

\end{enumerate}
\end{proof}
%%%%%%%%%%%%%%%%%%%%%%%%%%%%%%%%%%%%%%%%%%%%%%%%%%%%%%%%%%%%%%%%%%%%%%%%%%%%
\subsection*{Proof of Theorem \ref{th:interior}}
\begin{proof} 
The interior equilibrium $(x_1^*,x_2^*)$ of Model \eqref{pp-g} satisfies the following two equations:
\bae\label{interior}\begin{array}{lcl}
r_1 (1-\frac{x_1}{K_1})&=&\frac{ax_2}{1+h a x_1}\Leftrightarrow x_2=f(x_1)=\frac{r_1(K_1-x_1)(1+ahx_1)}{aK_1} \mbox{ subject to } 0< x_1<K_1\\\\
\frac{e a x_1}{1+ha x_1}&=&d-r_2 (1-\frac{x_2}{K_2})\Leftrightarrow  x_2=g(x_1)=\frac{K_2}{r_2}\left[\frac{e a x_1}{1+ha x_1}+(r_2-d)\right] \\
&& \mbox{ subject to } \frac{1}{a(\frac{e}{d-r_2}-h)}<x_1\mbox{ and } d<r_2+\frac{e a K_1}{1+haK_1}%\\\\

\end{array}\eae which implies that there is no interior equilibrium $(x_1^*,x_2^*)$ if $d>r_2+\frac{e a K_1}{1+haK_1}$.
The equations \eqref{interior} imply that $x_1^*$ is a root of  $fg(x_1)=g(x_1)-f(x_1)$: 
$$fg(x_1)=\frac{c_3 x^3_1+c_2x^2_1+c_1x_1+c0}{a r_2K_1(1+ahx_1)}=\frac{F(x_1)}{a r_2K_1(1+ahx_1)}$$
where
$$\begin{array}{lcl}
c_3=r_1r_2 (a h)^2&&c_2=ahr_1r_2(2-ahK_1)\\
c_1=a^2hK_1K_2\left(r_2+\frac{e}{h}-d\right)+r_1r_2(1-2ahK_1)&&c_0=a r_2K_1\left[K_2\left(1-\frac{d}{r_2}\right)-\frac{r_1}{a}\right]\end{array}$$ 

We classify interior equilibria for Model \eqref{pp-g} in the following cases:
\begin{enumerate}
\item If $c_0>0\Leftrightarrow K_2\left(1-\frac{d}{r_2}\right)-\frac{r_1}{a}>0 \Leftrightarrow \frac{r_1}{a K_2}<1-\frac{d}{r_2}$ (i.e., sufficient conditions for the locally asymptotical stability of $E_{01}$ according to Proposition \ref{p1:be}) , then we can conclude that $F(x_1)$ has either no interior equilibrium or two interior equilibria since $c_3>0$. This implies that Model \eqref{pp-g} may have either no interior equilibrium or two interior equilibria when $c_0>0$ (i.e., $E_{01}$ is locally asymptotically stable). 
\item If $c_0<0\Leftrightarrow  \frac{r_1}{a K_2}>1-\frac{d}{r_2}$ (i.e., $E_{01}$ is a saddle according to Proposition \ref{p1:be}), by using similar arguments for the case of $c_0>0$, then we can conclude that Model \eqref{pp-g} may have either one interior equilibrium or three interior equilibria when $c_0<0$.
\item Note that 
$$\begin{array}{lcl}
F'(x_1)&=&3c_3 x^2_1+2c_2x_1+c_1=3c_3 \left(x_1+\frac{c_2}{3c_3}\right)^2+c_1-\frac{c_2^2}{3c_3}\\\\
&=&3c_3 \left(x_1-\frac{ahK_1-2}{3 ah}\right)^2+a^2 r_2 K_1K_2\left(h+\frac{d}{r_2}-\frac{e}{r_2}\right)-\frac{r_1r_2(ahK_1+1)^2}{3}\end{array}.$$
Thus, we have the following two scenarios:
\begin{enumerate}

\item If $a^2 r_2 K_1K_2\left(h+\frac{d}{r_2}-\frac{e}{r_2}\right)<\frac{r_1r_2(ahK_1+1)^2}{3}$ (i.e., $c_1<\frac{c_2^2}{3c_3}$), then $F'(x_1)=0$ has two real roots 

$$x_c^1=\frac{-c_2-\sqrt{c_2^2-3c_1c_3}}{3c_3}<x_c^2=\frac{-c_2+\sqrt{c_2^2-3c_1c_3}}{3c_3}.$$ In this case, we can conclude that $F(x_1)$ has three positive roots if 
$$c_0<0, \, x_c^1>0,\,\,F(x_c^1)>0 \mbox { and } F(x_c^2)<0.$$
The equation $F(x_1)$ has one positive roots if 
$$c_0<0, \,  \mbox { and } F(x_c^2)<0.$$
The equation $F(x_1)$ has two positive roots if 
$$c_0>0, \, x_c^2>0,\,\, \mbox { and } F(x_c^2)<0.$$
\item If $a^2 r_2 K_1K_2\left(h+\frac{d}{r_2}-\frac{e}{r_2}\right)>\frac{r_1r_2(ahK_1+1)^2}{3}$  (i.e., $c_1>\frac{c_2^2}{3c_3}$), then $F'(x_1)>0$ which indicates that $F(x_1)$ is an increasing function of $x_1$. Therefore, $F(x_1)$ has one positive root if $c_0<0$ and has no positive roots if $c_0>0$.
\end{enumerate}
\item Now we discuss that the cases when Model \eqref{pp-g} has either no interior equilibrium or one interior equilibrium in terms of the function $f(x_1)$ and $g(x_1)$. First, we know that (1) The function $f(x_1)$ is a degree two polynomial with two roots $-\frac{1}{ah}$ and $K_1$. (2) $f(x_1)$ is a decreasing function whenever $x_1>\frac{K_1-1/ah}{2}$ while $g(x_1)$ is an increasing function with a unique positive root $x^*=\frac{d-r_2}{a[e-h(d-r_2)]}$ when $r_2<d<r_2+\frac{e}{h}$. Therefore, the functions $f(x_1)$ and $g(x_1)$ have no intercept for $x_1>0$ if 
$$\frac{d-r_2}{a[e-h(d-r_2)]}>K_1\Leftrightarrow d>r_2+\frac{eaK_1}{1+ahK_1}.$$
This also implies that when $r_2<d<r_2+\frac{e}{h}$ we have the following sufficient condition for $f(x_1)$ and $g(x_1)$ having a unique positive intercept:
$$\frac{d-r_2}{a[e-h(d-r_2)]}<K_1 \Leftrightarrow r_2<d<r_2+\frac{eaK_1}{1+ahK_1}.$$
If $r_2>d$, then the functions $f(x_1)$ and $g(x_1)$ have no intercept for $x_1>0$ in the following two cases:
\begin{enumerate}
\item  $\frac{K_1-1/ah}{2}<0$ and $g(0)>f(0)\Leftrightarrow K_2\left(1-\frac{d}{r_2}\right)>\frac{r_1}{a}.$
\item  $\frac{K_1-1/ah}{2}>0$ and $g(0)>f(\frac{K-1/ah}{2})\Leftrightarrow K_2\left(1-\frac{d}{r_2}\right)> \frac{r_1(1+ahK_1)^2}{4a^2hK_1}.$
\end{enumerate}
Now assume that  $K_1<\frac{1}{ah}$, then $f(x_1)$ is a decreasing function while $g(x_1)$ is an increasing function for $x_1>0$. Thus the functions $f(x_1)$ and $g(x_1)$ have an unique interior intercept for $x_1>0$ in the following two cases:
\begin{enumerate}
\item $r_2>d$ and $g(0)<f(0)\Leftrightarrow 0<K_2\left(1-\frac{d}{r_2}\right)<\frac{r_1}{a}.$
\item $r_2<d<r_2+\frac{e}{h}$ and $0<x^*=\frac{d-r_2}{a[e-h(d-r_2)]}<K_1$  where $g(x^*)=0$. 
\end{enumerate}
\end{enumerate}
The discussions above lead to the following summary regarding the interior equilibria of Model \eqref{pp-g} which is also listed in Table \ref{t2:interior}:
\begin{enumerate}
\item Model \eqref{pp-g} has no interior equilibrium if one of the following conditions is satisfied:
\begin{enumerate}
\item  $d>r_2+\frac{aeK_1}{ahK_1+1}$; or 
\item $\frac{r_1}{a K_2}<1-\frac{d}{r_2}$ and $a^2 r_2 K_1K_2\left(h+\frac{d}{r_2}-\frac{e}{r_2}\right)>\frac{r_1r_2(ahK_1+1)^2}{3}$; or 
\item $ K_1<\frac{1}{ah}$ and $\frac{r_1}{a K_2}<1-\frac{d}{r_2}$; or 
\item $ K_1>\frac{1}{ah}$ and $\frac{r_1(1+ahK_1)^2}{4a^2hK_1K_2}<1-\frac{d}{r_2}.$

\end{enumerate}
\item Model \eqref{pp-g} has an unique interior equilibrium if one of the following conditions is satisfied:
\begin{enumerate}
\item $r_2<d<r_2+\frac{aeK_1}{ahK_1+1}$; or 
\item $\frac{r_1}{a K_2}>1-\frac{d}{r_2}$ and $\left(h+\frac{d}{r_2}-\frac{e}{r_2}\right)>\frac{r_1(ahK_1+1)^2}{3a^2 K_1K_2}$;  or 
\item $\frac{r_1}{a K_2}>1-\frac{d}{r_2}$,\,$\left(h+\frac{d}{r_2}-\frac{e}{r_2}\right)<\frac{r_1(ahK_1+1)^2}{3a^2  K_1K_2}$ and $F(x^2_c)<0$.

\end{enumerate}
\item Model \eqref{pp-g} has two  interior equilibria if  $\frac{r_1}{a K_2}<1-\frac{d}{r_2}$,\,$\left(h+\frac{d}{r_2}-\frac{e}{r_2}\right)<\frac{r_1(ahK_1+1)^2}{3a^2  K_1K_2}$ and $F(x^2_c)<0$ with $x^2_c>0$.
\item Model \eqref{pp-g} has three interior equilibria if $\frac{r_1}{a K_2}>1-\frac{d}{r_2}$,\,$\left(h+\frac{d}{r_2}-\frac{e}{r_2}\right)<\frac{r_1(ahK_1+1)^2}{3a^2  K_1K_2}$, $F(x^1_c)>0$,\, and $F(x^2_c)<0$ with $x^1_c>0$.\end{enumerate} where $x_c^1=\frac{c_2}{3c_3}-\sqrt{c_1-\frac{c_2^2}{3c_3}}<x_c^2=\frac{c_2}{3c_3}+\sqrt{c_1-\frac{c_2^2}{3c_3}}$.\\
Let  $(x_1^*,x_2^*)$ be an interior equilibrium of Model \eqref{pp-g}, then its stability is determined by the eigenvalues of the following Jacobian matrix
\bae\label{Ji}
J_{(x_1^*,x_2^*)}^i=\left[\begin{array}{cc}
x_1^*\left(-\frac{r_1}{K_1}+\frac{e a x_2^*}{(1+a hx_1^*)^2}\right)&-\frac{a x_1^*}{1+h a x_1^*}\\
\frac{ea x^*_2}{(1+a hx_1^*)^2}&-\frac{r_2x_2^*}{K_2}
\end{array}\right].
\eae whose eigenvalues $\lambda_i, i=1,2$ satisfy the following equalities:
$$\lambda_1+\lambda_2=-\frac{r_1x_1^*}{K_1}-\frac{r_2x_2^*}{K_2}+\frac{ha^2x_1^*x_2^*}{(1+ahx_1^*)^2} <-\frac{r_1x_1^*}{K_1}-\frac{r_2x_2^*}{K_2}+\frac{a x^*_2}{2}$$\mbox{ and } $$\lambda_1\lambda_2=x_1^*x_2^*\left[\frac{r_1r_2}{K_1K_2}+\frac{ea^2}{(1+ahx_1^*)^3}-\frac{r_2a^2hx_2^*}{(1+ahx_1^*)^2}\right]\geq x_1^*x_2^*\left[\frac{r_1r_2}{K_1K_2}+\frac{ea^2}{(1+ah K_1)^3}-a^2hK_{2} (r_2+\frac{e }{h})\right].$$
This implies that the interior equilibrium $(x_1^*,x_2^*)$ is locally asymptotically stable if 
$$\frac{r_2}{K_2}>\frac{a}{2} \mbox{ and } \frac{r_1r_2}{K_1K_2}+\frac{ea^2}{(1+ah K_1)^3}>a^2hK_{2} (r_2+\frac{e }{h}).$$
Therefore, the statement of Theorem \ref{th:interior} holds.
\end{proof}

%%%%%%%%%%%%%%%%%%%%%%%%%%%%%%%%%%%%%%%%%%%%%%
\subsection*{Proof of Theorem \ref{th:global}}
\begin{proof} The proof of the global stability of $E_{01}$ is similar to the proof of the global stability of $E_{10}$, thus we only focus on the case of $E_{01}$.\\

According to Proposition \ref{p1:be}, we know that Model \eqref{pp-g} has only three boundary equilibria $E_{00}$, $E_{10}$ and  $E_{01}$ where $E_{00}$ is always a saddle; $E_{10}$ is locally asymptotically stable when the inequality $d>r_2+\frac{eaK_1}{1+ahK_1}$ hold while it is unstable if $d>r_2+\frac{eaK_1}{1+ahK_1}$; and   $E_{01}$ is locally asymptotically stable when the inequality $\frac{r_1}{a K_2}<1-\frac{d}{r_2}$ hold while it is unstable if $\frac{r_1}{a K_2}>1-\frac{d}{r_2}$. Proposition \ref{p1:be} also implies that if  $E_{10}$ is locally asymptotically stable then $E_{01}$ is unstable while $E_{01}$ is locally asymptotically stable then $E_{10}$ is unstable.\\

Now assume that $\frac{r_1}{a K_2}<1-\frac{d}{r_2}$ and one of the following conditions is satisfied
\begin{enumerate}
\item  $a^2 r_2 K_1K_2\left(h+\frac{d}{r_2}-\frac{e}{r_2}\right)>\frac{r_1r_2(ahK_1+1)^2}{3}$; or 
\item $ K_1<\frac{1}{ah}$; or 
\item $ K_1>\frac{1}{ah}$ and $\frac{r_1(1+ahK_1)^2}{4a^2hK_1K_2}<1-\frac{d}{r_2}.$\\
\end{enumerate}
Then according to the proof of Theorem \ref{th:interior}, Model \eqref{pp-g} has only three boundary equilibria $E_{00}$, $E_{10}$ and  $E_{01}$ where both $E_{00}$, $E_{10}$ are unstable and $E_{01}$ is locally asymptotically stable.\\

According to Lemma \ref{l1:dp},  Model \eqref{pp-g} has a compact global attractor. Thus, from an application of the Poincar\'e-Bendixson theorem (Guckenheimer and Holmes 1983 \cite{Guckenheimer1983}) we conclude that the trajectory starting at any initial condition living in the interior of $\mathbb R^2_+$ converges to one of the three boundary equilibria $E_{00}, E_{10}$ since Model \eqref{pp-g} has no interior equilibrium under the condition.  Since $E_{01}$ is the only locally asymptotically stable boundary equilibrium, therefore, every trajectory of Model \eqref{pp-g} converges to $E_{01}$ for any initial condition taken in the interior of $\mathbb R^2_+$. This implies that $E_{01}$ is global stable. Therefore, our statement holds.

\end{proof}
%%%%%%%%%%%%%%%%%%%%%%%%%%%%%%%%%%%%%%%%%%%%%%%%%%%%%%%%%%%%%
\subsection*{Proof of Theorem \ref{p2:be-ev}}
\begin{proof}
Let  $$a_v=\frac{\partial a(u,v)}{\partial v},\, a_u=\frac{\partial a(u,v)}{\partial u},\,\, a_{vv}=\frac{\partial^2 a(u,v)}{\partial v^2},\,\,  a_{uu}=\frac{\partial^2 a(u,v)}{\partial u^2},\,\,\mbox{ and }\,  a_{uv}=\frac{\partial^2 a(u,v)}{\partial v\partial u}=\frac{\partial^2 a(u,v)}{\partial u\partial v}.$$
Now we are exploring the local stability of $E_{x_10uv}=(K_1(u),0,u,v)$ and $E_{0x_2uv}=\left(0,K_2(v)\left(1-\frac{d(v)}{r_2}\right),u,v\right)$ where $u, v$ are the trait values that satisfy the following two equations:
$$\begin{array}{lcl}
\frac{r_1 x_1 K_1'(u)}{K_1^2(u)}&=&\frac{x_2\frac{\partial a(u,v)}{\partial u}}{(1+h a(u,v)x_1)^2}=\frac{x_2a_u(u,v)}{(1+h a(u,v)x_1)^2}\\\\
d'(v)&=&\frac{ex_1a_v(u,v)}{(1+h a(u,v)x_1)^2}+\frac{r_2 x_2 K_2'(v)}{K_2^2(v)}.
\end{array}
$$
Therefore, at the equilibrium $E_{x_10uv}=(K_1(u),0,u,v)$, we have 
$$K_1'(u)=0\,\, \mbox{ and }\,\, d'(v)=\frac{eK_1(u)a_v(u,v)}{(1+h a(u,v)K_1(u))^2}$$ 
with the following Jacobian matrix
\bae\label{JK1}
J_{(K_1(u), 0,u,v)}=\left[\begin{array}{cccc}
-r_1&-\frac{K_1(u)a(u,v)}{1+ha(u,v)K_1(u)}&0&0\\
0&\frac{ea(u,v)K_1(u)}{1+ha(u,v)K_1(u)}-d(v)+r_2&0&0\\
0&-\frac{\sigma_1^2 a_u(u,v)}{(1+h a(u,v)K_1(u))^2}&\frac{\sigma_1^2r_1 K_1''(u)}{K_1(u)}&0\\
\frac{\sigma_2^2 e a_v(u,v)(1-ha(u,v)K_1(u))}{(1+ha(u,v) K_1(u))^3}&\frac{\sigma_2^2 r_2  K_2'(v)}{K_2^2(v)}&b_{21}&b_{22}
\end{array}\right]
\eae with
$$b_{21}=\frac{\sigma_2^2eK_1(u)\left[a_{uv}(u,v)(1+haK_1(u))-2 K_1(u)a_u(u,v)a_v(u,v)\right]}{(1+h a(u,v)K_1(u))^3}$$
and
$$b_{22}=\frac{\sigma_2^2eK_1(u)\left[a_{vv}(u,v)(1+haK_1(u))-2 K_1(u)a^2_v(u,v)\right]}{(1+h a(u,v)K_1(u))^3}-\sigma_2^2d''(v).$$

The eigenvalues of \eqref{JK1} are 
$$-r_1,\,\,\frac{ea(u,v)K_1(u)}{1+ha(u,v)K_1(u)}-d(v)+r_2,\,\,\frac{\sigma_1^2r_1 K_1''(u)}{K_1(u)},\,\,\frac{\sigma_2^2eK_1(u)\left[a_{vv}(u,v)(1+haK_1(u))-2 K_1(u)a^2_v(u,v)\right]}{(1+h a(u,v)K_1(u))^3}-\sigma_2^2d''(v).$$
This implies that if the equilibrium $E_{x_10uv}=(K_1(u),0,u,v)$ exists, then it is locally asymptotically stable if 
$$0<\frac{ea(u,v)K_1(u)}{1+ha(u,v)K_1(u)}<d(v)-r_2,\,\,K_1''(u)<0,\,\, \mbox{ and } \frac{eK_1(u)\left[a_{vv}(u,v)(1+haK_1(u))-2 K_1(u)a^2_v(u,v)\right]}{(1+h a(u,v)K_1(u))^3}<d''(v).$$

Similarly, if the equalities $a_u(u,v)=0,\,d(v)<r_2$ and $d'(v)=\frac{r_2K'_2(v)}{K_2(v)}$ hold, then the equilibrium $E_{0x_2uv}=\left(0,K_2(v)\left(1-\frac{d(v)}{r_2}\right),u,v\right)$ exists  with the following Jacobian matrix
\bae\label{JK2}
J_{(0,K_2(v), u,v)}=\left[\begin{array}{cccc}
r_1-a(u,v) K_2(v)\left(1-\frac{d(v)}{r_2}\right)&0&0&0\\
ea(u,v)K_2(v)\left(1-\frac{d(v)}{r_2}\right)&d(v)-r_2&0&0\\
\frac{\sigma_1^2 r_1 K'_1(u)}{K^2_1(u)}&0&-\sigma_1^2K_2(v)\left(1-\frac{d(v)}{r_2}\right)a_{uu}(u,v)&b_{12}\\
\sigma_2^2 ea_v(u,v)&\frac{\sigma_2^2r_2\left(1-\frac{d(v)}{r_2}\right)K_2'(v)}{K_2(v)}&0&b_{22}\end{array}\right]
\eae where
$$b_{12}=-\sigma_1^2K_2(v)\left(1-\frac{d(v)}{r_2}\right)a_{uv}(u,v)$$ and 
$$b_{22}=\frac{\sigma_2^2r_2\left(1-\frac{d(v)}{r_2}\right)\left[K_2(v) K''_2(v)-2(K'_2(v))^2\right]}{K_2(v)^2}-\sigma_2^2d''(v).$$

The eigenvalues of \eqref{JK2} are 
$$r_1-a(u,v) K_2(v)\left(1-\frac{d(v)}{r_2}\right),\,\,-d(v)+r_2,\,\,-\sigma_1^2K_2(v)\left(1-\frac{d(v)}{r_2}\right)a_{uu}(u,v),\,\,$$ and $$\frac{\sigma_2^2r_2\left(1-\frac{d(v)}{r_2}\right)\left[K_2(v) K''_2(v)-2(K'_2(v))^2\right]}{K_2(v)^2}-\sigma_2^2d''(v).$$
This implies that the facultative-social-parasite-only equilibrium $E_{0x_2uv}=(0,K_2(v),u,v)$ exists  and it is locally asymptotically stable if the following inequalities hold
$$\frac{r_1}{1-\frac{d(v)}{r_2}}<a(u,v) K_2(v),\,\,r_2<d(v),\,\,a_{uu}(u,v)>0,\,\, \mbox{ and } \frac{r_2\left(1-\frac{d(v)}{r_2}\right)\left[K_2(v) K''_2(v)-2(K'_2(v))^2\right]}{K_2(v)^2}<d''(v).$$
\end{proof}
%%%%%%%%%%%%%%%%%%%%%%%%%%%%%%%%%%%%%%%%%%%%%%
\subsection*{Proof of Theorem \ref{th:be}}
\begin{proof}
First, we can check that the boundary equilibrium of \eqref{E-pp-g} can occur only if $u^*_1=u^*_2=0$. According to Proposition \ref{p1:be}, the boundary equilibria are: 
$$E_{0000}=(0,0,0,0),\,\, E_{x_1000}=(K_{01},0,0,0),\,\mbox{ and } E_{0x_200}=\left(0,K_{02}\left(1-\frac{d}{r_2}\right),0,0\right) \mbox{ provided that } r_2>d.$$
When $u_1^*=u^*_2=0$, the Jacobian matrix evaluated at an equilibrium $(x^*_1, x^*_2,0,0)$ can be represented as follows:
\bae\label{J0}
J_{(x^*_1, x^*_2,0,0)}=\left[\begin{array}{cccc}
a_{11}&a_{12}&0&0\\
a_{21}&a_{22}&0&0\\
0&0&b_{11}&b_{12}\\
0&0&b_{21}&b_{22}
\end{array}\right]
\eae whose eigenvalues are determined by the eigenvalues of matirices $A$ and $B$
\be\label{A}
A=\left[\begin{array}{cc}
a_{11}&a_{12}\\
a_{21}&a_{22}\end{array}\right]\\
=\left[\begin{array}{cc}
x_1^*\left(-\frac{r_1}{K_{01}}+\frac{e a_0 x_2^*}{(1+a_0 hx_1^*)^2}\right) + \left(r_1(1-\frac{x^*_1}{K_{01}})-\frac{a_0 x^*_2}{1+a_0 hx_1^*}\right)&-\frac{a_0 x_1^*}{1+h a_0 x_1^*}\\
\frac{ea_0 x^*_2}{(1+a_0 hx_1^*)^2}&-\frac{r_2x_2^*}{K_{02}}+ r_2(1-\frac{x^*_2}{K_{02}})+\frac{ea_0 x^*_1}{1+a_0 hx_1^*}-d
\end{array}\right]
\ee and 

\be\label{B}
B=\left[\begin{array}{cc}
b_{11}&b_{12}\\
b_{21}&b_{22}\end{array}\right]\\
=\left[\begin{array}{cc}
\sigma_1^2\left(-\frac{r_1x_1^*}{\sigma_{K_1}^2K_{01}}+\frac{a_0x_2^*}{\sigma_a^2(1+a_0 hx_1^*)^2}\right)&-\frac{\sigma_1^2 a_0 x_2^*}{\sigma_a^2(1+a_0 hx_1^*)^2}\\
\frac{\sigma_2^2 e a_0 x_1^*}{\sigma_a^2(1+a_0 hx_1^*)^2}&-\sigma_2^2\left(\frac{r_2x_2^*}{\sigma_{K_2}^2K_{02}}+\frac{a_0ex_1^*}{\sigma_a^2(1+a_0 hx_1^*)^2}\right) \end{array}\right].
\ee Therefore, we can conclude the following cases:
\begin{enumerate}
\item $E_{0000}$: In this case, the eigenvalues of $A$ are $r_1$ and $r_2-d$ while the eigenvalues of $B$ are zeros. Thus,  $E_{0000}$ is always a saddle.
\item $E_{x_1000}=(K_{01},0,0,0)$: In this case, the eigenvalues of $A$ are $-r_1$ and $r_2-d+\frac{ea_0K_{01}}{1+ahK_{01}}$ while the eigenvalues of $B$ are $-\frac{r_1\sigma_1^2}{\sigma_{K_1}^2}$ and $-\frac{\sigma_2^2ea_0K_{01}}{\sigma^2_a(1+a_0hK_{01})^2}$. Thus,  $E_{x_1000}$ is locally asymptotically stable if $r_2+\frac{ea_0K_{01}}{1+ahK_{01}}<d$ (i.e., $E_{10}$ is locally asymptotically stable for the ecological model \eqref{pp-g}) and $E_{x_1000}$ is a saddle if $r_2+\frac{ea_0K_{01}}{1+ahK_{01}}>d$.
\item $E_{0x_200}=(0,K_{02}\left(1-\frac{d}{r_2}\right),0,0)$ exists if $r_2>d$:  In this case, the eigenvalues of $A$ are $r_1+\frac{a_0d K_{02}}{r_2}-a_0K_{02}$ and $d-r_2$ while the eigenvalues of $B$ are $-\frac{\sigma_1^2K_{02}a_0(d-r_2)}{\sigma_{a}^2r_2}$ and $\frac{\sigma_2^2(d-r_2)}{\sigma^2_{K_2}}$. Thus,  $E_{0x_200}$ is always a saddle.
\item When $u^*_1=u^*_2=0$, the existence conditions regarding interior equilibria $(x_1^*,x_2^*,0,0)$ of the evolutionary model \eqref{E-pp-g} are the same as interior equilibria $(x_1^*,x_2^*)$ of the ecological model \eqref{pp-g} by letting $K_1=K_{01},\, K_2=K_{02} \mbox{ and } a=a_0$. Thus, according to Theorem \ref{th:interior}, Model \eqref{E-pp-g} can have one, two, or three interior equilibria $(x_1^*,x_2^*,0,0)$ depending on the values of $r_i, d, h, e, a_0$ and $K_{0i}$ where $i=1,2$. Sufficient condition of an interior equilibrium $(x_1^*,x_2^*,0,0)$ being locally asymptotically stable is that both matrices $A$ and $B$ have negative eigenvalues.

According to Theorem \ref{th:interior},  the eigenvalues of $A$ are negative if inequalities \eqref{stable} hold, i.e.,
$$\frac{r_2}{K_{02}}>\frac{a_0}{2} \mbox{ and } \frac{r_1r_2}{K_{01}K_{02}}+\frac{ea_0^2}{(1+a_0h K_{01})^3}>a_0^2hK_{02} (r_2+\frac{e }{h}).$$
Let $\lambda_1$ and $\lambda_2$ be eigenvalues of $B$. Then we have as follows: 
$$\begin{array}{lcl}
\lambda_1+\lambda_2&=&-\frac{r_1\sigma_1^2x_1^*}{K_{01}\sigma^2_{K_1}}-\frac{r_2\sigma_2^2x_2^*}{K_{02}\sigma^2_{K_2}}+\frac{a_0(\sigma_1^2x_1^*-e\sigma_2^2x_2^*)}{(1+a_0hx_1^*)^2}\\
&<&-\left(\frac{r_1}{K_{01}\sigma^2_{K_1}}-a_0\right)\sigma_1^2x^*_1-\frac{r_2\sigma_2^2x_2^*}{K_{02}\sigma^2_{K_2}}-\frac{a_0e\sigma_2^2x_2^*}{(1+a_0hx_1^*)^2}
\end{array} $$ and
$$\begin{array}{lcl}\lambda_1\lambda_2&=&\sigma_1^2\sigma_2^2\left[\frac{r_1r_2x^*_1x^*_2}{K_{01}K_{02}}+
\frac{a_0r_1e(x_1^*)^2}{\sigma_{K_1}^2K_{01}\sigma_a^2(1+ahx_1^*)^2}  -\frac{a_0r_2(x^*_2)^2}{\sigma_{K_2}^2K_{02}\sigma_a^2(1+ahx_1^*)^2}\right]\\
&=&\sigma_1^2\sigma_2^2\left[\frac{r_1r_2x^*_1x^*_2}{K_{01}K_{02}}+\frac{a_0\left(\frac{r_1e(x_1^*)^2}{\sigma_{K_1}^2K_{01}}-\frac{r_2(x^*_2)^2}{\sigma_{K_2}^2K_{02}}\right)}{\sigma_a^2(1+ahx_1^*)^2}\right]\\
&>&\sigma_1^2\sigma_2^2\left[\frac{r_1r_2x^*_1x^*_2}{K_{01}K_{02}}+\frac{a_0r_1e(x_1^*)^2}{\sigma_{K_1}^2K_{01}\sigma_a^2(1+ahx_1^*)^2}  -\frac{a_0r_2(x^*_2)^2}{\sigma_{K_2}^2K_{02}\sigma_a^2}\right]\\
&=&\sigma_1^2\sigma_2^2\left[\frac{r_1r_2x^*_1x^*_2}{K_{01}K_{02}}  -\frac{a_0r_2(x^*_2)^2}{\sigma_{K_2}^2K_{02}\sigma_a^2}+\frac{a_0r_1e(x_1^*)^2}{\sigma_{K_1}^2K_{01}\sigma_a^2(1+ahx_1^*)^2}\right]\\
&=&\sigma_1^2\sigma_2^2\left[\frac{r_2x^*_2}{K_{02}}\left(\frac{r_1x^*_1}{K_{01}}  -\frac{a_0x^*_2}{\sigma_{K_2}^2\sigma_a^2}\right)+\frac{a_0r_1e(x_1^*)^2}{\sigma_{K_1}^2K_{01}\sigma_a^2(1+ahx_1^*)^2}\right]\\
\end{array}$$
Therefore, we can conclude that if $\frac{r_1}{K_{01}\sigma^2_{K_1}}>a_0$ and one of the following two inequalities hold
\bae\label{evo1}
\begin{array}{lcl}\frac{r_1x^*_1}{K_{01}}  &>&\frac{a_0x^*_2}{\sigma_{K_2}^2\sigma_a^2}\Leftrightarrow \frac{r_1\sigma_{K_2}^2\sigma_a^2}{a_0K_{01}}  >\frac{x^*_2}{x^*_1} \mbox{ or }\\
$$\frac{r_1e\sigma_{K_2}^2K_{02}}{r_2\sigma_{K_1}^2K_{01}}&>&\frac{(x^*_2)^2}{(x_1^*)^2}\Leftrightarrow \frac{\sigma_{K_2}^2}{\sigma_{K_1}^2}\frac{r_1eK_{02}}{r_2K_{01}}>\frac{(x^*_2)^2}{(x_1^*)^2}\end{array},\eae then we have $\lambda_1+\lambda_2<0$ and  $\lambda_1\lambda_2>0$, which implies that $B$ has two negative eigenvalues. Thus, we can conclude that if the interior equilibrium $(x_1^*,x_2^*,0,0)$ exists, then it is locally asymptotical stable if one of the two inequalities in \eqref{evo1} holds and the following inequalities hold:
$$\min\Big\{\frac{2r_2}{K_{02}},\frac{r_1}{\sigma^2_{K_1}K_{01} }\Big\}>a_0 \mbox{ and } \frac{r_1r_2}{K_{01}K_{02}}+\frac{ea_0^2}{(1+a_0h K_{01})^3}>a_0^2hK_{02} (r_2+\frac{e }{h}).$$
\end{enumerate}

Therefore, the statement of Theorem \ref{th:be} holds.

\end{proof}

%%%%%%%%%%%%%%%%%%%%%%%%%%%%%%%%%%%%%%%%%%%%%%%%%%%%%%%%%%%%%
\subsection*{Proof of Theorem \ref{th:ess-b}}
\begin{proof}
According to the definition of ESS, $E_{x_1000}$ is an ESS if it is an ESE and it satisfies the ESS maximum principle \eqref{ESS-max}. According to Theorem \ref{th:be}, $E_{x_1000}$ is an ESE if the inequality $r_2+\frac{ea_0K_{01}}{1+ahK_{01}}<d$ holds. Notice that
$$\begin{array}{lcl}\max_{v_1\in \mathbb R}\{G_1(v_1, 0,0,K_{01},0)\}&=&\max_{v_1\in \mathbb R}\Big\{r_1\left(1-\frac{K_{01}}{K_1(v_1)}\right)\Big\}\\
&=&r_1\left(1-\frac{K_{01}}{K_1(v_1)}\right)\Big\vert_{v_1=0}=G_1(0, 0,0,K_{01},0)=0\end{array}$$ and 
$$\begin{array}{lcl}\max_{v_2\in \mathbb R}\{G_2(v_2, 0,0,K_{01},0)\}&=&\max_{v_2\in \mathbb R}\Big\{\frac{e\, a(0,v_2) K_{01}}{1+h\,a(0,v_2)K_{01}}-d+r_2 \Big\}\\
&=&\frac{e\, a(0,v_2) K_{01}}{1+h\,a(0,v_2)K_{01}}-d+r_2\Big\vert_{v_2=0}=G_2(0, 0,0,K_{01},0)<0\end{array}.$$
Therefore, $E_{x_1000}$ is an ESS whenever the inequality $r_2+\frac{ea_0K_{01}}{1+ahK_{01}}<d$ holds.\\

Similarly, the interior equilibrium $(x^*_1, x^*_2, 0,0)$ is an ESS if it is an ESE and it satisfies the ESS maximum principle \eqref{ESS-max}. According to Theorem \ref{th:be}, $(x^*_1, x^*_2, 0,0)$ is an ESE if the inequalities \eqref{stable-evo} hold. Notice that
$K_i(v), a(0,v), a(v,0), i=1,2$ are decreasing functions with respect to $\vert v\vert$, and the following equations
$$\begin{array}{lcl}
G_2(v_2, 0,0,x^*_1, x^*_2)&=&\frac{e\, a(0,v_2) x^*_1}{1+h\,a(0,v_2)x^*_1}-d+r_2 \left(1-\frac{x^*_2}{K_2(v_2)}\right)\\
\frac{\partial G_2(v_2, 0,0,x^*_1, x^*_2)}{\partial v_2}&=&-v_2\sigma_2^2\left[\frac{r_2 x_2^* }{\sigma_{K_2}^2K_2(v_2)}+\frac{ex_1^*a(0,v_2)}{\sigma_a^2(1+h a(0,v_2)x_1^*)^2}\right].
\end{array}$$Therefore, we can conclude that 
$$\begin{array}{lcl}\max_{v_2\in \mathbb R}\{G_2(v_2, 0,0,x^*_1, x^*_2)\}&=&\max_{v_2\in \mathbb R}\Big\{\frac{e\, a(0,v_2) x^*_1}{1+h\,a(0,v_2)x^*_1}-d+r_2 \left(1-\frac{x^*_2}{K_2(v_2)}\right)\Big\}\\
&=&\frac{e\, a_0 x_1^*}{1+ha_0x_1^*}-d+r_2 \left(1-\frac{x^*_2}{K_{02}}\right)=G_2(0, 0,0,x^*_1, x^*_2)=0\end{array}.$$
On the other hand, we have
$$\begin{array}{lcl}
G_1(v_1, 0,0,x^*_1, x^*_2)&=&r_1\left(1-\frac{x_1^*}{K_1(v_1)}\right)-\frac{a(v_1,0) x_2^*}{1+ h\,a(v_1,0)x_1^*}\\
\frac{\partial G_1(v_1, 0,0,x^*_1, x^*_2)}{\partial v_1}&=&\sigma_1^2v_1\left[-\frac{r_1 x_1^*}{\sigma_{K_1}^2K_1(v_1)}+\frac{x_2^*a(v_1,0)}{\sigma_a^2(1+h a(v_1,0)x_1^*)^2}\right].
\end{array}$$
Since the inequalities \eqref{stable-evo}  hold, then we have
$$\begin{array}{lcl}
-\frac{r_1 x_1^*}{\sigma_{K_1}^2K_1(v_1)}+\frac{x_2^*a(v_1,0)}{\sigma_a^2(1+h a(v_1,0)x_1^*)^2}&\leq&-\frac{r_1 x_1^*}{\sigma_{K_1}^2K_1(v_1)}+\frac{x_1^*\max\Big\{\frac{\sigma_{K_2}}{\sigma_{K_1}}\sqrt{\frac{r_1eK_{02}}{r_2K_{01}}}, \frac{r_1\sigma_{K_2}^2\sigma_a^2}{a_0K_{01}}\big\}a(v_1,0)}{\sigma_a^2(1+h a(v_1,0)x_1^*)^2}\\
&\leq&x_1^*\left[-\frac{r_1 }{\sigma_{K_1}^2K_{01}}+\frac{\max\Big\{\frac{\sigma_{K_2}}{\sigma_{K_1}}\sqrt{\frac{r_1eK_{02}}{r_2K_{01}}}, \frac{r_1\sigma_{K_2}^2\sigma_a^2}{a_0K_{01}}\big\}a_0}{\sigma_a^2}\right].
\end{array}$$This implies that if $\frac{r_1 }{\sigma_{K_1}^2K_{01}}>\frac{\max\Big\{\frac{\sigma_{K_2}}{\sigma_{K_1}}\sqrt{\frac{r_1eK_{02}}{r_2K_{01}}}, \frac{r_1\sigma_{K_2}^2\sigma_a^2}{a_0K_{01}}\big\}a_0}{\sigma_a^2}$, then we have
$$\frac{\partial G_1(v_1, 0,0,x^*_1, x^*_2)}{\partial v_1}<0 \mbox{ if } v_1>0; \frac{\partial G_1(v_1, 0,0,x^*_1, x^*_2)}{\partial v_1}>0 \mbox{ if } v_1<0.$$ Therefore, we can conclude that
$$\begin{array}{lcl}\max_{v_1\in \mathbb R}\{G_1(v_1, 0,0,x^*_1, x^*_2)\}&=&\max_{v_1\in \mathbb R}\Big\{r_1\left(1-\frac{x_1^*}{K_1(v_1)}\right)-\frac{a(v_1,0) x_2^*}{1+ h\,a(v_1,0)x_1^*}\Big\}\\
&=&G_1(0, 0,0,x^*_1, x^*_2)=0\end{array}$$ if the following inequality holds
$$\frac{r_1 }{\sigma_{K_1}^2K_{01}}>\frac{\max\Big\{\frac{\sigma_{K_2}}{\sigma_{K_1}}\sqrt{\frac{r_1eK_{02}}{r_2K_{01}}}, \frac{r_1\sigma_{K_2}^2\sigma_a^2}{a_0K_{01}}\big\}a_0}{\sigma_a^2}.$$

\end{proof}

%%%%%%%%%%%%%%%%%%%%%%%%%%%%%%%%%%%%%%%%%%%%%%%%%%%%%%%%%%%%%
\subsection*{Proof of Theorem \ref{th:unique}}
\begin{proof}
According to Theorem \ref{th:ess-b}, we know that the boundary equilibrium $E_{x_1000}$ is an ESS of the co-evolutionary host-parasite model \eqref{E-pp-g} whenever the inequality $r_2+\frac{ea_0K_{01}}{1+a_0hK_{01}}<d$ holds. Under this condition, Model \eqref{E-pp-g} has only two boundary equilibria $E_{0000}$ and $E_{x_1000}$ where $E_{0000}$ is always a saddle.

According to Lemma \ref{l1:dp}, we know that $\limsup_{t\rightarrow\infty} x_1(t)\leq K_{01}$. Since $a(u_1,u_2)\leq a_0$ and the function $\frac{e a(u_1, u_2) x_1}{1+ha(u_1,u_2)x_1}$ increases with respect to $ax_1$, thus for time large enough, we have
$$\frac{dx_2}{x_2dt}=\frac{e a(u_1, u_2) x_1}{1+ha(u_1,u_2)x_1}-d+r_2 (1-\frac{x_2}{K_2(u_2)})\leq \frac{ea_0K_{01}}{1+a_0hK_{01}}-d+r_2<0.$$
This implies that Model \eqref{E-pp-g} has no interior equilibrium. Therefore, when the inequality $r_2+\frac{ea_0K_{01}}{1+ahK_{01}}<d$ holds, Model \eqref{E-pp-g} has only a unique ESE $E_{x_1000}$ which is also an ESS.

If $K_{02} a_0<\frac{r_1}{1-\frac{d}{r_2}}$ and $K_{01}a_0<\frac{1}{h}$ holds, then for any given trait $u=(u_1,u_2)$, we have $K_2 a<\frac{r_1}{1-\frac{d}{r_2}}$. According to Proposition \ref{p1:be} and Theorem \ref{th:interior}, the ecological model \eqref{pp-g} has a unique interior equilibrium $(x_1^*, x_2^*)$ and three boundary equilibria $E_{00}, E_{10}, E_{01}$. Thus the evolutionary model \eqref{E-pp-g} has at least one interior equilibrium $(x_1^*, x_2^*, 0,0)$ and three boundary equilibria $E_{0000}, E_{x_1000}, E_{0x_200}$ which are all saddles. According to Theorem \ref{th:ess-b}, the interior equilibrium $(x_1^*, x_2^*, 0,0)$ is an ESS if \eqref{ess-unique} holds.\\

Now we should show that the evolutionary model \eqref{E-pp-g} has no interior equilibrium $(x_1^*, x_2^*, u^*_1,u^*_2)$ with $u^*_1u^*_2>0$ when \eqref{ess-unique} holds. Assume that this is not true, then let $(x_1^*, x_2^*, u^*_1,u^*_2)$ be the interior equilibrium, then it satisfies the equations\eqref{eq-u1} and \eqref{eq-u2} which gives the following equality:
\bae\label{x1-x2}
\frac{r_1 x^*_1 u^*_1\sigma_{K_2}^2K_2(u^*_2)}{r_2 x^*_2 u^*_2\sigma_{K_1}^2K_1(u^*_1)}=\frac{x_2}{ex_1}\Leftrightarrow x_2^*=\sqrt{\frac{e r_1\sigma_{K_2}^2 u^*_1K_2(u^*_2)}{r_2 \sigma_{K_1}^2u^*_2K_1(u^*_1)}}x^*_1=\frac{\sigma_{K_2}}{\sigma_{K_1}}\sqrt{\frac{e r_1 u^*_1K_{02}e^{-\frac{(u^*_2)^2}{2\sigma_{K_2}^2}+\frac{(u^*_1)^2}{2\sigma_{K_1}^2}}}{r_2 u^*_2K_{01}}}x^*_1.\eae
In addition, the equations\eqref{eq-u1} and \eqref{eq-u2} requires $u^*_1, u^*_2$ and $u^*_1-u^*_2$ having the same sign. This implies that $\vert u^*_1\vert> \vert u^*_2\vert. $ Thus, according to \eqref{x1-x2}, if $\sigma_{K_2}>\sigma_{K_1}$, then we have
$$x_2^*=\frac{\sigma_{K_2}}{\sigma_{K_1}}\sqrt{\frac{e r_1 u^*_1K_{02}e^{-\frac{(u^*_2)^2}{2\sigma_{K_2}^2}+\frac{(u^*_1)^2}{2\sigma_{K_1}^2}}}{r_2 u^*_2K_{01}}}x^*_1>\frac{\sigma_{K_2}}{\sigma_{K_1}}\sqrt{\frac{e r_1K_{02}e^{-\frac{(u^*_2)^2}{2\sigma_{K_2}^2}+\frac{(u^*_2)^2}{2\sigma_{K_1}^2}}}{r_2 K_{01}}}x^*_1>\frac{\sigma_{K_2}}{\sigma_{K_1}}\sqrt{\frac{e r_1K_{02}}{r_2 K_{01}}}x^*_1.$$This is a contradiction to the assumption that $\frac{\sigma_{K_2}}{\sigma_{K_1}}\sqrt{\frac{r_1eK_{02}}{r_2K_{01}}}>\frac{x^*_2}{x_1^*}$ in \eqref{ess-unique}. Thus, Model \eqref{E-pp-g} has no interior equilibrium $(x_1^*, x_2^*, u^*_1,u^*_2)$ with $u^*_1u^*_2>0$ when \eqref{ess-unique} holds. This indicates that the interior equilibrium $(x_1^*, x_2^*, 0,0)$ is the unique ESS if \eqref{ess-unique} holds. 

\end{proof}

%%%%%%%%%%%%%%%%%%%%%%%%%%%%%%%%%%%%%%%%%%%%%%%%%%%%%%%%%%%%%
\subsection*{Proof of Proposition \ref{p3:ess} }
\begin{proof}
 If $u^*_1u^*_2\neq 0$,

then from equations \eqref{eq-u1} and \eqref{eq-u2}, we must have $x^*_1x^*_2>0$. Let $g(u_1,u_2)=\sqrt{\frac{e r_1\sigma_{K_2}^2 u_1K_2(u_2)}{r_2 \sigma_{K_1}^2u_2K_1(u_1)}}$, then according to the equation \eqref{x1-x2}, i.e.
$$
\frac{r_1 x_1 u_1\sigma_{K_2}^2K_2(u_2)}{r_2 x_2 u_2\sigma_{K_1}^2K_1(u_1)}=\frac{x_2}{ex_1}\Leftrightarrow x_2=\sqrt{\frac{e r_1\sigma_{K_2}^2 u_1K_2(u_2)}{r_2 \sigma_{K_1}^2u_2K_1(u_1)}}x_1=g(u_1,u_2)x_1$$
which gives the following equality when it combines with the equation \eqref{eq-x1}:
$$r_1 \left(1-\frac{x_1}{K_1(u_1)}\right)=\frac{a(u_1,u_2) \sqrt{\frac{e r_1 u_1\sigma_{K_2}^2K_2(u_2)}{r_2 u_2\sigma_{K_1}^2K_1(u_1)}}x_1}{1+h a(u_1,u_2)x_1}=\frac{a(u_1,u_2)g(u_1,u_2)x_1}{1+h a(u_1,u_2)x_1}.$$Thus we have a unique solution $x_1^*$ that is a function of traits $u_i, i=1,2$:
\bae\label{solution-x1}
\begin{array}{lcl}
x^*_1&=&f_1(u_1,u_2)=\frac{K_1}{2}-\frac{r_1+a K_1 g}{2r_1a h}+\frac{\sqrt{\left(r_1+a K_1g\right)^2+r_1a h K_1 \left(r_1 a h K_1+2 r_1-2  K_1g\right)}}{2r_1a h}
\end{array}\eae where $g=g(u_1,u_2)=\sqrt{\frac{e r_1\sigma_{K_2}^2u_1K_2(u_2)}{r_2 \sigma_{K_1}^2u_2K_1(u_1)}}$. Similarly, Equation \eqref{x1-x2} combines with \eqref{eq-x2} gives the following equality:
$$r_2 \left(1-\frac{x_1 g(u_1,u_2)}{K_2(u_2)}\right)=-\frac{e a(u_1, u_2) x_1}{1+ha(u_1,u_2)x_1}$$ which also gives a unique solution
$x_1^*$ that is a function of traits $u_i, i=1,2$:
\bae\label{solution2-x1}
\begin{array}{lcl}
x^*_1&=&f_2(u_1,u_2)=\frac{K_2}{2g}+\frac{a e K_2}{2r_2a hg}-\frac{1}{2a h}+\frac{\sqrt{\left(r_2g-a e K_2\right)^2+r_2a h K_2 \left(r_2 a h K_2+2 r_2g+2 e a K_2\right)}}{2r_1a h}.
\end{array}\eae Therefore, we are able to solve for $u^*_1$ and $u^*_2$ by letting
$$f_1(u^*_1,u^*_2)=f_2(u^*_1,u^*_2) \mbox{ and } \frac{r_1 u_1^*}{\sigma_{K_1}^2K_1(u_1^*)}=\frac{(u_1^*-u_2^*)g(u_1^*,u^*_2)a(u^*_1,u^*_2)}{\sigma_a^2\left(1+h a(u^*_1,u^*_2)f_1(u^*_1,u^*_2)\right)^2}.$$
Based on the discussion above, we can conclude that the necessary condition for $(u^*_1,u^*_2, x^*_1, x^*_2)$ being an interior equilibrium of Model \eqref{E-pp-g} is that the following equalities hold:
$$x_1^*=f_1(u^*_1,u^*_2)=\frac{x^*_2}{g(u^*_1,u^*_2)}=f_2(u^*_1,u^*_2)\mbox{ and } \frac{r_1 u_1^*}{\sigma_{K_1}^2K_1(u_1^*)}=\frac{(u_1^*-u_2^*)g(u_1^*,u^*_2)a(u^*_1,u^*_2)}{\sigma_a^2\left(1+h a(u^*_1,u^*_2)f_1(u^*_1,u^*_2)\right)^2}.$$

\end{proof}

%%%%%%%%%%%%%%%%%%%%%%%%%%%%%%%%%%%%%%%%%%%%%%%%%%%%%%%%%%%%%
%%%%%%%%%%%%%%%%%%%%%%%%%%%%%%%%%%%%%%%%%%%%%%%%%%%%%%%%%%%%%
%%%%%%%%%%%%%%%%%%%%%%%%%%%%%%%%%%%%%%%%%%%%%%%%%%%%%%%%%%%%%

%\newpage
\section*{Acknowledgements}

This research is partially supported by Simons Collaboration Grants for Mathematicians (208902) to Y. K., NSF DMS (1313312) to Y. K. and J. H. F., and Arizona State University College of Letters and Sciences.

%\section*{References}
\bibliographystyle{elsarticle-harv}

%%%%%%%%%%%%%%%%%%%%%%%%%%%%%%%%%%%%%%%%%%%%%%%%%%%%%%%
%%%%%%%%%%%%%%%%%%%%%%%%%%%%%%%%%%%%%%%%%%%%%%%%%%%%%%%

%%%%%%%%%%%%%%%%%%%%%%%%%%%%%%%%%%%%%%%%%%%%%%%%%%%%%%%%%
%%%%%%%%%%%%%%%%%%%%%%%%%%%%%%%%%%%%%%%%%%%%%%%%%%%%%%%%%

\end{document}